\documentclass[10pt, a4paper]{article}

\newif\ifexternalize
\externalizefalse

\usepackage[utf8]{inputenc}
\usepackage[T1]{fontenc}
\usepackage{lmodern}

\usepackage{geometry}
\usepackage{setspace}
\usepackage{parskip}
\usepackage{fancyhdr}
\usepackage{titlesec}
\usepackage{abstract}
\usepackage{graphicx}
\usepackage{url}

\usepackage{amsmath}
\usepackage{amssymb}
\usepackage{bm}
\usepackage{mathtools}
\usepackage{amsthm}

\usepackage{xcolor}
\usepackage{hyperref}
\hypersetup{
    colorlinks=true,
    linkcolor=rose,
    urlcolor=rose,
    citecolor=rose,
}

\usepackage{cleveref}

\titleformat{\section}{\large\bfseries}{\thesection}{1em}{}
\titleformat{\subsection}{\normalsize\bfseries}{\thesubsection}{1em}{}

\theoremstyle{plain}
\newtheorem{theorem}{Theorem}[section]
\crefname{theorem}{theorem}{theorems}
\Crefname{theorem}{Theorem}{Theorems}

\theoremstyle{definition}

\crefname{definition}{definition}{definitions}
\Crefname{definition}{Definition}{Definitions}

\theoremstyle{remark}
\newtheorem{remark}[theorem]{Remark}
\crefname{remark}{remark}{remarks}
\Crefname{remark}{Remark}{Remarks}

\renewcommand{\vec}[1]{\bm{#1}}  
\renewcommand{\d}[1][]{\,\mathrm{d}#1} 

\newcommand{\npe}{_{n+1}}

\newcommand{\npeh}{_{n+1/2}}

\newcommand{\n}{_{n}}
\newcommand{\h}{^\mathrm{h}}
\newcommand{\e}{^e}

\newcommand{\nel}{n_{\mathrm{e}}}
\newcommand{\nnode}{n_{\mathrm{n}}}
\newcommand{\transp}{^\top}
\renewcommand{\epsilon}{\varepsilon}

\newcommand{\domainspace}{\Omega}
\newcommand{\nodal}[1]{\hat{\vec{#1}}}
\newcommand{\state}{\vec{x}}
\newcommand{\operator}[1]{\mathcal{#1}}
\newcommand{\R}{\mathbb{R}}
\newcommand{\diag}{\mathrm{diag}}

\newcommand{\inv}{^{-1}}
\newcommand{\rotmat}{\vec{R}}
\renewcommand{\centerline}{\vec{\varphi}}
\newcommand{\vphi}{\vec{v}_\varphi}
\newcommand{\vphis}{\vec{v}_{\varphi,s}}
\newcommand{\dotvphi}{\dot{\vec{v}}_\varphi}
\newcommand{\nodevphi}{\nodal{v}_\varphi}
\newcommand{\nodedotvphi}{\dot{\nodal{v}}_\varphi}
\newcommand{\nodecenterline}{\nodal{\varphi}}
\newcommand{\finaltime}{t_{\mathrm{end}}}

\newcommand{\sdot}[2][7mu]{\dot{#2\mkern#1}\mkern-#1}

\newcommand{\inner}[2]{\left(#1, #2 \right)_\domainspace }
\newcommand{\innerBoundary}[2]{\left[#1, #2 \right]_{\partial \domainspace} }

\usepackage{blindtext}

\usepackage{dashbox}%

\newcommand\colordashedph[1][H]{\setlength{\fboxsep}{0.3pt}\setlength{\dashlength}{2.2pt}\setlength{\dashdash}{1.1pt} \dbox{\colorbox{lightgray!40!Lavender!40!}{\phantom{#1}}}}

\newcommand\colordashedphsmall{%
    \setlength{\fboxsep}{0.2pt}%
    \setlength{\dashlength}{1.6pt}%
    \setlength{\dashdash}{0.8pt}%
    \dbox{%
        \colorbox{lightgray!40!Lavender!40!}{%
            \phantom{\rule{4pt}{4pt}}%
        }%
    }%
}

\renewcommand{\square}{\colordashedph}
\newcommand{\smallsquare}{\colordashedphsmall}

\newcommand{\zeros}{\color{lightgray!40!Lavender!40!} \vec{0} \color{black}}
\newcommand{\zero}{\color{lightgray!40!Lavender!40!} 0 \color{black}}

\usepackage{circuitikz}
\usepackage{pgfplots, pgfplotstable}
\pgfplotsset{compat=newest}

\usetikzlibrary{external}
\ifexternalize
    \tikzset{external/optimize=false}
\fi
\DeclareRobustCommand{\plotref}[1]{%
    \tikzexternaldisable\ref{#1}\tikzexternalenable}

\newenvironment{authcontrib}[1]{%
    \subsection*{\textnormal{\textbf{Author Contributions}}}%
    \noindent #1}%
{}%
\newenvironment{acks}[1]{%
    \subsection*{\textnormal{\textbf{Acknowledgements}}}%
    \noindent #1}%
{}%
\newenvironment{dci}[1]{%
    \subsection*{\textnormal{\textbf{Declaration of competing interests}}}%
    \noindent #1}%
{}%
\newenvironment{code}[1]{%
    \subsection*{\textnormal{\textbf{Code}}}%
    \noindent #1}%
{}%

\usepackage{booktabs} 

\usepackage[svgnames]{xcolor}

\definecolor{darkblue}{RGB}{51,34,136}
\definecolor{lightblue}{RGB}{136,204,238}
\definecolor{teal}{RGB}{68,170,153}
\definecolor{green}{RGB}{17,119,51}
\definecolor{pistacciogreen}{RGB}{153,153,51}
\definecolor{mustard}{RGB}{221,204,119}
\definecolor{rose}{RGB}{204,102,119}
\definecolor{burgundy}{RGB}{136,34,85}
\definecolor{pink}{RGB}{170,68,153}

\colorlet{color1}{burgundy}
\colorlet{color2}{pistacciogreen}
\colorlet{color3}{darkblue}
\colorlet{color4}{rose}
\colorlet{color5}{mustard}
\colorlet{color6}{lightblue}
\colorlet{color7}{pink}
\colorlet{color8}{teal}
\colorlet{color9}{green}

\colorlet{crv01}{YellowGreen}
\colorlet{crv02}{color1}
\colorlet{crv03}{color6}
\colorlet{crv04}{color2}
\colorlet{crv05}{color3}
\colorlet{crv06}{color4}
\colorlet{crv07}{color5}
\colorlet{crv08}{DarkSeaGreen}
\colorlet{crv09}{Fuchsia}
\colorlet{crv10}{Brown}
\colorlet{crv11}{DarkKhaki}
\colorlet{crv12}{Lavender}
\colorlet{crv13}{Magenta}
\colorlet{crv14}{NavyBlue}
\colorlet{crv15}{Plum}

\pgfplotscreateplotcyclelist{lineplotscheme}{%
    color1,color2,color3,color4,color5,color6,color7,color8,color9}
\pgfplotsset{cycle list name=lineplotscheme}

\def\plotlinewidth{1.0pt}
\def\thickplotlinewidth{2pt}
\def\thinplotlinewidth{0.8pt}
\newlength{\plotwidth}
\newlength{\plotheight}
\usepgfplotslibrary{patchplots} 
\usepgfplotslibrary{colormaps}  

\pgfplotsset{
    legend image code/.code={
            \draw[mark repeat=2,mark phase=2]
            plot coordinates {
                    (0cm,0cm)
                    (0.2cm,0cm)        
                    (0.5cm,0cm)         
                };%
        }
}

\pgfplotsset{
    colormap={redgrad}{
            color(0cm)=(color1!10)
            color(1cm)=(color1)
        }
}

\DeclareRobustCommand{\redgradrule}{%
    \includegraphics{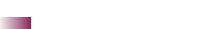}
}

\usepackage[backend=biber, style=numeric-comp,giveninits=true, maxbibnames=8, sortcites=true, sorting=none, url=false, isbn=false, date=year]{biblatex}
\title{Mixed formulation and structure-preserving discretization of Cosserat rod dynamics in a port-Hamiltonian framework}

\author{
\begin{tabular}{ccccc}
Philipp L. Kinon\footnote{Corresponding author: \href{mailto:philipp.kinon@kit.edu}{philipp.kinon@kit.edu}} \textsuperscript{,}\footnote{Karlsruhe Institute of Technology (KIT)} \textsuperscript{,}\footnote{Eindhoven University of Technology (TU/e)} & &
Simon R. Eugster\footnotemark[3] & &  Peter Betsch\footnotemark[2] 
\end{tabular}
}

\date{ }

\begin{document}

\maketitle

\begin{abstract}
    \noindent
    An energy-based modeling framework for the  nonlinear dynamics of spatial Cosserat rods undergoing large displacements and rotations is proposed.
The mixed formulation features independent displacement, velocity and stress variables and is further objective and locking-free. Finite rotations are represented using a director formulation that avoids singularities and yields a constant mass matrix. This results in an infinite-dimensional nonlinear port-Hamiltonian (PH) system governed by partial differential-algebraic equations with a quadratic energy functional. Using a time-differentiated compliance form of the stress-strain relations allows for the imposition of kinematic constraints, such as inextensibility or shear-rigidity.
A structure-preserving finite element discretization leads to a finite-dimensional system with PH structure, thus
facilitating the design of an energy-momentum consistent integration scheme.
Dissipative material behavior (via the generalized-Maxwell model) and non-standard actuation approaches (via pneumatic chambers or tendons) integrate naturally into the framework.
As illustrated by selected numerical examples, the present framework establishes a new approach to energy-momentum consistent formulations in computational mechanics involving finite rotations.

\textbf{Keywords:} Simo-Reissner beam --- Kirchhoff beam --- port-Hamiltonian systems --- differential-algebraic equations ---  structure-preserving discretization --- mixed finite elements
\end{abstract}


\section{Introduction}
\label{sec_intro}
Flexible multibody systems consisting of rigid bodies and flexible components \cite{bauchau_2011_flexible} have a plethora of applications, as for example in the fields of soft robotics \cite{albu-schaffer_2008_soft} or cable dynamics \cite{dorlich_2018_flexible}.
These systems are inherently complex, driven by geometric nonlinearities and intricate material behavior, thus requiring accurate and robust simulation methods \cite{hairer_2006_geometric}.

To this end, the present work provides a novel modeling framework for the dynamics of \emph{Cosserat rods} \cite{antman_2005_nonlinear}, which are a
cornerstone in modeling the nonlinear dynamics of thin, flexible structures.
Hereto, the governing equations are recast as a \emph{port-Hamiltonian} (PH) system \cite{duindam_2009_modeling}.
This approach inherently provides a mixed formulation in terms of independent displacement, velocity and stress variables and facilitates the \emph{structure-preserving discretization} in space and time.

\paragraph{Cosserat rods.}
Cosserat rods are also referred to as Simo-Reissner beams \cite{simo_1985_finitea,reissner_1972_onedimensional} or geometrically exact beams \cite{betsch_2003_constrained,romero_2002_objective}.\footnote{Throughout this work we use the terms \emph{Cosserat rod} and \emph{beam} interchangeably.}
These models have become highly influential \cite{betsch_2003_constrained,romero_2002_objective,crisfield_1999_objectivity,lang_2012_numerical,eugster_2014_directorbased,leyendecker_2006_objective,harsch_2023_nonunit,sonneville_2014_geometrically,herrmann_2024_relativekinematic,wasmer_2024_projectionbased,herrmann_2025_mixed,debeurre_2025_quaternionbased}, as they avoid linear or higher order approximations.
Lately, these models appear also more frequently in soft robotics \cite{rucker_2011_statics,till_2019_realtime,eugster_2022_soft,alessi_2024_rod}.

The complexity of this model is mostly due to its geometric nonlinearity originating from finite rotations, whose parametrization can be approached in various ways \cite{bauchau_2011_flexible}. Possible options to parameterize the rotational degrees of freedom include, but are not limited to Euler angles, (projected) unit quaternions or vectors of relative rotation.
In this work, a singularity-free parametrization will be considered, which avoids coordinates of rotational nature altogether and yields a constant and diagonal mass matrix --- the director framework.
Due to its beneficial properties, the director formulation of Cosserat rods has been a popular choice in the design of finite element methods. The corresponding finite element implementations rely either on the imposition of algebraic nodal constraints \cite{betsch_2003_constrained,eugster_2014_directorbased,leyendecker_2006_objective,betsch_2002_frameindifferent,betsch_2002_dae}
or on the introduction of nodal rotational degrees of freedom
\cite{romero_2002_objective,gruttmann_2000_theory,leyendecker_2008_discrete,huang_2022_electromechanically,eugster_2023_family}.
More details on the choice of rotational parameters can be found in the works \cite{ibrahimbegovic_1997_choice,betsch_1998_parametrization,bauchau_2014_interpolation,romero_2004_interpolation}.

\paragraph{Port-Hamiltonian systems.}

The port-Hamiltonian (PH) framework provides a unifying approach to modeling, simulation, and control of interconnected dynamical systems \cite{duindam_2009_modeling}, ideally suited for flexible multibody systems. The PH description takes into account dissipation and interaction with the environment, and can thus be considered an extension of classical Hamiltonian dynamics.

Originally introduced for finite-dimensional systems in \cite{maschke_1992_portcontrolled}, the extension of the PH modeling paradigm to infinite-dimensional systems can be traced back to \cite{vanderschaft_2002_hamiltonian} and has been subject to extensive research in various areas \cite{rashad_2020_twenty}.
In the context of structural dynamics,
PH formulations have been developed for geometrically exact strings \cite{thoma_2022_porthamiltonian,kinon_2023_porthamiltonian,kinon_2024_generalized}, plates \cite{brugnoli_2019_porthamiltonian,brugnoli_2019_porthamiltoniana} and beams restricted to small or moderate deformations \cite{warsewa_2021_porthamiltonian,macchelli_2004_modeling,brugnoli_2021_mixed,ponce_2024_porthamiltonian,ponce_2025_constrained}.
Nonlinear strain-based beam models have been addressed in \cite{macchelli_2007_portbased,macchelli_2013_stabilisation,ayala_2023_energybased}. Using a material formulation entirely in a body-fixed frame, these formulations are closely related to Hodges' \emph{intrinsic} description of the Cosserat rod \cite{hodges_1990_mixed}. The equivalence to the description commonly used in structural mechanics has been shown in \cite{rodriguez_2022_control} and the related PH representation is detailed in \cite{artola_2022_modalbased}.

\paragraph{Structure-preserving discretization.}
PH formulations inherently provide energy-consistency, thus being ideally suited for the design of structure-preserving discretization methods \cite{hairer_2006_geometric},
which guarantee that the underlying power balance is respected at the discrete level. This avoids numerical dissipation and energy inconsistencies, which is particularly critical in long-term simulations and in control-oriented problems \cite{ortega_2002_interconnection,caasenbrood_2022_energyshaping}.

Numerical schemes that preserve both energy and momentum are commonly referred to as EM schemes.
Although \cite{simo_1995_nonlinear} was the first EM method for spatial Cosserat rods, the important property of objectivity
was violated, as has been shown in \cite{crisfield_1999_objectivity}.
The first EM method for spatial Cosserat rods that retains the objectivity is due to \cite{romero_2002_objective}, see also \cite{betsch_2002_dae,betsch_2003_constrained}, making use of the above-mentioned director framework.

In the context of infinite-dimensional PH systems,
the spatial discretization is likewise achieved through finite element (FE) methods \cite{cardoso-ribeiro_2021_partitioned} yielding finite dimensional PH systems, see also
\cite{brugnoli_2021_porthamiltonian,warsewa_2021_porthamiltonian,brugnoli_2021_mixed,thoma_2022_porthamiltonian,kinon_2023_porthamiltonian,kinon_2025_energymomentumconsistent} for applications in structural dynamics.
It is important to note that the typical use of energy variables in the PH approach inherently gives rise to mixed FE formulations, thus incorporating also independent stress and/or strain fields.
Similarly, time discretization plays a pivotal role
for accurate simulations with energetic consistency. Energy-consistent integration schemes for PH systems
include discrete-gradient formulations \cite{kinon_2023_discrete,kinon_2023_porthamiltonian,goren-sumer_2008_gradient,kinon_2026_discrete} and projection-based approaches \cite{giesselmann_2024_energyconsistent}.

\subsection{Research gaps}
\label{subsec_gaps}

To the best of the authors' knowledge,
mixed formulations and FE methods for the dynamics of Cosserat rods are still quite rare.
Despite the extensive research on this beam formulation and the PH framework, their integration has received only limited attention. The existing strain-based PH formulations
\cite{macchelli_2007_portbased,artola_2022_modalbased,ayala_2023_energybased} require reconstructing absolute position and orientation as it is also frequently the case for contributions outside the PH community \cite{zupan_2003_finiteelement,cesarek_2012_kinematically}. Moreover, the Lie-group setting in \cite{macchelli_2007_portbased} introduces additional geometric complexity, and neither discretizations in \cite{macchelli_2007_portbased,artola_2022_modalbased,ayala_2023_energybased} employ finite elements.
Contrarily, proposed formulations from the computational mechanics community are predominantly restricted to the static scenario \cite{kim_1998_new,wackerfuss_2009_mixed,greco_2024_objective,santos_2010_hybrid,herrmann_2025_mixed,marino_2017_lockingfree}.

Another challenge in the simulation of slender mechanical structures is
transverse shear locking. This numerical stiffening effect arises during discretization, leading to unrealistically small deformations.
While selectively reduced integration \cite{betsch_2002_frameindifferent,romero_2002_objective}
is a popular solution, its treatment using mixed approaches is rarely discussed, see exceptions \cite{kim_1998_new,santos_2010_hybrid,herrmann_2025_mixed} and \cite{brugnoli_2020_porthamiltonian,kinon_2025_energymomentumconsistent} in the PH context.

Linear constitutive relations can be naturally included into the PH formulation using a compliance equation based on a Hellinger-Reissner variational principle \cite{herrmann_2025_mixed,greco_2024_objective,kim_1998_new,santos_2010_hybrid}.
So-called \emph{velocity-stress} formulations, which are popular in the PH community \cite{thoma_2024_velocitystress,brugnoli_2021_mixed,thoma_2022_porthamiltonian}, are related to this concept and consider the compliance equation in time-differentiated form. However, the connection to the Hellinger-Reissner principle has not been explicitly made in these works.
For Cosserat rods, such formulations also bring the advantage of easily preventing certain deformation modes such that also shear-rigid or inextensible beams may be simulated without further modifications of the formulation \cite{herrmann_2025_mixed}. However, an application to the dynamics of Cosserat rods, especially in combination with a PH formulation, is still open.

Lastly, both the description of visco-elasticity with the generalized-Maxwell model \cite{ferri_2023_efficient,weeger_2022_mixed,lestringant_2020_discrete} and the actuation using pneumatic chambers and tendons
\cite{rucker_2011_statics,till_2019_realtime,eugster_2022_soft,till_2019_realtime,alessi_2024_rod}
are central in continuum robotics. Their inclusion into a
PH framework remains an open topic.

\subsection{Contributions}
\label{subsec_contributions}

The findings in the previous works \cite{kinon_2024_generalized,kinon_2025_energymomentumconsistent} indicate that a mixed PH formulation can advance the development of structure-preserving methods for nonlinear structural dynamics. Therein, the specific PH framework \cite{mehrmann_2019_structurepreserving,beattie_2018_linear,mehrmann_2023_control,kinon_2023_discrete} has been proven useful also in the infinite-dimensional setting.
In the present work, we adopt this formulation for spatial Cosserat rods, address the specific challenges associated with finite rotations, and bridge the above-mentioned gaps as follows:
\begin{itemize}
    \item[C1)] We propose a mixed, infinite-dimensional PH formulation for spatial Cosserat rods, which naturally accommodates independent displacement, velocity and stress variables, while ensuring intrinsic energy-momentum consistency. The employed director-parametrization of finite rotations provides constant and diagonal mass matrices.
    \item[C2)] We show the structure-preserving spatial and temporal discretization using mixed finite elements and second-order time-stepping methods, which are facilitated by the mixed PH formulation. Our approach ensures exact discrete representations of the power and angular momentum balance.
    \item[C3)] We demonstrate that the FE approach, induced by the PH framework, mitigates shear locking and straightforwardly includes shear-rigid and inextensible beam formulations as special cases. This can be easily realized in the context of a compliance equation in time-differentiated form.
    \item[C4)] We extend the given PH model with respect to visco-elastic material behavior in terms of a generalized-Maxwell model, and actuation by pneumatic chambers and tendons. These extensions fit well into the proposed PH formulation.
\end{itemize}

\subsection{Notation}
We follow the following notation rules.
Scalars $f,g,n,m,s \in \mathbb{R}$ are represented in italic fonts.
Vectors
and $d$-tuples are treated equally as one-column matrices,
e.g., $\vec{a},\vec{b},\vec{c} \in \R^{d \times 1}$, using bold fonts.
All 3-dimensional vectors representing physical quantities are expressed directly in coordinates, that is, as triples in $\mathbb{R}^3$ relative to a fixed inertial orthonormal basis. Accordingly, the basis vectors of this inertial frame have the canonical coordinate representations $\vec{e}_1=(1,0,0)$, $\vec{e}_2=(0,1,0)$ and $\vec{e}_3=(0,0,1)$.
Making use of the Einstein summation convention for repeated indices, i.e., for Latin indices $i,j,k \in \{1,2,3\}$ and Greek indices $\alpha, \beta \in \{1, 2\}$,
a 3-dimensional vector can be written as
$ \vec{a} = \sum_{i=1}^3 a_i \vec{e}_i = a_i \vec{e}_i = (a_1, a_2, a_3)$.
For convenience, a short bracket notation is used to express $[\vec{a}\transp,\vec{b}\transp]\transp =:(\vec{a},\vec{b})$.
The scalar product is denoted by $\vec{a}\transp \vec{b} = a_i b_i$ and the cross-product in $\mathbb{R}^3$ reads $\vec{c} = \vec{a} \times \vec{b} = \epsilon_{ijk} a_j b_k \vec{e}_i = \tilde{\vec{a}}\vec{b}$, where we made use of the Levi-Civita symbol $\epsilon_{ijk}$ and $\widetilde{\square}$ denotes the mapping between $\R^3$ and the space of skew-symmetric matrices $\mathfrak{so}(3) \subset \mathbb{R}^{3\times 3}$. The Kronecker delta $\delta_{ij}$ yields one for $i=j$ and zero for $i\neq j$.
Matrices $\vec{A} \in \mathbb{R}^{n \times m}$ use bold faces and
differential matrix operators involving derivatives are denoted as calligraphic uppercase letters $\mathcal{A}$.
Matrix transposition is denoted by $\vec{A}\transp$ and the formal adjoint of an operator $\mathcal{A}^*$.
Matrix-vector multiplication and the action of an operator are denoted similarly as $\vec{A} \vec{b}$ and $\mathcal{A} \vec{b}$.
Identity matrices are represented as $\vec{I}$ and matrices full of zeros as $\vec{0}$, their dimension being clear from the context.
The square $\square$ represents a placeholder.
The partial derivative of a function $f$ is denoted $\partial_{\square} f$.
We introduce further
$\square_{,s} := \partial_s \square$ and $\dot{\square} := \partial_t \square$.
The inner product on a one-dimensional domain $\domainspace$, parameterized by $s$, is denoted
$(\vec{a}, \vec{b})_\domainspace := \int_\domainspace \vec{a}\transp \vec{b} \d s$ or $(f, g)_\domainspace := \int_\domainspace f \, g \d s$,
thus taking into account the correct dimension.

\subsection{Outline}
The remainder of this work is structured as follows: In Sec.~\ref{sec_problem} we introduce the problem statement connected to spatial Cosserat rod dynamics. In Sec.~\ref{sec_model} the model is reformulated as a first-order system and its
PH structure is shown, inducing naturally a power balance equation.
Facilitated by the PH formulation, this power balance equation can be carried over to the discrete level by means of structure-preserving discretization as shown in Sec.~\ref{sec_discretization}.
Thereafter, Sec.~\ref{sec_extension} deals with the extension of the formulation to visco-elasticity and two actuation mechanisms, which are relevant in continuum robotics applications.
Ultimately, Sec.~\ref{sec_examples} contains numerical examples, which verify the properties of the proposed approach. These are further recapitulated in the conclusion in Sec.~\ref{sec_conclusion}. Derivations and details, which are not crucial to grasp the main arguments, are comprised in Appendices~\ref{app_math}-\ref{sec_appendix_system_matrices}.

\section{Cosserat rod model}
\label{sec_problem}

\subsection{Kinematics} \label{sec_kinematic}

\begin{figure}[tb]
    \centering
    \def\svgwidth{0.7\textwidth}
    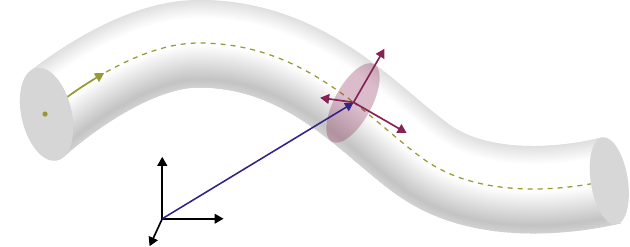
    \caption{Schematic depiction of a Cosserat rod configuration.}
    \label{fig_cosserat_sketch}
\end{figure}

\paragraph{Configuration.}
The following concepts are based on the models in \cite{betsch_2002_frameindifferent,romero_2002_objective}.
Consider a one-dimensional material beam domain $\Omega = [0, L] \subset \R$, with length $L \in \mathbb{R}$, such that the motion of each material point on the centerline (also referred to as \emph{backbone} or \emph{line of centroids}) is described via $\centerline(s,t) \in \mathbb{R}^3$. Here, the material coordinate $s \in \Omega$ is the arc-length parametrization of the centerline of the stress-free initial configuration $\domainspace_0 = \{ \centerline(s,0) = s L \vec{e}_3 | s \in \domainspace \}$ and $t \in [0, \finaltime] \subset \mathbb{R}$ denotes time.  Allowing a geometrically exact description of the kinematics, every cross-section remains plane and rigid during motion and its orientation is described by means of a
rotation matrix
\begin{align} \label{eq_rotmat}
    \rotmat(s,t) = \big[ \vec{d}_1(s,t) \quad \vec{d}_2(s,t) \quad \vec{d}_3(s,t)\big]= \vec{d}_k(s,t) \vec{e}_k\transp ,
\end{align}
from $SO(3)$, the special orthogonal group. The columns of $\rotmat$ form a
local frame of orthonormal directors $\vec{d}_i(s,t) \in \R^3$ for $i=1,2,3$. See Fig.~\ref{fig_cosserat_sketch} for a depiction of the kinematics.
The director $\vec{d}_3(s,t)$ is assumed to be normal to each cross section and $\vec{d}_1(s,t), \vec{d}_2(s,t)$ are transverse to it, such that $\vec{d}_3(s,t) = \vec{d}_1(s,t) \times \vec{d}_2(s,t)$.
For simplicity, we further assume that $\vec{d}_1(s,t)$ and $\vec{d}_2(s,t)$ coincide with the principal axes of inertia of the cross-section at all times.
In the following, the explicit dependency on space and time is omitted and only mentioned where deemed necessary.
Abusing notation, we also interpret $\rotmat=\rotmat(\vec{d})$ as function of the vertical concatenation of the directors $\vec{d}=(\vec{d}_1,\vec{d}_2,\vec{d}_3) \in \R^9$. Furthermore, we call $(\centerline,\vec{d})$ \emph{displacement} quantities, as they uniquely characterize the displacement of each material point on the beam configuration.
Since $\rotmat(\vec{d}) \in SO(3)$, it holds that $\det(\vec{R}(\vec{d}))=+1$ and $\rotmat(\vec{d})\transp \rotmat(\vec{d}) = \rotmat(\vec{d}) \rotmat(\vec{d})\transp = \vec{I}$.
Consequently, there are 6 orthonormality conditions
\begin{align} \label{eq:ortho_constraint}
    \vec{0} = \vec{g}(\vec{d}) \in \mathbb{R}^6 \quad \Leftrightarrow \quad \frac{1}{2} \left( \vec{d}_i \transp \vec{d}_j - \delta_{ij} \right) = 0 ,
\end{align}
for $(i,j) \in \{(1,1),(2,2),(3,3),(1,2),(2,3),(3,1)\}$.

\paragraph{Strains.}
We make use of the objective, material strain measures $ \vec{\Gamma}=\rotmat\transp \centerline_{,s}$ to describe shear deformation and dilatation and $\widetilde{\vec{K}} = \rotmat\transp \rotmat_{,s}$ for curvature and torsion,
which can be expressed equivalently as
\begin{equation} \label{eq:strain_measures}
    \vec{\Gamma} = (\vec{d}_k \transp \centerline_{,s})\vec{e}_k = \begin{bmatrix}
        \vec{d}_1 \transp \centerline_{,s} \\
        \vec{d}_2 \transp \centerline_{,s} \\
        \vec{d}_3 \transp \centerline_{,s}
    \end{bmatrix}
    , \qquad \qquad
    \vec{K}
    = \frac{1}{2} \epsilon_{ijk} (\vec{d}_k \transp \vec{d}_{j,s}) \vec{e}_i
    = \frac{1}{2}\begin{bmatrix}
        \vec{d}_{3}\transp \vec{d}_{2,s} - \vec{d}_{2}\transp \vec{d}_{3,s}  \\
        \vec{d}_{1}\transp \vec{d}_{3,s} -  \vec{d}_{3}\transp \vec{d}_{1,s} \\
        \vec{d}_{2}\transp \vec{d}_{1,s} - \vec{d}_{1}\transp \vec{d}_{2,s}
    \end{bmatrix},
\end{equation}

where we switched to the axial vector of $\widetilde{\vec{K}}$.
Due to the bilinearity of both $\vec{\Gamma}(\vec{d},\centerline_{,s})$ and $\vec{K}(\vec{d},\vec{d}_{,s})$, one can compactly rewrite them using two different representations each, i.e.,
\begin{align} \label{eq:strain_measures_short}
    \vec{\Gamma} = \vec{R}(\vec{d})\transp \centerline_{,s} = \vec{J}_{\vec{N}}(\centerline_{,s})\transp\vec{d} , \qquad  \qquad
    \vec{K} = \vec{L}(\vec{d}_{,s})\vec{d} = - \vec{L}(\vec{d})\vec{d}_{,s} .
\end{align}
The rotation matrix \eqref{eq_rotmat} as well as the newly introduced matrices
\begin{align} \label{eq_JN_L}
    \vec{J}_{\vec{N}}(\centerline_{,s})  = \begin{bmatrix}
                                               \centerline_{,s} &  & \zeros           &  & \zeros           \\
                                               \zeros           &  & \centerline_{,s} &  & \zeros           \\
                                               \zeros           &  & \zeros           &  & \centerline_{,s} \\
                                           \end{bmatrix} \in \mathbb{R}^{9 \times 3} ,  \qquad
    \vec{L}(\vec{d})                                   = \frac{1}{2}
    \begin{bmatrix}
        \zeros              &  & -\vec{d}_{3}\transp &  & \vec{d}_{2}\transp  \\
        \vec{d}_{3}\transp  &  & \zeros              &  & -\vec{d}_{1}\transp \\
        -\vec{d}_{2}\transp &  & \vec{d}_{1}\transp  &  & \zeros
    \end{bmatrix} \in \mathbb{R}^{3 \times 9}
\end{align}
are linear functions of their arguments. Note that $\vec{L}(\vec{d})\vec{d}=\vec{0}$ due to \eqref{eq:ortho_constraint}.

\begin{remark} \label{rem_initial_conf}
    We assume straight and stress-free initial configurations that serve as the reference configuration. Thus, the cross-section normal $\vec{d}_3(s,t=0)$ aligns with the centerline tangent vector $\centerline_{,s}(s,t=0)$ implying the homogeneous reference strains $\vec{\Gamma}_0 :=  \vec{\Gamma}(s,t=0) = (0,0,1) $ and $\vec{K}_0 := \vec{K}(s,t=0) = \vec{0}$. For pre-curved beams, the reference strains would be different, which does not affect the generality of the theoretical framework developed in the remainder of this work.
    Moreover, $\vec{\Gamma}$ describes the rate of change of the centerline (i.e., the tangent vector), while
    $\vec{K}$ describes the rate of change of orientation when moving along the beam. The fact that $\vec{K}$ is derived from a skew-symmetric matrix $\widetilde{\vec{K}}$ can be traced back to taking the derivative of the orthonormality constraints \eqref{eq:ortho_constraint} with respect to $s$.
\end{remark}

\subsection{Constitutive modeling}
Starting from the general case of hyperelastic material models based on a scalar-valued strain energy function $W = \breve{W}(\vec{\Gamma},\vec{K})$ being a function of the strain measures, the corresponding work-conjugated stress quantities are derived as
\begin{align} \label{eq_material_standard}
    \vec{N} = \partial_{\vec{\Gamma}} \breve{W} , \qquad \qquad \vec{M} = \partial_{\vec{K}} \breve{W} .
\end{align}
Here, indicated by an uppercase notation, $\vec{N} = N_i \vec{e}_i$ are the material normal and shear forces and $\vec{M} = M_i \vec{e}_i$ are the material bending and torsion moments. The spatial stress quantities $\vec{n}$ and $\vec{m}$ can be linked to their material counterparts via
\begin{align} \label{trafo_spat_forces}
    \vec{n} = \rotmat(\vec{d}) \vec{N} = N_i \vec{d}_i , \qquad  \qquad \vec{m} = \rotmat(\vec{d}) \vec{M} = M_i \vec{d}_i .
\end{align}
In the following, we assume linear constitutive laws such that $W$ has quadratic form. This is a common assumption in the literature since locally only small strains are expected.
We thus write
\begin{align} \label{quadratic_stored_energy}
    \breve{W}(\vec{\Gamma},\vec{K}) = \frac{1}{2} (\vec{\Gamma}-\vec{\Gamma}_0)\transp \vec{C}_\Gamma (\vec{\Gamma}-\vec{\Gamma}_0) + \frac{1}{2} (\vec{K}-\vec{K}_0)\transp \vec{C}_K (\vec{K}-\vec{K}_0)
\end{align}
with constant elasticity matrices $\vec{C}_\Gamma = \diag(k_{\mathrm{s1}}, k_{\mathrm{s2}}, k_{\mathrm{e}}) \in \R^{3\times 3}$ and
$\vec{C}_K = \diag(k_{\mathrm{b1}},k_{\mathrm{b2}},k_{\mathrm{t}}) \in \R^{3\times 3}$ as well as $\vec{\Gamma}_0$ and $\vec{K}_0$ according to Remark~\ref{rem_initial_conf}.

\paragraph{Compliance form.}
A quadratic hyperelastic strain energy function can also be expressed in terms of the stresses as \emph{complementary energy function} $W = W^\star(\vec{N},\vec{M})$ such that
\begin{align} \label{eq_comp_strain}
    W^\star(\vec{N},\vec{M}) = \frac{1}{2} \vec{N}\transp \vec{C}_N \vec{N} + \frac{1}{2} \vec{M}\transp \vec{C}_M \vec{M}
\end{align}
with compliance matrices $\vec{C}_N = \diag(k_{\mathrm{s1}}^{-1}, k_{\mathrm{s2}}^{-1}, k_{\mathrm{e}}^{-1}) \in \R^{3\times 3}$, $\vec{C}_M = \diag(k_{\mathrm{b1}}^{-1},k_{\mathrm{b2}}^{-1},k_{\mathrm{t}}^{-1}) \in \R^{3\times 3}$. Note that these matrices are allowed to be singular if some stiffnesses tend to infinity, e.g., $k_{\mathrm{s1}}, k_{\mathrm{s2}} \rightarrow \infty$ and thus $\vec{C}_N = \diag(0,0,k_{\mathrm{e}}^{-1})$ for shear-stiff beam theories, like the well-known \emph{Kirchhoff}-beam model. Similar to \eqref{eq_material_standard}, the strain measures can be derived by differentiation, such that
\begin{align} \label{eq:material}
    \vec{\Gamma} - \vec{\Gamma}_0 = \partial_{\vec{N}} W^\star = \vec{C}_N \vec{N}  , \qquad  \qquad \vec{K}-\vec{K}_0 = \partial_{\vec{M}} W^\star = \vec{C}_M \vec{M} .
\end{align}
In the case of hyperelastic stress-strain relations that cannot be inverted, a Hu-Washizu type approach \cite{kinon_2023_porthamiltonian,kinon_2025_energymomentumconsistent} might be required.

\subsection{Equations of motion} \label{sec_eom}
The partial differential equations governing the dynamics of the Cosserat rod
can be derived from the principle of virtual work, see Appendix~\ref{virtual_work_appendix}, as
\begin{equation} \label{EoM_Antman_both}
    \rho A \ddot{\centerline} = \vec{n}_{,s}  + \bar{\vec{n}}  , \qquad  \qquad
    \vec{I}_\rho \dot{\vec{\omega}} + \vec{\omega} \times \vec{I}_\rho \vec{\omega}  = \vec{m}_{,s} + \centerline_{,s} \times \vec{n} + \bar{\vec{m}} ,
\end{equation}
stating the local balance of linear and angular momentum in spatial form, respectively \cite{simo_1985_finitea,antman_2005_nonlinear}.
Here, we introduced the uniformly distributed mass density per unit length $\rho \geq 0$ and cross section area $A > 0$,
$\vec{\omega}(s,t) \in \mathbb{R}^3$ is the spatial angular velocity,
$\vec{I}_\rho \in \mathbb{R}^{3 \times 3}$ is the spatial representation of the inertia tensor,
and $\bar{\vec{n}}(s,t), \bar{\vec{m}}(s,t) \in \mathbb{R}^3$ are external distributed forces and moments per unit length.

In contrast to \eqref{EoM_Antman_both}, we target a formulation in terms of directors as introduced above. Considering the derivations from Appendix~\ref{virtual_work_appendix},
one obtains the equivalent representation of the Cosserat rod dynamics in the form
\begin{subequations} \label{eq:balance_ang_momentum_1}
    \begin{align}
        \rho A \ddot{\centerline}   & = \vec{n}_{,s}  + \bar{\vec{n}}  ,                            \label{eq:balance_ang_momentum_1_0}                                                          \\
        \vec{M}_\rho \ddot{\vec{d}} & = - \vec{J}_{\vec{N}} \vec{N} - \operator{J}_{\vec{M}}\vec{M} - \vec{G}_d\transp \vec{\lambda} + \vec{T}\bar{\vec{m}}, \label{eq:balance_ang_momentum_1_1} \\
        \vec{g}(\vec{d})            & =\vec{0}, \label{eq:balance_ang_momentum_2}
    \end{align}
\end{subequations}
where $\vec{J}_{\vec{N}}$ has been introduced in \eqref{eq_JN_L}. Further, we have introduced a constant and diagonal mass matrix $\vec{M}_\rho =\diag(M_\rho^{11}\vec{I},M_\rho^{22}\vec{I},\vec{0}) \in \R^{9\times 9}$,
and the operator matrix
$\operator{J}_{\vec{M}}$, which is defined together with its adjoint $\operator{J}_{\vec{K}}=\operator{J}_{\vec{M}}^*$
as
\begin{align} \label{eq_BM_BK}
    \operator{J}_{\vec{M}}(\vec{d},\vec{d}_{,s})\square  =
    \vec{L}(\vec{d}_{,s})\transp \square + \partial_s(\vec{L}(\vec{d})\transp \square \,),  \qquad \qquad
    \operator{J}_{\vec{K}}(\vec{d},\vec{d}_{,s})\square     =  \vec{L}(\vec{d}_{,s}) \square - \vec{L}(\vec{d})\partial_s \square \ .
\end{align}
Here, $\vec{L}$ is given in \eqref{eq_JN_L}.
Additionally in \eqref{eq:balance_ang_momentum_1_1}, we have introduced
$\vec{T}\in \R^{9\times 3}$ defined as
\begin{equation} \label{eq_T}
    \vec{T}(\vec{d}) = -\frac{1}{2} \begin{bmatrix}
        \tilde{\vec{d}}_1 \\ \tilde{\vec{d}}_2 \\ \tilde{\vec{d}}_3
    \end{bmatrix} , \qquad \text{with} \qquad \vec{T}(\vec{d})\transp = \frac{1}{2} \big[
        \tilde{\vec{d}}_1 \quad \tilde{\vec{d}}_2 \quad \tilde{\vec{d}}_3
        \big] ,
\end{equation}
which connects angular velocities and director velocities, see \eqref{eq_tmp2} in Appendix~\ref{app_math}. Lastly,
the constraint gradient matrix
\begin{align} \label{constraint_jacobian}
    \vec{G}_d(\vec{d}) = \partial_{\vec{d}}\vec{g}(\vec{d}) =
    \begin{bmatrix}
        \vec{d}_1\transp             &  & \zeros                       &  & \zeros                       \\
        \zeros                       &  & \vec{d}_2\transp             &  & \zeros                       \\
        \zeros                       &  & \zeros                       &  & \vec{d}_3\transp             \\
        \frac{1}{2} \vec{d}_2\transp &  & \frac{1}{2} \vec{d}_1\transp &  & \zeros                       \\
        \zeros                       &  & \frac{1}{2} \vec{d}_3\transp &  & \frac{1}{2} \vec{d}_2\transp \\
        \frac{1}{2} \vec{d}_3\transp &  & \zeros                       &  & \frac{1}{2} \vec{d}_1\transp \\
    \end{bmatrix}  \in \R^{6 \times 9}
\end{align}
and the vector containing all Lagrange multipliers $\vec{\lambda} \in \R^6$
encode constraint forces in \eqref{eq:balance_ang_momentum_1_1}.
These forces ensure the fulfillment of the
orthonormality constraints \eqref{eq:balance_ang_momentum_2}, which are equivalent to \eqref{eq:ortho_constraint}.
Note that $\vec{N}$ and $\vec{M}$ are calculated by means of $\centerline$ and $\vec{d}$ via \eqref{eq_material_standard} and \eqref{eq:strain_measures_short}. The related initial boundary value problem is completed by suitable initial and boundary conditions.

The purely displacement-based formulation \eqref{eq:balance_ang_momentum_1} can be prone to locking and more tedious to deal with in the context of energy-momentum consistent discretization. We thus target a mixed, port-Hamiltonian description in the following section.

\section{Mixed formulation as a port-Hamiltonian system}
\label{sec_model}

\subsection{Local form as first order system} \label{sec_first_order}

Our goal is to rewrite the problem at hand as a set of first-order partial-differential equations with algebraic constraints due to the orthonormality of the directors. To this end, we introduce independent stress and velocity variables, linked via a compliance equation in time-differentiated form, and also consider the orthonormality constraints on velocity level.

Let us introduce independent centerline velocities $\vphi = \dot{\centerline}$ and director-related velocities $\vec{v}_d = \dot{\vec{d}}$ as well as
the displacement and velocity vectors $\vec{q}, \vec{v} \in \R^{12}$ such that
\begin{align} \label{eq:kinematics}
    \vec{q} =
    \begin{bmatrix}
        \centerline \\
        \vec{d}
    \end{bmatrix}
    , \quad
    \vec{v} = \begin{bmatrix}
                  \vphi \\
                  \vec{v}_d
              \end{bmatrix} \qquad \text{with} \qquad
    \dot{\vec{q}} = \vec{v} .
\end{align}
By making use of the newly introduced velocities,
we can rewrite the balance equations \eqref{eq:balance_ang_momentum_1_0} and \eqref{eq:balance_ang_momentum_1_1} as
\begin{equation}
    \label{eom_rewritten}
    \rho A \dotvphi  = \partial_s (\rotmat\vec{N}) + \bar{\vec{n}} ,   \qquad  \qquad
    \vec{M}_\rho \dot{\vec{v}}_d = - \vec{J}_{\vec{N}} \vec{N} - \operator{J}_{\vec{M}}\vec{M} - \vec{G}_d\transp \vec{\lambda} + \vec{T}\bar{\vec{m}} .
\end{equation}
Next, taking the time derivative of the compliance equation \eqref{eq:material} together with \eqref{eq:strain_measures_short}
we obtain an evolution equation for the stresses as
\begin{equation} \label{eq:strain_evo_stress_resultants}
    \vec{C}_N \dot{\vec{N}}  = \rotmat(\vec{d})\transp \vphis + \vec{J}_{\vec{N}}(\centerline_{,s})\transp \vec{v}_d , \qquad  \qquad
    \vec{C}_M \dot{\vec{M}}  = \operator{J}_{\vec{K}}(\vec{d},\vec{d}_{,s}) \vec{v}_d .
\end{equation}
Note that we have used the operator matrix $\operator{J}_{\vec{K}}(\vec{d},\vec{d}_{,s})$ in \eqref{eq:strain_evo_stress_resultants}, which is defined in \eqref{eq_BM_BK}.
Lastly,
the use of directors as independent quantities requires the explicit enforcement of 6 orthonormality constraints
\eqref{eq:balance_ang_momentum_2}, which we will also consider in
time-differentiated form
as
\begin{align} \label{eq:constraint_director}
    \vec{G}(\vec{q})\vec{v} = \begin{bmatrix} \vec{0} & \vec{G}_d(\vec{d}) \end{bmatrix} \begin{bmatrix} \vphi \\ \vec{v}_d \end{bmatrix}
    = \vec{0} ,
\end{align}
where the constraint matrix $\vec{G}_d(\vec{d})$ is defined in \eqref{constraint_jacobian}.
Eventually, together with suitable initial and boundary conditions,
equations \eqref{eq:kinematics}, \eqref{eom_rewritten}, \eqref{eq:strain_evo_stress_resultants} and \eqref{eq:constraint_director} are the basis for a mixed, PH formulation of the Cosserat rod dynamics.
For the sake of brevity, let us introduce
\begin{equation} \label{eq_abbrev}
    \begin{aligned}
        \vec{\sigma} = \begin{bmatrix}
                           \vec{N} \\ \vec{M}
                       \end{bmatrix} , \ \
        \vec{C} = \begin{bmatrix}
                      \vec{C}_N & \vec{0}   \\
                      \vec{0}   & \vec{C}_M
                  \end{bmatrix}
        , \ \
        \vec{M}_o = \begin{bmatrix}
                        \rho A \vec{I} & \vec{0}      \\
                        \vec{0}        & \vec{M}_\rho
                    \end{bmatrix} , \ \
        \vec{B}_{\domainspace} = \begin{bmatrix}
                                     \vec{I} & \vec{0} \\
                                     \vec{0} & \vec{T}
                                 \end{bmatrix} , \ \
        \vec{u}_\Omega = \begin{bmatrix} \bar{\vec{n}} \\ \bar{\vec{m}} \end{bmatrix} , \\
        \operator{J}_{\vec{\sigma}}(\vec{q}) = \begin{bmatrix}
                                                   \rotmat(\vec{d})\transp \partial_s \square & \vec{J}_{\vec{N}}(\centerline_{,s})\transp   \\
                                                   \vec{0}                                    & \operator{J}_{\vec{K}}(\vec{d},\vec{d}_{,s})
                                               \end{bmatrix} ,
        \ \
        \operator{J}_{\vec{v}}(\vec{q}) = \begin{bmatrix}
                                              -\partial_s (\rotmat(\vec{d}) \square) & \vec{0}                                      \\
                                              \vec{J}_{\vec{N}}(\centerline_{,s})    & \operator{J}_{\vec{M}}(\vec{d},\vec{d}_{,s})
                                          \end{bmatrix} , \
    \end{aligned}
\end{equation}
to compactly rewrite \eqref{eom_rewritten} and \eqref{eq:strain_evo_stress_resultants} in the following section.

\subsection{Port-Hamiltonian representation} \label{sec_pH_conti}
The PH representation of the governing partial differential-algebraic equations at hand comprises equations
\eqref{eq:kinematics}, \eqref{eom_rewritten}, \eqref{eq:strain_evo_stress_resultants} and \eqref{eq:constraint_director}. Using the definitions \eqref{eq_abbrev},
we can rewrite them as
\begin{align}
    \label{compliance_form}
    \begin{bmatrix}
        \vec{I} & \zeros    & \zeros  & \zeros  \\
        \zeros  & \vec{M}_o & \zeros  & \zeros  \\
        \zeros  & \zeros    & \vec{C} & \zeros  \\
        \zeros  & \zeros    & \zeros  & \vec{0}
    \end{bmatrix}
    \begin{bmatrix}
        \dot{\vec{q}}      \\
        \dot{\vec{v}}      \\
        \dot{\vec{\sigma}} \\
        \dot{\vec{\lambda}}
    \end{bmatrix} & =
    \begin{bmatrix}
        \zeros    & \vec{I}                              & \zeros                           & \zeros                   \\
        - \vec{I} & \zeros                               & -\operator{J}_{\vec{v}}(\vec{q}) & -\vec{G}(\vec{q})\transp \\
        \zeros    & \operator{J}_{\vec{\sigma}}(\vec{q}) & \zeros                           & \zeros                   \\
        \zeros    & \vec{G}(\vec{q})                     & \zeros                           & \zeros
    \end{bmatrix}
    \begin{bmatrix}
        \vec{0}      \\
        \vec{v}      \\
        \vec{\sigma} \\
        \vec{\lambda}
    \end{bmatrix}
    +
    \begin{bmatrix}
        \zeros                 \\
        \vec{B}_{\domainspace} \\
        \zeros                 \\
        \zeros
    \end{bmatrix}
    \vec{u}_\Omega ,
\end{align}
where the so-called \emph{descriptor} matrix $\vec{E}=\diag(\vec{I},\vec{M}_o,\vec{C},\vec{0})$ is constant, symmetric and rank deficient by at least 9 due to the 6 algebraic constraints and the singular mass matrix.
If certain stiffnesses are assumed to tend to infinity (e.g., by describing a shear-rigid beam), the rank deficiency increases accordingly due to a singular compliance matrix $\vec{C}$.
In short form one can rewrite \eqref{compliance_form} in the PH framework by \cite{beattie_2018_linear,mehrmann_2019_structurepreserving,mehrmann_2023_control}\footnote{Also see the adaptions in the context of infinite-dimensional systems in \cite{kinon_2025_energymomentumconsistent,kinon_2023_porthamiltonian,kinon_2024_generalized}.} as
\begin{equation} \label{phDAE_compact}
    \begin{aligned}
        \vec{E}\dot{\vec{x}} & = \operator{J}(\vec{x})\vec{z}(\vec{x}) + \vec{B}(\vec{x})\vec{u}_\Omega , \\
        \vec{y}_\domainspace & = \vec{B}(\vec{x})\transp \vec{z}(\vec{x}),
    \end{aligned}
\end{equation}
where the state vector $\state \in \mathbb{R}^{36}$ comprises the primary unknowns
\begin{equation} \label{eq_state}
    \state = (\vec{q},\vec{v},\vec{\sigma},\vec{\lambda}) .
\end{equation}
The other newly introduced quantities in \eqref{phDAE_compact} can be identified by comparison with \eqref{compliance_form}. Moreover,
$\operator{J}(\vec{x})$ is a differential operator matrix, termed \emph{structure} operator matrix, which involves partial derivatives.
Further, the costate function $\vec{z}$ has been introduced and assumes the linear form $\vec{z}(\vec{x}) =
    \vec{Q} \vec{x}$ with $\vec{Q} = \diag( \vec{0}, \vec{I}, \vec{I}, \vec{I} )$.
In \eqref{phDAE_compact}, we have appended the distributed output relation, yielding $\vec{y}_\domainspace = (\vphi, \vec{\omega})$,
i.e., the centerline velocity and the spatial angular velocity $\vec{\omega}$. By definition this distributed output is power-conjugated to the distributed input forces and torques $\vec{u}_\domainspace=(\bar{\vec{n}},\bar{\vec{m}})$.

\paragraph{Port-Hamiltonian structure.}
At this point, it remains to show that the above equations possess a PH structure. This requires
that (i) the costate function $\vec{z}$ is properly linked to the underlying Hamiltonian of the system, and that (ii) the structure operator matrix $\operator{J}$ is formally skew-adjoint, i.e., $\operator{J}(\vec{x}) = -\operator{J}(\vec{x})^*$. These two conditions ensure that the PH system has an underlying Stokes-Dirac structure and consequently satisfies a power-balance equation.

Condition (i) requires that
\begin{align} \label{eq:costate}
    \vec{E}\transp \vec{z}(\vec{x}) = \delta_{\vec{x}} H ,
\end{align}
where $\delta_{\vec{x}} H$ is the variational derivative of a suitable Hamiltonian $H$, see \cite{beattie_2018_linear,mehrmann_2019_structurepreserving,mehrmann_2023_control}. With the above definitions of $\vec{E}$ and $\vec{z}(\vec{x})$ this is indeed the case for the Hamiltonian given
by
\begin{align} \label{hamiltonian_director_stress}
    H(\vec{x}) & = \int_\domainspace h(\state) \d s  = \frac{1}{2} \int_\domainspace \vec{x}\transp \vec{Q}\transp \vec{E} \vec{x} \d s = \frac{1}{2} \int_\domainspace \left(\vec{v}\transp \vec{M}_o \vec{v} +  \vec{\sigma}\transp \vec{C} \vec{\sigma}\right) \d s ,
\end{align}
which represents the total mechanical energy of the system, i.e., the sum of kinetic energy and the complementary energy stored in the beam due to deformation.
Consequently, $ \delta_{\vec{x}} H = \partial_{\vec{x}} h = \vec{Q}\transp \vec{E} \vec{x}$ and \eqref{eq:costate} holds true since
$\vec{Q}\transp \vec{E} = \vec{E}\transp \vec{Q}$.
Note that the Hamiltonian does not depend on $\vec{q}$ and $\vec{\lambda}$, since these do not contribute to the energy of the system. These properties are reflected by the corresponding zero blocks in the matrices $\vec{Q}$ and $\vec{E}$, respectively.

Concerning condition (ii),
the property of formal skew-adjointness of $\operator{J}$ in \eqref{compliance_form} can be regarded as the extension of skew-symmetry for matrices to differential matrix operators. Loosely speaking, this means that the computation of $\inner{\vec{z}}{\operator{J}\vec{z}}$ may only yield boundary terms after integration by parts\footnote{
    For an operator $\operator{A}$ its adjoint $\operator{A}^*$ is defined by $\inner{\vec{a}}{\operator{A}\vec{b}} = \inner{\operator{A}^* \vec{a}}{\vec{b}}$ for all $\vec{a},\vec{b}$ in the domain of $\operator{A}$.
    A formally skew-adjoint operator satisfies $\inner{\vec{a}}{\operator{J}\vec{b}} = - \inner{\operator{J}\vec{a}}{\vec{b}}$ for vanishing boundary terms for all $\vec{a},\vec{b}$ in the domain of $\operator{J}$.
    For more details, see \cite{brugnoli_2020_porthamiltonian} or Chap.~4 in \cite{duindam_2009_modeling}.
}.
For elements that do not involve spatial derivatives, such as $\vec{I}$ and $\vec{G}$, this is straightforward by checking the skew-symmetry. However, for the differential operators $\operator{J}_{\vec{\sigma}}$ and $\operator{J}_{\vec{v}}$ this requires some more attention.
Eventually, the result
\begin{equation} \label{boundary_terms}
    \inner{\vec{z}}{\operator{J} \vec{z}} = \left[\vphi \transp \vec{n} + \vec{\omega} \transp \vec{m} \right]_{\partial \domainspace} =: \vec{y}_{\partial \domainspace}\transp \vec{u}_{\partial \domainspace}
\end{equation}
is obtained along the lines of
Appendix~\ref{app_skew_adjoint}.
This product encodes the power flow through the boundary $\partial \domainspace = \{0, L \}$ exerted by means of the outputs $\vec{y}_{\partial \domainspace}$ and boundary inputs $\vec{u}_{\partial \domainspace}$. The latter are given by the boundary values of the initial boundary value problem.
For example, in the case of the pure Neumann boundary conditions with prescribed quantities $\square_{\mathrm{N}}$ on the Neumann boundary $\partial \domainspace_{\mathrm{N}} = \{0,L\}$, one has
\begin{align} \label{boundary_input_output}
    \vec{u}_{\partial \domainspace} = \begin{bmatrix}
                                          -\vec{n}_{\mathrm{N},0}(t) \\
                                          -\vec{m}_{\mathrm{N},0}(t) \\
                                          \vec{n}_{\mathrm{N},L}(t)  \\
                                          \vec{m}_{\mathrm{N},L}(t)  \\
                                      \end{bmatrix}, \qquad \qquad \vec{y}_{\partial \domainspace} = \begin{bmatrix}
                                                                                                         \vphi(0,t)        \\
                                                                                                         \vec{\omega}(0,t) \\
                                                                                                         \vphi(L,t)        \\
                                                                                                         \vec{\omega}(L,t)
                                                                                                     \end{bmatrix} .
\end{align}

\paragraph{Power balance equation.}
Since both conditions have been verified, the governing equations \eqref{phDAE_compact} have PH structure. As a consequence,
an underlying power balance equation with respect to the Hamiltonian \eqref{hamiltonian_director_stress} is induced as
\begin{align} \label{power_balance_cont}
    \dot{H} =
    \int_\domainspace \delta_{\vec{x}} H\transp \dot{\vec{x}} \d s =
    \int_\domainspace \vec{z}\transp \vec{E} \dot{\vec{x}} \d s =
    \int_\domainspace \vec{z}\transp \left(  \operator{J}\vec{z} + \vec{B}\vec{u}_\Omega \right) \d s =
    \vec{y}_{\partial \domainspace}\transp \vec{u}_{\partial \domainspace} + \int_\domainspace \vec{y}_\domainspace\transp \vec{u}_\domainspace \d s ,
\end{align}
where we have used the chain rule, \eqref{phDAE_compact}, \eqref{eq:costate} and \eqref{boundary_terms} and omitted the state-dependency of the respective quantities.
Note that the power balance \eqref{power_balance_cont} automatically ensures conservation of energy for closed systems (i.e., for $\vec{u}_\domainspace = \vec{0}$, $\vec{u}_{\partial \domainspace} = \vec{0}$) and demonstrates the property of passivity. This property lies at the core of PH systems and can be exploited when designing energy-momentum consistent methods.
Thus, when it comes to the numerical discretization, we aim to preserve the PH structure to retain this property in simulations.

\begin{remark}
    The underlying mathematical notion is that of a Stokes-Dirac structure, which defines so-called \emph{flow} and \emph{effort} variables on respective dual spaces. As pointed out for finite-dimensional systems of the same structure \cite{beattie_2018_linear}, the present framework is consistent with that geometric viewpoint on PH systems. In fact, we can define flows $\vec{f} = \vec{E} \dot{\vec{x}}$ and efforts $\vec{e} = \vec{z}(\vec{x}) = \vec{Q}\vec{x}$. Here, the Stokes-Dirac structure is modulated by the state, which can be seen from the state-dependency of the structure matrix operator $\operator{J}(\state)$. An in-depth explanation of that viewpoint is omitted here for the sake of conciseness. For more detailed explanations, the interested reader is referred to Chap.~4 in \cite{duindam_2009_modeling}, or the instructive works by Brugnoli and co-authors \cite{brugnoli_2019_porthamiltonian,brugnoli_2021_porthamiltonian,brugnoli_2020_porthamiltonian}.
\end{remark}

\subsection{Weak form} \label{sec_weak_form}

The weak form of the problem can be derived in a straightforward manner from the differential form \eqref{compliance_form} and serves as the basis for the following FE discretization. Accordingly, premultiplying on both sides by appropriate test functions
$(\vec{M}_o \delta \vec{v}, \delta \vec{q}, \delta \vec{\sigma}, \delta \vec{\lambda})$, integrating over the spatial domain $\domainspace$ and applying integration by parts yields the detailed weak relations
\begin{subequations} \label{weak_form_directors}
    \begin{align}
        \inner{\rho A \delta \vphi}{\dot{\centerline}- \vphi}            & =0 , \label{weak_form_directors_1}       \\
        \inner{\vec{M}_\rho \delta \vec{v}_d}{\sdot{\vec{d}}- \vec{v}_d} & =0 , \label{weak_form_directors_2}       \\
        \inner{\delta \centerline}{\rho A \dotvphi - \bar{\vec{n}}}
        +  \inner{\delta \centerline_{,s}}{ \rotmat\vec{N}}
        - \innerBoundary{\delta \centerline}{\vec{n}_N}                  & =0 , \label{weak_form_directors_3}       \\
        \inner{\delta \vec{d}}{\vec{M}_\rho \dot{\vec{v}}_d +  \vec{J}_{\vec{N}} \vec{N} + \vec{G}_d\transp \vec{\lambda}
            - \vec{T}\bar{\vec{m}}}
        + \inner{\operator{J}_{\vec{K}}\delta \vec{d}}{ \vec{M}}
        - \innerBoundary{\delta \vec{d}}{\vec{T}\vec{m}_N}               & =0 , \label{weak_form_directors_4}       \\
        \inner{\delta \vec{N}}{\vec{C}_N \dot{\vec{N}} -
            \rotmat\transp \vphis - \vec{J}_{\vec{N}}\transp \vec{v}_d
        }                                                                & =0 ,     \label{weak_form_directors_5}   \\
        \inner{\delta \vec{M}}{\vec{C}_M \dot{\vec{M}} -
            \operator{J}_{\vec{K}} \vec{v}_d
        }                                                                & =0 ,       \label{weak_form_directors_6} \\
        \inner{\delta \vec{\lambda}}{\vec{G}_d \vec{v}_d
        }                                                                & =0 ,\label{weak_form_directors_7}
    \end{align}
\end{subequations}
for all sufficiently smooth test functions.
This weak form can also be linked to a Hellinger-Reissner principle, as shown in Appendix~\ref{sec_HR}.

\begin{remark}\label{rem_BCs}
    The weak form \eqref{weak_form_directors} is valid if the Dirichlet boundary conditions are considered in the space of admissible test functions. However, if the Dirichlet boundary conditions are enforced in a weak sense, the weak form must be slightly modified, e.g., by introducing additional Lagrange multipliers and accounting for the prescribed boundary values. In general, enforcing non-homogeneous mixed boundary conditions
    for PH systems requires additional care \cite{brugnoli_2020_partitioned, thoma_2022_explicit,brugnoli_2022_explicit,thoma_2022_porthamiltonian}.
\end{remark}

\begin{remark} \label{rem_angmom}
    Notably the presented formulation also satisfies the fundamental balance of total angular momentum, which can be shown by means of the weak form \eqref{weak_form_directors}, as proven in Appendix~\ref{sec_ang_mom}.
\end{remark}

\subsection{Objectivity} \label{sec_obj}
Let us discuss how the presented formulation, with weak form \eqref{weak_form_directors}, behaves under superimposed rigid body motions, or, for that matter equivalently, a change of observer. To this end, consider $\vec{c}(t) \in \mathbb{R}^3$ and $\vec{\Lambda}(t) \in SO(3)$ that prescribe superimposed finite translations and rotations of the beam. Thus, the kinematic quantities change as
\begin{equation}
    \centerline^\# = \vec{c} + \vec{\Lambda} \centerline, \qquad \vec{d}_i^\# = \vec{\Lambda} \vec{d}_i, \qquad \vphi^\# = \dot{\vec{c}} + \vec{\Lambda} \vphi + \dot{\vec{\Lambda}}\centerline , \qquad \vec{v}_{d,i}^\# = \vec{\Lambda}  \vec{v}_{d,i} + \dot{\vec{\Lambda}}  \vec{d}_i
\end{equation}
for $i \in \{1, 2, 3\}$, Here, $\square^\#$ denotes a transformed quantity taking into account the superimposed rigid body motion.
The external distributed loads transform as
\begin{equation}
    \bar{\vec{n}}^\# = \vec{\Lambda} \bar{\vec{n}}, \qquad \bar{\vec{m}}^\# = \vec{\Lambda} \bar{\vec{m}} .
\end{equation}
The material strains \eqref{eq:strain_measures} remain invariant under that transformation, see also \cite[p.~1781]{betsch_2002_frameindifferent} and \cite[p.~1689]{romero_2002_objective}. Correspondingly, it can be stated that
\begin{equation}
    \vec{\Gamma}(\vec{d}^\#,\partial_s\centerline^\#) = \vec{\Gamma}(\vec{d},\partial_s \centerline) , \qquad \vec{K}(\vec{d}^\#,\partial_s\vec{d}^\#) = \vec{K}(\vec{d},\partial_s\vec{d}) .
\end{equation}
Using $\dot{\vec{\Lambda}}\transp\vec{\Lambda} + \vec{\Lambda}\transp\dot{\vec{\Lambda}} = \vec{0}$, it can also be shown that the corresponding strain rates $\dot{\vec{\Gamma}}$ and $\dot{\vec{K}}$ are unaltered such that the compliance equations in time-differentiated form \eqref{eq:strain_evo_stress_resultants} are identically satisfied by transformed quantities. This ensures the objectivity of the stress measures and their time-derivatives.
Further, the orthonormality constraints on velocity level \eqref{eq:constraint_director} are not affected by such superimposed rigid body motions.
Lastly, pre-multiplying the kinematic relation \eqref{eq:kinematics} and the balance equations \eqref{eom_rewritten} with $\vec{\Lambda}$ and $\diag(\vec{\Lambda},\vec{\Lambda},\vec{\Lambda})$, respectively, yields again the same relations expressed in the transformed quantities.

This proof may alternatively be done by means of the weak form \eqref{weak_form_directors} with the appropriately transformed test functions. Consequently, the present mixed formulation for PH geometrically exact beams is objective.
Notably, this property is not altered by the discretization methods proposed in the subsequent section, cf. Remark~\ref{rem_obj_discrete}.

\section{Structure-preserving discretization}
\label{sec_discretization}

\subsection{Mixed finite element discretization in space}
We pursue a structure-preserving discretization in space using a mixed finite element (FE) method  as for the 2D case in \cite{kinon_2025_energymomentumconsistent}. The PH system for the Cosserat rod dynamics \eqref{compliance_form} is readily available for spatial discretization via its weak form \eqref{weak_form_directors}.

\paragraph{Discretization ansatz.}
We discretize the domain into $\nel$ finite elements, i.e., $\domainspace = \cup_{e=1}^{\nel} \domainspace\e$, according to general standard procedures \cite{hughes2012finite}. Depending on the choice of ansatz functions this corresponds to $\nnode$ nodes on the FE grid. Displacement and velocity quantities $(\vec{q}, \vec{v})$ are approximated with $C^0$-continuous, Lagrangian polynomial ansatz functions of order $p$, while stress type quantities $\vec{\sigma}$ use discontinuous Lagrangian polynomial ansatz functions of order $p-1$.
The corresponding spaces are denoted
\begin{equation}
    \begin{aligned}
        \mathcal{V}\h   & = \big\{ \vec{q}\h \in H^1(\domainspace, \R^{12}) \, \big| \, \vec{q}\h|_{\domainspace\e} \in \mathcal{P}^p (\domainspace\e, \R^{12}) \, \forall \, \domainspace\e \in \domainspace, \vec{q}\h|_{\partial \domainspace_{\mathrm{D}}} = \vec{q}_{\mathrm{D}} \big\}      , \\
        \mathcal{V}\h_0 & = \big\{ \vec{q}\h \in H^1(\domainspace, \R^{12}) \, \big| \, \vec{q}\h |_{\domainspace\e} \in \mathcal{P}^p (\domainspace\e, \R^{12}) \, \forall \, \domainspace\e \in \domainspace, \vec{q}\h|_{\partial \domainspace_{\mathrm{D}}} = \vec{0} \big\}      ,             \\
        \mathcal{S}\h   & = \big\{ \vec{\sigma}\h \in L^2(\domainspace, \R^{6}) \ \, \big| \, \vec{\sigma}\h|_{\domainspace\e} \in \mathcal{P}^{p-1}(\domainspace\e, \R^{6}) \, \forall \, \domainspace\e \in \domainspace \big\} ,
    \end{aligned}
\end{equation}
for $p \in \{1 ,2\}$, where $\square\h$ denotes an approximation. Moreover, $H^1$ denotes the Sobolev space of square-integrable $C^0$-continuous functions, $L^2$ the space of square-integrable functions and $\mathcal{P}^p$ is the space of polynomial functions of order $p$. By this choice, we have $\nnode = \nel \, p + 1$.

The Lagrange multipliers are evaluated collocation-wise, i.e., with discrete values at the nodes and zero elsewhere, their approximated space denoted as $\mathcal{L}\h$. This can be formally realized via delta functions \cite[Chap.~1.11]{arfken_2013_mathematical}.
Correspondingly,
the orthonormality of the directors \eqref{eq:constraint_director}
is relaxed to nodal points of the FE mesh. This facilitates the use of Lagrange polynomials for $(\vec{q}, \vec{v})$ and is a common approach in the literature
\cite{betsch_2002_frameindifferent,betsch_2002_dae,romero_2002_objective}.\footnote{Alternatively, a weak enforcement can be employed, see \cite{wasmer_2024_projectionbased,harsch_2021_finite}.}
Concerning the interpolation between nodes,
\lq\lq \emph{the lack of orthonormality of the director frame can be regarded as a discretization error, which
    diminishes if the number of elements is increased}\rq\rq\, \cite{eugster_2014_directorbased}.
As employed later, an interpolation using quadratic ansatz functions yields practically indistinguishable interpolations compared to an $SE(3)$-finite element accounting for the correct manifold, see Fig.~2 in \cite{eugster_2023_family}.

The discrete weak formulation using a Bubnov-Galerkin approach is then stated as follows. Find $(\vec{q}\h,\vec{v}\h,\vec{\sigma}\h, \vec{\lambda}\h) \in \mathcal{V}\h \times \mathcal{V}\h \times \mathcal{S}\h \times \mathcal{L}\h$, such that \eqref{weak_form_directors} holds for all $(\delta\vec{q}\h, \delta\vec{v}\h,\delta\vec{\sigma}\h, \delta\vec{\lambda}\h) \in \mathcal{V}_0\h \times \mathcal{V}_0\h \times \mathcal{S}\h \times \mathcal{L}\h$.
To this end, we formulate the approximative quantities by means of ansatz matrices $\vec{\Phi},\vec{\Psi},\vec{\Xi}$ as
\begin{equation} \label{eq_ansatz}
    \begin{aligned}
        \vec{q}\h(s,t)      & = \vec{\Phi}(s) \hat{\vec{q}}(t) \in \mathcal{V}\h ,      & \quad  \vec{v}\h(s,t)       & = \vec{\Phi}(s) \hat{\vec{v}}(t) \in \mathcal{V}\h ,      \\
        \vec{\sigma}\h(s,t) & = \vec{\Psi}(s) \hat{\vec{\sigma}}(t) \in \mathcal{S}\h , & \quad  \vec{\lambda}\h(s,t) & = \vec{\Xi}(s) \hat{\vec{\lambda}}(t) \in \mathcal{L}\h ,
    \end{aligned}
\end{equation}
where $\hat{\vec{q}}, \hat{\vec{v}} \in \mathbb{R}^{12(\nel \, p + 1)}$, $\hat{\vec{\sigma}} \in \mathbb{R}^{6 \nel p}$ and $\hat{\vec{\lambda}} \in \mathbb{R}^{6 (\nel \, p + 1)}$ contain the unknowns from all finite elements and the respective nodes. Moreover, for the approximation of the corresponding test function $\delta\square\h$ consider the same basis as for the ansatz $\square\h$. In this context, the above ansatz matrices contain the polynomial ansatz functions, e.g., $  \vec{q}\h(s,t) = \vec{\Phi}(s) \hat{\vec{q}}(t) = \sum_{i=1}^{\nnode} N_i(s) \nodal{q}_i(t)$.
Therein, $N_i$ are basis functions for $\mathcal{P}^p(\domainspace)$.
For further details see Appendix A.2 in \cite{kinon_2025_energymomentumconsistent} for the planar case.
We further follow common FE procedures, see for instance \cite{hughes2012finite}, which also includes an assembly of the global system matrices and vectors accounting for the respective degrees of freedom. Additionally, we perform all computations on element level and consider the isoparametric concept to transform all calculations to a standard reference element with $\xi \in \hat{\domainspace} = [-1,1]$.

\paragraph{Finite-dimensional PH system.}
The discrete system dynamics are obtained by approximating the weak form \eqref{weak_form_directors} with the ansatz spaces \eqref{eq_ansatz}.
Taking into account the arbitrariness of the test function nodal values, we obtain the PH differential-algebraic equations (PH-DAE)
\begin{align}
    \label{compliance_form_discrete}
    \begin{bmatrix}
        \vec{I} & \zeros      & \zeros    & \zeros  \\
        \zeros  & \nodal{M}_o & \zeros    & \zeros  \\
        \zeros  & \zeros      & \nodal{C} & \zeros  \\
        \zeros  & \zeros      & \zeros    & \vec{0}
    \end{bmatrix}
    \begin{bmatrix}
        \dot{\nodal{q}}      \\
        \dot{\nodal{v}}      \\
        \dot{\nodal{\sigma}} \\
        \dot{\nodal{\lambda}}
    \end{bmatrix} & =
    \begin{bmatrix}
        \zeros    & \vec{I}                               & \zeros                                        & \zeros                       \\
        - \vec{I} & \zeros                                & -\nodal{J}_{\vec{\sigma v}}(\nodal{q})\transp & -\nodal{G}(\nodal{q})\transp \\
        \zeros    & \nodal{J}_{\vec{\sigma v}}(\nodal{q}) & \zeros                                        & \zeros                       \\
        \zeros    & \nodal{G}(\nodal{q})                  & \zeros                                        & \zeros
    \end{bmatrix}
    \begin{bmatrix}
        \vec{0}        \\
        \nodal{v}      \\
        \nodal{\sigma} \\
        \nodal{\lambda}
    \end{bmatrix}
    +
    \begin{bmatrix}
        \zeros                   & \zeros               \\
        \nodal{B}_{\domainspace} & \nodal{B}_{\partial} \\
        \zeros                   & \zeros               \\
        \zeros                   & \zeros
    \end{bmatrix}
    \begin{bmatrix}
        \nodal{u}_\Omega \\
        \vec{u}_\partial
    \end{bmatrix} .
\end{align}
The discrete system matrices can be found in Appendix~\ref{sec_appendix_system_matrices}.
Analogously to the infinite-dimensional version \eqref{phDAE_compact}, we can recast the semi-discrete set of equations in a PH framework \cite{mehrmann_2023_control,mehrmann_2019_structurepreserving,beattie_2018_linear} given by
\begin{equation} \label{discrete_ph}
    \begin{aligned}
        \nodal{E} \dot{\nodal{\state}} & = \nodal{J}(\nodal{\state}) \nodal{z}(\nodal{\state}) + \nodal{B}(\nodal{x})\nodal{u} , \\
        \nodal{y}                      & = \nodal{B}(\nodal{x})\transp \nodal{z}(\nodal{\state}) ,
    \end{aligned}
\end{equation}
where we have appended the discrete output equations for $ \nodal{y} = (\nodal{y}_\domainspace,\nodal{y}_\partial)$,
containing both distributed nodal as well as boundary terms. These quantities have the physical interpretation of velocities, that are power-conjugated to the nodal distributed and boundary input forces and torques. See Appendix~\ref{sec_appendix_system_matrices} for more details.
Here, the finite-dimensional version of the state \eqref{eq_state} contains the degrees of freedom of all finite elements and is given by
\begin{equation}
    \hat{\state} = (\hat{\vec{q}}, \hat{\vec{v}}, \hat{\vec{\sigma}}, \hat{\vec{\lambda}})  .
\end{equation}
Inserting the ansatz \eqref{eq_ansatz} into
\eqref{hamiltonian_director_stress} gives rise to the approximated Hamiltonian
\begin{align} \label{discrete_hamiltonian}
    H(\state\h) = \hat{H}(\nodal{\state}) =  \frac{1}{2} \nodal{v}\transp \nodal{M}_o \nodal{v} + \frac{1}{2} \nodal{\sigma}\transp \nodal{C}\nodal{\sigma} = \frac{1}{2} \nodal{\state}\transp \nodal{Q} \nodal{E} \nodal{\state} ,
\end{align}
where $\nodal{Q}= \diag(\vec{0},\vec{I},\vec{I},\vec{I})$, $\nodal{E} = \diag(\vec{I}, \nodal{M}_o, \nodal{C}, \vec{0})$ and further system matrices can be found in Appendix~\ref{sec_appendix_system_matrices}.

\paragraph{PH structure and power balance.}
We can again check the PH structure of the equations at hand.
Firstly, corresponding to condition (i) from the infinite-dimensional setting, it can be stated that
\begin{align} \label{discrete_const_law}
    \nodal{E}\transp \nodal{z}(\nodal{\state}) = \nabla \hat{H}(\nodal{\state})
\end{align}
holds true for the semi-discrete system because of the discrete Hamiltonian \eqref{discrete_hamiltonian} and the definitions of $\nodal{E}$ and $\nodal{z}(\nodal{x})=\nodal{Q}\nodal{x}$ in the PH-DAE \eqref{compliance_form_discrete}.
Secondly, the formally skew-adjoint nature of the structure operator $\operator{J}$ has been retained as skew-symmetry of the structure matrix $\nodal{J}(\nodal{\state}) = - \nodal{J}(\nodal{\state})\transp$ in \eqref{compliance_form_discrete}. Note that the state-dependence of $\nodal{J}$ highlights the geometric nonlinearity of the problem at hand. This corresponds to condition (ii) from the infinite-dimensional setting.

As already discussed, one core property of the PH formulation is that it satisfies a formalized power balance equation, see \eqref{power_balance_cont}. Since \eqref{discrete_ph} is a discrete PH-DAE system, it satisfies the semi-discrete power balance equation
\begin{equation} \label{energy_balance_discrete_space}
    \dot{\hat{H}} = \nabla \hat{H}(\nodal{x})\transp \dot{\nodal{x}}
    = \nodal{z}(\nodal{\state})\transp \nodal{E} \dot{\nodal{x}} = \nodal{z}(\nodal{\state})\transp \left( \nodal{J}(\nodal{\state}) \nodal{z}(\nodal{\state}) + \nodal{B}(\nodal{x})\nodal{u}\right)
    = \nodal{y}\transp\nodal{u} ,
\end{equation}
where we have made use of \eqref{discrete_ph}, \eqref{discrete_const_law} and the skew-symmetry of $\nodal{J}(\nodal{\state})$.
The equation resembles a discrete version of \eqref{power_balance_cont}, showcasing that the property of passivity and losslessness has been preserved within the discretization procedure.

\begin{remark} \label{rem_obj_discrete}
    We have outlined how to prove the objectivity of the present beam formulation in Section~\ref{sec_obj}. The same property is retained by the present space discretization and the proof is conducted in a completely analogous manner as before. To this end, keep in mind that spatial approximations are mere linear combinations of the nodal unknowns, e.g.,
    \begin{equation}
        (\centerline\h)^\# = \sum_{i=1}^{n_{\text{n}}} N_i(s) \nodal{\varphi}_i^\#(t),
    \end{equation}
    such that discretization and superimposition of rigid body motions commute. Consequently, the objectivity is inherited from the pure displacement-based formulations \cite{betsch_2002_frameindifferent,romero_2002_objective}, see also numerical verifications in the cited references.
\end{remark}

\subsection{Energy-momentum consistent discretization in time}

We further aim at a structure-preserving time discretization of the PH-DAE \eqref{discrete_ph}. To this end, we apply the implicit mid-point rule. Let $\nodal{\state}\n \approx \nodal{\state}(t\n)$ at time instance $t\n$ and consider an equidistant time-grid such that $[0, \finaltime] = \cup_{n=0}^N [ t\n, t\npe]$ with $N \geq 1$ time steps of constant size $h = t\npe-t\n > 0$.
The time-stepping scheme we propose reads
\begin{equation}
    \label{timediscrete_EoM}
    \begin{array}{rcl}
        \nodal{E}\left(\hat{\vec{x}}\npe-\hat{\vec{x}}\n\right) & = & h \, \nodal{J}(\hat{\vec{x}}\npeh) \hat{\vec{z}}(\nodal{x}\npeh) + h \, \nodal{B}(\hat{\vec{x}}\npeh)\hat{\vec{u}}\npeh , \\
        \hat{\vec{y}}\npeh                                      & = & \nodal{B}(\hat{\vec{x}}\npeh)\transp \hat{\vec{z}}(\nodal{x}\npeh) ,
    \end{array}
\end{equation}
for $n=0,\ldots,N-1$,
where the midpoint state $\hat{\vec{x}}\npeh:=\frac{1}{2}(\hat{\vec{x}}\npe+\hat{\vec{x}}\n)$ and $\hat{\vec{u}}\npeh = \hat{\vec{u}}(t\npeh)$.
Correspondingly, $\hat{\vec{y}}\npeh \approx \hat{\vec{y}}(t\npeh)$.
Additionally, the time-discrete costate vector $\hat{\vec{z}}(\nodal{x}\npeh)$ verifies
\begin{equation} \label{eq_co_state}
    \nodal{E}\transp\hat{\vec{z}}(\nodal{x}\npeh) = \nabla \hat{H}(\hat{\vec{x}}\npeh) ,
\end{equation}
since this is a mere pointwise evaluation of \eqref{discrete_const_law}.

The structure-preserving scheme is energy-consistent and second-order accurate. Let us demonstrate that an exact energy balance in discrete time is obtained.
For at most quadratic Hamiltonians, like here \eqref{discrete_hamiltonian}, the midpoint-evaluated gradient satisfies the directionality property\footnote{This property can also be exploited for more general Hamiltonians when using discrete gradient methods \cite{hairer_2006_geometric,kinon_2023_porthamiltonian,kinon_2023_discrete}.}, which mimics the chain rule in the discrete time setting such that
$\nabla\hat{H}(\hat{\vec{x}}\npeh)\transp\left(\hat{\vec{x}}\npe-\hat{\vec{x}}\n\right) = \hat{H}(\hat{\vec{x}}\npe) - \hat{H}(\hat{\vec{x}}\n)$, see also Remark~\ref{rem_meanvalue}.
Exploiting this property we obtain
\begin{equation} \label{energy_balance_discrete_time}
    \begin{aligned}
        \hat{H}(\hat{\vec{x}}\npe) - \hat{H}(\hat{\vec{x}}\n)
         & = \nabla\hat{H}(\hat{\vec{x}}\npeh)\transp\left(\hat{\vec{x}}\npe-\hat{\vec{x}}\n\right)                                                                        \\
         & = \nodal{z}(\nodal{\state}\npeh)\transp \nodal{E} \left(\hat{\vec{x}}\npe-\hat{\vec{x}}\n\right)                                                                \\
         & = h \nodal{z}(\nodal{\state}\npeh)\transp \left( \nodal{J}(\nodal{\state}\npeh) \nodal{z}(\nodal{\state}\npeh) + \nodal{B}(\nodal{x}\npeh)\nodal{u}\npeh\right) \\
         & = h\, (\nodal{y}\npeh) \transp \nodal{u}\npeh ,
    \end{aligned}
\end{equation}
where we have used \eqref{timediscrete_EoM}, \eqref{eq_co_state} and the skew-symmetry of $\nodal{J}$, analogously to the continuous-time derivation \eqref{energy_balance_discrete_space}. This proves that the present time-stepping scheme exhibits passivity and losslessness (which includes energy-conservation in the case of vanishing inputs).

\begin{remark}
    The presented time discretization also respects the balance of total angular momentum, see Appendix~\ref{sec_ang_mom}.
\end{remark}

\begin{remark}\label{rem_meanvalue}
    The above directionality condition can be traced back to the mean value theorem. In fact, given an at most quadratic polynomial function ${f} \in C^1(\R^d, \R)$,
    $
        {f}(\vec{b}) - {f}(\vec{a})
        = \nabla f ( \frac{\vec{a} + \vec{b}}{2})\transp(\vec{b}-\vec{a}) ,
    $
    for any $\vec{a},\vec{b}\in \R^d$. Likewise, this property extends to bilinear, vector-valued functions of two arguments $\vec{f} \in C^1(\R^d \times \R^d, \R^l)$ such that
    $
        \vec{f}(\vec{b},\vec{d}) - \vec{f}(\vec{a},\vec{c})
        = \vec{f}(\frac{\vec{a} + \vec{b}}{2}, \vec{d}-\vec{c}) +
        \vec{f}( \vec{b} - \vec{a}, \frac{\vec{c} + \vec{d}}{2}) ,
    $
    for any $\vec{a},\vec{b},\vec{c},\vec{d}\in \R^d$.
\end{remark}

\subsection{Kinematic consistency in discrete time} \label{sec_kinematic_consistency}
This section is concerned with two aspects regarding the approximated kinematics of the proposed model.
Firstly,
the orthonormality constraint \eqref{eq:balance_ang_momentum_2} is satisfied by the proposed time integration approach despite using the velocity-level constraint \eqref{eq:constraint_director}.
To show this, consider the first line of \eqref{compliance_form_discrete}, which is discretized in time by means of \eqref{timediscrete_EoM} such that we obtain the
discrete time kinematics for the directors
\begin{align}
    (\nodal{d}_k)\npe - (\nodal{d}_k)\n & = h (\nodal{v}_{d,k})\npeh \label{discrete_1} ,
\end{align}
for each FE node $k$. Note that $\nodal{d}_k \in \R^9$ refers to all three director unknowns at node $k$.
Furthermore, the fourth line of \eqref{compliance_form_discrete}, which is discretized in time by means of \eqref{timediscrete_EoM}, yields
\begin{align}
    \vec{G}_d((\nodal{d}_k)\npeh) (\nodal{v}_{d,k})\npeh & = \vec{0} , \label{discrete_2}
\end{align}
for each node $k$, see further details in relations~\eqref{sifting} and \eqref{sifting_2} in the Appendix.
We can now make use of Remark~\ref{rem_meanvalue}, since the orthonormality constraints \eqref{eq:ortho_constraint} are quadratic and \eqref{discrete_2} contains the constraint Jacobian $\vec{G}_d(\vec{d}) = \partial_{\vec{d}} \vec{g}(\vec{d})$.
Substituting \eqref{discrete_1} in \eqref{discrete_2} for the midpoint velocities, one thus obtains
\begin{align}
    \vec{0} =  \partial_{\vec{d}} \vec{g}((\nodal{d}_k)\npeh) ((\nodal{d}_k)\npe - (\nodal{d}_k)\n) = \vec{g}((\nodal{d}_k)\npe) - \vec{g}((\nodal{d}_k)\n)
\end{align}
for all $k$. Eventually, given a consistent initial choice of orthonormal directors $(\nodal{d}_k)_0$ the constraints
\begin{align} \label{eq_proof_no_drift}
    \vec{g}((\nodal{d}_k)\n) = \vec{0}
\end{align}
are exactly fulfilled
for all $k$, $n\geq1$. Thus, the orthonormality constraint is satisfied at each FE node in each discrete point in time. This is due to the structure-preserving midpoint time discretization of the velocity-level constraints \eqref{eq:constraint_director} pertaining to the present PH model.

\paragraph{Consistent strain approximation.}
Further, we can show that the discrete strains calculated from the independent stresses via \eqref{eq:material}
\begin{equation} \label{discrete_strains}
    \vec{\Gamma}\h\n := \vec{C}_N \vec{N}\h\n + \vec{\Gamma}_0, \qquad \qquad \vec{K}\h\n := \vec{C}_M \vec{M}\h\n + \vec{K}_0
\end{equation}
are variationally consistent with the discrete strains computed in terms of the displacement quantities \eqref{eq:strain_measures_short}, i.e., $\vec{\Gamma}(\centerline\h\n,\vec{d}\h\n)$ and $\vec{K}(\vec{d}\h\n)$.
Note our Remark~\ref{rem_initial_conf} on the assumed reference strains $\vec{\Gamma}_0$, $\vec{K}_0$.
This consistency ensures a physically valid approximation in each time instance despite the fact that the PH model is based on the evolution equations \eqref{eq:strain_evo_stress_resultants}. Hereto, we evaluate \eqref{eq:strain_measures_short} in the different time instances and use the findings from Remark~\ref{rem_meanvalue} due to the bilinearity of the strain functions. This yields
\begin{equation}\label{eq_discrete_strain_disp}
    \begin{aligned}
        \vec{\Gamma}(\centerline\h\npe,\vec{d}\h\npe) - \vec{\Gamma}(\centerline\h\n,\vec{d}\h\n) & =  \rotmat(\vec{d}\h\npeh)\transp \partial_s (\centerline\h\npe - \centerline\h\n)
        + \vec{J}_{\vec{N}}(\partial_s \centerline\h\npeh)\transp (\vec{d}\h\npe - \vec{d}\h\n) ,                                                                                      \\
        \vec{K}(\vec{d}\h\npe) - \vec{K}(\vec{d}\h\n)                                             & =
        \vec{L}( \partial_s \vec{d}\h\npeh)(\vec{d}\h\npe - \vec{d}\h\n)
        - \vec{L}(\vec{d}\h\npeh) \partial_s(\vec{d}\h\npe - \vec{d}\h\n) .
    \end{aligned}
\end{equation}
On the other hand, consider the third line of \eqref{compliance_form_discrete}, which is discretized in time by means of \eqref{timediscrete_EoM}.
The equivalent discrete space and time versions of the weak forms \eqref{weak_form_directors_5} and \eqref{weak_form_directors_6} are given by
\begin{equation} \label{eq_discrete_strain_stress}
    \begin{aligned}
        \inner{\delta \vec{N}\h}{\vec{\Gamma}\h\npe - \vec{\Gamma}\h\n -
            h \rotmat(\vec{d}\h\npeh)\transp \partial_s (\vphi\h)\npeh - h \vec{J}_{\vec{N}}(\vec{\centerline}\h\npeh)\transp (\vec{v}_d\h)\npeh
        } & =0 , \\
        \inner{\delta \vec{M}\h}{\vec{K}\h\npe - \vec{K}\h\n -
            h \operator{J}_{\vec{K}}(\vec{d}\h\npeh) (\vec{v}_d\h)\npeh
        } & =0 .
    \end{aligned}
\end{equation}
We may now use $\centerline\h\npe - \centerline\h\n = h (\vphi\h)\npeh$ and $\vec{d}\h\npe - \vec{d}\h\n = h (\vec{v}_d\h)\npeh$, which are equivalent to the time-discretized first line of \eqref{compliance_form_discrete}, and account for the definition \eqref{eq_BM_BK} of $\operator{J}_{\vec{K}}$. Eventually, taking into account \eqref{eq_discrete_strain_disp}, \eqref{eq_discrete_strain_stress} yields
\begin{equation} \label{eq_proof_consistent_strain}
    \begin{aligned}
        \inner{\delta \vec{N}\h}{\vec{\Gamma}\h\npe - \vec{\Gamma}\h\n - \left( \vec{\Gamma}(\centerline\h\npe,\vec{d}\h\npe) - \vec{\Gamma}(\centerline\h\n,\vec{d}\h\n) \right)
        } & =0 , \\
        \inner{\delta \vec{M}\h}{\vec{K}\h\npe - \vec{K}\h\n - \left( \vec{K}(\vec{d}\h\npe) - \vec{K}(\vec{d}\h\n) \right)
        } & =0 ,
    \end{aligned}
\end{equation}
for all $n\geq0$.
We now assume initial values for the stress variables
$\vec{N}_0^h$, $\vec{M}_0^h$ that,
through~\eqref{discrete_strains}, induce consistent strain measures satisfying
$\vec{\Gamma}\h_0 = \vec{\Gamma}(\centerline\h_0,\vec{d}\h_0)$ and $\vec{K}\h_0 = \vec{K}(\vec{d}\h_0)$. This
yields the variational identities
\begin{equation} \label{eq_initialization_strains}
    \begin{aligned}
        \inner{\delta \vec{N}\h}{\vec{\Gamma}\h_0 - \vec{\Gamma}(\centerline\h_0,\vec{d}\h_0)
        } & =0 , \\
        \inner{\delta \vec{M}\h}{\vec{K}\h_0 - \vec{K}(\vec{d}\h_0)
        } & =0 .
    \end{aligned}
\end{equation}
Eventually, \eqref{eq_proof_consistent_strain} yields the result
\begin{equation} \label{eq_proof_consistent_strain_new}
    \begin{aligned}
        \inner{\delta \vec{N}\h}{\vec{\Gamma}\h\npe - \vec{\Gamma}(\centerline\h\npe,\vec{d}\h\npe)
        } & =0 , \\
        \inner{\delta \vec{M}\h}{\vec{K}\h\npe - \vec{K}(\vec{d}\h\npe)
        } & =0 ,
    \end{aligned}
\end{equation}
for all $n\geq0$. Consequently, despite the original rate form of the underlying PH formulation, consistent values for $\vec{\Gamma}\h\npe$ and $\vec{K}\h\npe$ are recovered at discrete times $t\npe$. By \eqref{discrete_strains}, this implies variationally consistent mixed stress quantities $\vec{N}\h\npe$ and $\vec{M}\h\npe$.

\section{Extension of the model}
\label{sec_extension}
In practical applications, such as soft robotics, visco-elastic material behavior and actuation via pneumatic chambers and tendons are highly relevant.
To this end, we show in this section that the inclusion of such concepts into the PH formulation \eqref{compliance_form} is seamlessly possible.
These extensions do not interfere with the discretization procedure from the previous chapter.

\subsection{Extension to visco-elasticity} \label{sec_visco}

We target an extension to visco-elastic material behavior using the generalized-Maxwell model. It represents an \emph{internal variable model} and has been discussed in, e.g., \cite{ferri_2023_efficient,weeger_2022_mixed,lestringant_2020_discrete} in the context of Cosserat rods, or \cite{kinon_2024_generalized} for geometrically exact strings. While we apply this approach to each strain measure separately, some works only include damping for the bending behaviour \cite{lestringant_2020_discrete}.

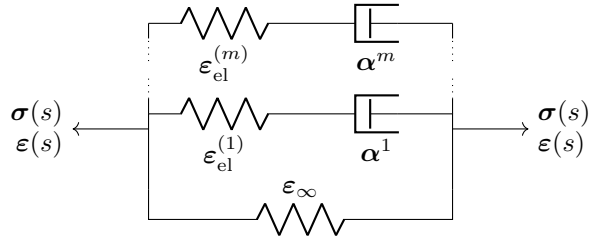
\begin{figure}[tb]
    \centering
    \begin{circuitikz}
    \draw[<-]  (0,0.8) node[left,text width=1.3 cm, align=right] {$\vec{\sigma}(s)$ $\vec{\epsilon}(s)$} to (1,0.8);
    \draw (1,0) to (1,-0.4) to[R=$\vec{\epsilon}_\infty$] (5,-0.4) to (5,0);
    \draw[->] (5,0.8) to (6,0.8) node[right,text width=1.3 cm] {$\vec{\sigma}(s)$ $ \vec{\epsilon}(s)$};
    \draw (5,0) to (5,1) to [damper=$\vec{\alpha}^1$] (3,1) to[R=$ \vec{\epsilon}_\mathrm{el}^{(1)}$] (1,1) to (1,0);
    \draw (1,1) to (1,1.2);
    \draw[dotted] (1,1.2) to (1,1.5);
    \draw (5,1) to (5,1.2);
    \draw[dotted] (5,1.2) to (5,1.5);
    \draw[dotted] (1,1.7) to (1,2);
    \draw[dotted] (5,1.7) to (5,2);
    \draw (1,2) to (1,2.2);
    \draw (5,2) to (5,2.2);
    \draw (5,2.2) to [damper=$\vec{\alpha}^m$] (3,2.2) to[R=$ \vec{\epsilon}_\mathrm{el}^{(m)}$] (1,2.2);
\end{circuitikz}
    \caption{Schematic rheological model for the generalized-Maxwell approach.
    }
    \label{fig:gen_max}
\end{figure}

The model consists of an elastic element and $m \geq 1$ viscous Maxwell elements in parallel, see the rheological model in Fig.~\ref{fig:gen_max}. It is also referred to as \emph{generalized relaxation model}, see \cite[Chap.~10]{simo_2006_computational}. In this context, each additional branch introduces new material parameters, allowing for a more precise modeling of the time-dependent material behavior compared to, for example, the common Kelvin-Voigt model \cite{linn_2013_geometrically}. For the sake of the deductions in this section, consider $\vec{\epsilon} = (\vec{\Gamma},\vec{K})$ as well as the definition of the stress quantities $\vec{\sigma}$ and the compliance matrix $\vec{C}$, see \eqref{eq_abbrev}.

%
Due to the structure of the model,
the stress quantities of all branches add up to the total stress,
\begin{align} \label{eq:total_stress}
    \vec{\sigma} = \vec{\sigma}_\infty + \sum_{i=1}^{m} \vec{\sigma}_\mathrm{el}^{(i)} ,
\end{align}
where $\vec{\sigma}_\infty $ is the long-term elastic stress and $\vec{\sigma}_\mathrm{el}^{(i)}$ are the individual stress contributions of each branch. Note that in each viscous branch the elastic stress is equal to the viscous stress in the Maxwell element $\vec{\sigma}_\mathrm{el}^{(i)} = \vec{\sigma}_\mathrm{visc}^{(i)}$.
Moreover, the total strain is equal to the strain in each branch, i.e.,
\begin{align} \label{eq:total_strain}
    \vec{\epsilon}(\centerline,\vec{d}) = \vec{\epsilon}_\infty = \vec{\epsilon}_\mathrm{el}^{(i)} + \vec{\alpha}^{(i)} ,
\end{align}
for each $i$, where $\vec{\alpha}^{(i)}$ is the viscous strain in the damper element of the respective branch and $\vec{\epsilon}_\infty$ is the strain in the purely elastic branch. A common assumption is that the viscous stress
\begin{align} \label{eq:viscous_stress}
    \vec{\sigma}_\mathrm{visc}^{(i)} = \vec{V}^{(i)} \dot{\vec{\alpha}}^{(i)}
\end{align}
is linearly related to the viscous strain rate via a symmetric, positive semi-definite viscosity matrix $\vec{V}^{(i)}$.

To account for this linearly visco-elastic behavior in the PH model, we first introduce additional state variables $\vec{\sigma}_\mathrm{el}^{(i)} = \vec{\sigma}_\mathrm{visc}^{(i)}$ for each branch. Next, we include additional stress contributions in the balance of momentum \eqref{compliance_form} according to \eqref{eq:total_stress}. Lastly, we account for evolution equations
\begin{equation}
    \begin{aligned}
        \vec{C}_\infty \dot{\vec{\sigma}}_\infty           & = \dot{\vec{\epsilon}}(\centerline,\vec{d}) ,                                                                       \\
        \vec{C}^{(i)} \dot{\vec{\sigma}}_\mathrm{el}^{(i)} & = \dot{\vec{\epsilon}}(\centerline,\vec{d}) - \vec{V}^{(i)-1} \vec{\sigma}_\mathrm{el}^{(i)} , \quad i=1,\ldots,m ,
    \end{aligned}
\end{equation}
for the long term elastic stress $\vec{\sigma}_\infty$ and for each individual elastic stress contribution $\vec{\sigma}_\mathrm{el}^{(i)}$ in a compliance form based on the time derivative of the two equalities in \eqref{eq:total_strain}.
We have further accounted for the constitutive relation \eqref{eq:viscous_stress} as well as $\vec{C}_\infty \vec{\sigma}_\infty = \vec{\epsilon}_\infty$ and $\vec{C}^{(i)} \vec{\sigma}_\mathrm{el}^{(i)} = \vec{\epsilon}_\mathrm{el}^{(i)}$. Note that for consistency with a purely elastic beam model with elasticity matrices $\vec{C}_\Gamma$ and $\vec{C}_K$, see Eq.~\ref{quadratic_stored_energy}, the elasticity matrices of all branches must add up to those matrices.

\begin{remark} \label{rem_V}
    A common approach \cite{linn_2013_geometrically,ferri_2023_efficient} is to assume diagonal viscosity matrices $\vec{V}^{(i)}$ based on the same structure as the elasticity matrices. To this end, one can introduce dynamic viscosities $\eta_G, \eta_E \geq 0$ for shear and torsion as well as dilatation and bending, respectively, such that
    \begin{align}
        \tau_G = \frac{\eta_G}{G}, \qquad  \qquad \tau_E = \frac{\eta_E}{E}
    \end{align}
    are the corresponding characteristic relaxation times. Here, $E, G \geq 0$ denote Young's modulus and the shear modulus, respectively. This eventually gives rise to the viscosity matrices
    \begin{align}
        \vec{V}^{(i)} = \begin{bmatrix}
                            \vec{V}_\Gamma^{(i)} & \zeros          \\
                            \zeros               & \vec{V}_K^{(i)}
                        \end{bmatrix}
        \quad \text{with} \quad
        \vec{V}_\Gamma^{(i)} = \begin{bmatrix}
                                   \tau_G^{(i)} & \zero        & \zero        \\
                                   \zero        & \tau_G^{(i)} & \zero        \\
                                   \zero        & \zero        & \tau_E^{(i)}
                               \end{bmatrix} \vec{C}_\Gamma^{(i)} , \quad \qquad  \quad
        \vec{V}_K^{(i)}      = \begin{bmatrix}
                                   \tau_E^{(i)} & \zero        & \zero        \\
                                   \zero        & \tau_E^{(i)} & \zero        \\
                                   \zero        & \zero        & \tau_G^{(i)}
                               \end{bmatrix} \vec{C}_K^{(i)} ,
    \end{align}
    where $\vec{C}_\Gamma^{(i)}$ and $\vec{C}_K^{(i)}$ are elasticity matrices consistent with the compliance matrix $\vec{C}^{(i)}$ of the $i$-th viscous branch.
\end{remark}

\begin{remark}
    The Kelvin-Voigt dissipation model, as for instance used in \cite{linn_2013_geometrically,herrmann_2024_relativekinematic}, is contained in the present framework if $m=1$ dissipative Maxwell branch is considered. In this branch, the compliance of the elastic element has to be zero, such that it behaves rigidly.
\end{remark}

\paragraph{Dissipative PH model}
The dissipative system can again be viewed from a PH perspective. To this end, we concatenate all stresses of the Maxwell branches in $\vec{\sigma}_\mathrm{vis} = (\vec{\sigma}_\mathrm{el}^{(1)}, \dots, \vec{\sigma}_\mathrm{el}^{(m)})$ to define the state and costate of the new PH system as
\begin{align}
    \state = (\centerline, \vec{v}, \vec{\sigma}_\infty, \vec{\sigma}_\mathrm{vis}, \vec{\lambda}) , \qquad \qquad
    \vec{z}(\state) = (\vec{0}, \vec{v}, \vec{\sigma}_\infty, \vec{\sigma}_\mathrm{vis}, \vec{\lambda}) .
\end{align}
together with the Hamiltonian
\begin{align} \label{hamiltonian_director_stress_diss}
    H(\vec{x}) & = \int_\domainspace \left( \frac{1}{2} \vec{v}\transp \vec{M}_o \vec{v} + \frac{1}{2} \vec{\sigma}_\infty\transp \vec{C}_\infty \vec{\sigma}_\infty + \frac{1}{2}  \vec{\sigma}_\mathrm{vis} \transp \vec{C}_\mathrm{vis}  \vec{\sigma}_\mathrm{vis} \right) \d s ,
\end{align}
comprising the kinetic energy as of before together with all complementary strain energy terms pertaining to the elastic elements of the generalized-Maxwell model with $\vec{C}_\mathrm{vis} = \diag(\vec{C}^{(1)}, \dots,  \vec{C}^{(m)})$.
Similarly, we define $ \vec{V}\inv = \diag((\vec{V}^{(1)})\inv, \dots, (\vec{V}^{(m)})\inv)$,
$ \operator{J}_{\vec{v}}^\mathrm{vis} = \left[\operator{J}_{\vec{v}} \dots  \operator{J}_{\vec{v}} \right]$ and
$\operator{J}_{\vec{\sigma}}^\mathrm{vis} = \left[ \operator{J}_{\vec{\sigma}}\transp \dots  \operator{J}_{\vec{\sigma}}\transp \right]\transp $ with $m$ copies of the respective quantity. Thus, we write the dynamics as
\begin{align}
    \label{compliance_form_diss}
    \begin{bmatrix}
        \vec{I} & \zeros    & \zeros         & \zeros               & \zeros  \\
        \zeros  & \vec{M}_o & \zeros         & \zeros               & \zeros  \\
        \zeros  & \zeros    & \vec{C}_\infty & \zeros               & \zeros  \\
        \zeros  & \zeros    & \zeros         & \vec{C}_\mathrm{vis} & \zeros  \\
        \zeros  & \zeros    & \zeros         & \zeros               & \vec{0}
    \end{bmatrix}
    \begin{bmatrix}
        \dot{\vec{q}}                   \\
        \dot{\vec{v}}                   \\
        \dot{\vec{\sigma}}_\infty       \\
        \dot{\vec{\sigma}}_\mathrm{vis} \\
        \dot{\vec{\lambda}}
    \end{bmatrix} & =
    \begin{bmatrix}
        \zeros    & \vec{I}                                  & \zeros                  & \zeros                               & \zeros          \\
        - \vec{I} & \zeros                                   & -\operator{J}_{\vec{v}} & -\operator{J}_{\vec{v}}^\mathrm{vis} & -\vec{G}\transp \\
        \zeros    & \operator{J}_{\vec{\sigma}}              & \zeros                  & \zeros                               & \zeros          \\
        \zeros    & \operator{J}_{\vec{\sigma}}^\mathrm{vis} & \zeros                  & -\vec{V}\inv                         & \zeros          \\
        \zeros    & \vec{G}                                  & \zeros                  & \zeros                               & \zeros
    \end{bmatrix}
    \begin{bmatrix}
        \vec{0}                   \\
        \vec{v}                   \\
        \vec{\sigma}_\infty       \\
        \vec{\sigma}_\mathrm{vis} \\
        \vec{\lambda}
    \end{bmatrix}
    +
    \begin{bmatrix}
        \zeros                 \\
        \vec{B}_{\domainspace} \\
        \zeros                 \\
        \zeros                 \\
        \zeros
    \end{bmatrix}
    \vec{u}_\Omega .
\end{align}
In short, this fits again into the PH framework \cite{mehrmann_2019_structurepreserving,mehrmann_2023_control,beattie_2018_linear}
\begin{equation} \label{phDAE_compact_diss}
    \begin{aligned}
        \vec{E}\dot{\vec{x}} & = (\operator{J}(\vec{x}) - \operator{R}) \vec{z}(\vec{x}) + \vec{B}(\vec{x})\vec{u}_\Omega , \\
        \vec{y}_\domainspace & = \vec{B}(\vec{x})\transp \vec{z} ,
    \end{aligned}
\end{equation}
where $\operator{R}=\diag(\vec{0},\vec{0},\vec{0},\vec{V}\inv,\vec{0})$ is a constant, symmetric, positive semi-definite matrix\footnote{In general, this notion can be extended to a state-dependent, self-adjoint differential operator.}. The positive semi-definite property can be checked by means of the block-diagonal structure and $(\vec{V}^{(i)})\inv$ being symmetric and positive semi-definite, see Remark~\ref{rem_V}.

The power-balance equation induced by this PH system extends the one by the purely elastic model \eqref{power_balance_cont}
such that
\begin{align} \label{power_balance_cont_diss}
    \dot H & =  \vec{y}_{\partial \domainspace}\transp \vec{u}_{\partial \domainspace} + \int_\domainspace \vec{y}_\domainspace\transp \vec{u}_\domainspace \d s - D    .
\end{align}
Here, we have introduced the dissipated power due to viscous deformations (also called \emph{dissipation function} \cite{simo_2006_computational})
\begin{equation}
    D = \int_\domainspace \vec{z} \transp \operator{R} \vec{z} \d s =  \int_\domainspace (\vec{\sigma}_\mathrm{vis})\transp \vec{V}\inv \vec{\sigma}_\mathrm{vis} \d s = \int_\domainspace \sum_{i=1}^{m} (\vec{\sigma}_\mathrm{el}^{(i)})\transp (\vec{V}^{(i)})\inv \vec{\sigma}_\mathrm{el}^{(i)} \d s \geq 0 .
\end{equation}
Consequently,
the power balance equations \eqref{power_balance_cont_diss} immediately leads to the
dissipation inequality
\begin{equation}
    \dot H \leq \vec{y}_{\partial \domainspace}\transp \vec{u}_{\partial \domainspace} + \int_\domainspace \vec{y}_\domainspace\transp \vec{u}_\domainspace \d s ,
\end{equation}
proving the passivity of the dissipative infinite-dimensional PH system.

\begin{remark}
    As shown in \cite{kinon_2024_generalized} for geometrically exact strings, we can obtain the model \eqref{compliance_form_diss} alternatively by power-preserving interconnection. In this context, all viscous Maxwell branches are interpreted as purely dissipative PH models and the elastic branch as its own model similar to \eqref{compliance_form}. The overall model is obtained by means of a gyrator interconnection of passive subsystems.
\end{remark}

\subsection{Extension to pneumatic and tendon actuation} \label{sec_pneumatic_tendon}

In the context of soft robotic manipulators, two state-of-the-art actuation approaches are given by
elastic cables acting as so-called tendons \cite{rucker_2011_statics,till_2019_realtime} and pneumatic chambers \cite{eugster_2022_soft,till_2019_realtime}. Both can be included seamlessly in the present beam model in a unified way.
We follow the approach from \cite{eugster_2022_soft,till_2019_realtime}, see also \cite{alessi_2024_rod}, and include additional stress contributions acting as distributed actuation terms\footnote{Equivalently, one can account for these effects using the distributed forces and moments per unit length $(\bar{\vec{n}},\bar{\vec{m}})$ combined with related Neumann boundary terms. This approach is however more tedious to implement in an FE framework.}.

\begin{figure}[tb]
    \centering
    \def\svgwidth{0.7\textwidth}
    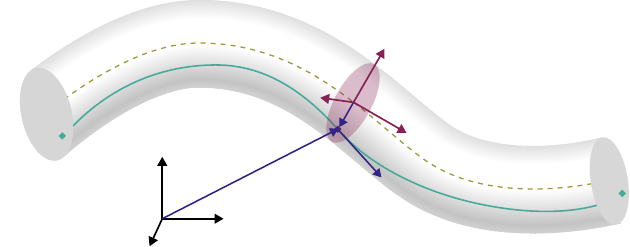
    \caption{Schematic depiction for the actuation of a Cosserat rod, where the $k$-th actuation element and related kinematic quantities are displayed.}
    \label{fig_actuation_sketch}
\end{figure}

Thus,
we extend the equations of motion \eqref{compliance_form} by introducing additional terms as
\begin{align} \label{eq_add_stresses_replacement}
    \vec{N} & \rightarrow \vec{N} + \rotmat(\vec{d})\transp \vec{n}_u , \qquad \qquad \vec{M} \rightarrow \vec{M} + \rotmat(\vec{d})\transp \vec{m}_u ,
\end{align}
where $\vec{n}_u, \vec{m}_u \in \mathbb{R}^3$ are additional internal forces and moments emerging from input actuation.
For both tendon and pressure chamber actuation they assume the form
\begin{equation} \label{eq_add_stresses}
    \vec{n}_u = \sum_{k=1}^{N} \tau_k \vec{t}_k , \qquad \qquad \vec{m}_u = \sum_{k=1}^{N} \vec{\rho}_k \times \tau_k \vec{t}_k ,
\end{equation}
where $k$ is the index of the actuation element, i.e., tendon or pressure chamber, and
\begin{equation}
    \tau_k = \begin{cases}
        p_k(t) A_k \leq 0 & \text{for pressure actuation} , \\
        T_k(t) \geq 0     & \text{for tendon actuation} ,
    \end{cases}
\end{equation}
the actuation force magnitude.
Therein, $p_k$ is the pressure and $A_k$ is the cross section area of the pressurized chamber, while $T_k$ is the tension force in a tendon. Both $p_k(t)$ and $T_k(t)$ are regarded as known input quantities, which will be related to a distributed forcing term.
Moreover in \eqref{eq_add_stresses}, $\vec{t}_k$ is the corresponding unit direction vector and
$\vec{\rho}_k(s,t) = \rho^\alpha_k \vec{d}_\alpha(s,t)$,
for $\alpha = 1,2$, describes the vector from the centerline of the beam to the line of centroids of the $k$-th actuation element. Here, $\rho^\alpha_k$ are material coordinates, which are fixed.
Correspondingly,
$\vec{r}_k = \centerline + \vec{\rho}_k $.
The only difference between the two actuation methods is that tendons induce traction in tangential direction of the tendons while pressure chambers induce traction in cross-section normal direction.
This is realized via the unit direction vector $\vec{t}_k$ such that
\begin{equation}\label{eq:unit-direction-vector}
    \vec{t}_k = \begin{cases}
        \vec{d}_3                                                      & \text{for pressure actuation} , \\
        \frac{\vec{r}_{k,s}}{\lvert\lvert \vec{r}_{k,s} \rvert \rvert} & \text{for tendon actuation} .
    \end{cases}
\end{equation}
For a schematic depiction for the case of tendon actuation see Fig.~\ref{fig_actuation_sketch}.
In the strong form \eqref{compliance_form}, these terms result in additional distributed inputs of the form
$\vec{B}_{u}(\state) \vec{\tau}$, where $\vec{\tau}=(\tau_1, \dots, \tau_N) \in \R^N$ are additional input quantities.
To account for \eqref{eom_rewritten}, \eqref{eq_add_stresses_replacement} and \eqref{eq_add_stresses}, we introduce the distributed port matrix
\begin{align}
    \vec{B}_{u}(\state) & = \begin{bmatrix}
                                (\vec{t}_{k,s})_{k=1}^{N} \\
                                -(\vec{J}_{\vec{N}}\vec{R}\transp \vec{t}_k + \operator{J}_{\vec{M}} (\vec{R}\transp  \tilde{\vec{\rho}}_k \vec{t}_k))_{k=1}^{N}
                            \end{bmatrix} \in \R^{6 \times N} ,
\end{align}
where each column of the matrix corresponds to the influence of one actuation element. Note that the state-dependency of the port matrix does not only emerge from the operators but also from the unit direction vector $\vec{t}_k$ introduced in \eqref{eq:unit-direction-vector}. Likewise, this also leads to an additional input term in the differential form of the PH dynamics \eqref{phDAE_compact}, eventually extended to
\begin{align}
    \vec{E}\dot{\vec{x}} = \operator{J}(\vec{x})\vec{z}(\vec{x}) + \vec{B}(\vec{x})\vec{u}_\Omega  + \vec{B}_{\vec{\tau}}(\vec{x})\vec{\tau} , \qquad
    \text{where}
    \qquad
    \vec{B}_{\vec{\tau}}(\vec{x}) & = \begin{bmatrix}
                                          \zeros              \\
                                          \vec{B}_{u}(\state) \\
                                          \zeros              \\
                                          \zeros
                                      \end{bmatrix} .
\end{align}
Similarly, these additional contributions are accounted for in the weak form analogously to the internal forces and moments via additional terms
$\inner{\delta \centerline_{,s}}{\vec{n}_u}$ and $\inner{\operator{J}_{\vec{K}}\delta \vec{d} }{\vec{R}\transp \vec{m}_u}$
in \eqref{weak_form_directors_3} and \eqref{weak_form_directors_4}, respectively.

\section{Numerical examples}
\label{sec_examples}
Let us apply the proposed framework to representative numerical examples.  Newton's method is used to solve the nonlinear discrete time equations \eqref{timediscrete_EoM} in every time step with an absolute tolerance $\epsilon$ on the norm of the residual. The spatial integrals are approximated by means of a three-point Gauss quadrature and for the polynomials ansatz functions we choose $p=2$. The Dirichlet boundary conditions for the displacements and velocities are incorporated directly into the
ansatz functions, see Remark~\ref{rem_BCs}, and no gravitational effects are considered.
The computations have been performed using the finite element code \texttt{moofeKIT}\footnote{\url{https://github.com/kit-ifm/moofeKIT}}, and the generated data are published in \cite{kinon_2025_repository} and archived under \cite{kinon_2025_18007058}.
The first example is used to assess the structure-preservation of the employed discretization approach. Additionally, good agreements with an analytical solution are shown.
In the second benchmark, we find good agreement of our results with data reported in the literature and verify kinematic consistency properties according to Sec.~\ref{sec_kinematic_consistency}. This example also deals with the visco-elastic extension from Sec.~\ref{sec_visco} and underlines the usability of the present framework for inextensible Kirchhoff beams.
The third example studies the formulation in a quasistatic scenario showing that the present formulation matches an analytical solution and avoids locking.
Lastly, we use our formulation to simulate the dynamic motion of a soft robotic arm actuated by means of pneumatic chambers, as introduced in Sec.~\ref{sec_pneumatic_tendon}.

\begin{table}[tb]
    \centering
    \caption{Flying spaghetti: Physical and numerical parameters.}
    \begin{tabular}{l l l l l l l l l l}
        \toprule
        $h$   & $\finaltime$ & $\nel$ & $\epsilon$ & $L$  & $\rho A$ & $M_\rho^{11} = M_\rho^{22}$ & $k_{\mathrm{b1}} = k_{\mathrm{b2}} = k_{\mathrm{t}}$ & $k_{\mathrm{s1}}=k_{\mathrm{s2}}=k_{\mathrm{e}}$ \\
        \midrule
        $0.1$ & $15$         & $10$   & $10^{-11}$ & $10$ & $ 1 $    & $ 10 $                      & $ 10^3 $                                             & $10^4$                                           \\
        \bottomrule
    \end{tabular}
    \label{tab_freeflight}
\end{table}%
\begin{figure}[t]
    \setlength{\plotwidth}{0.9\textwidth}
    \setlength{\plotheight}{0.5\plotwidth}
    \centering
    \tikzsetnextfilename{flyingspaghetti_snapshots}
    \input{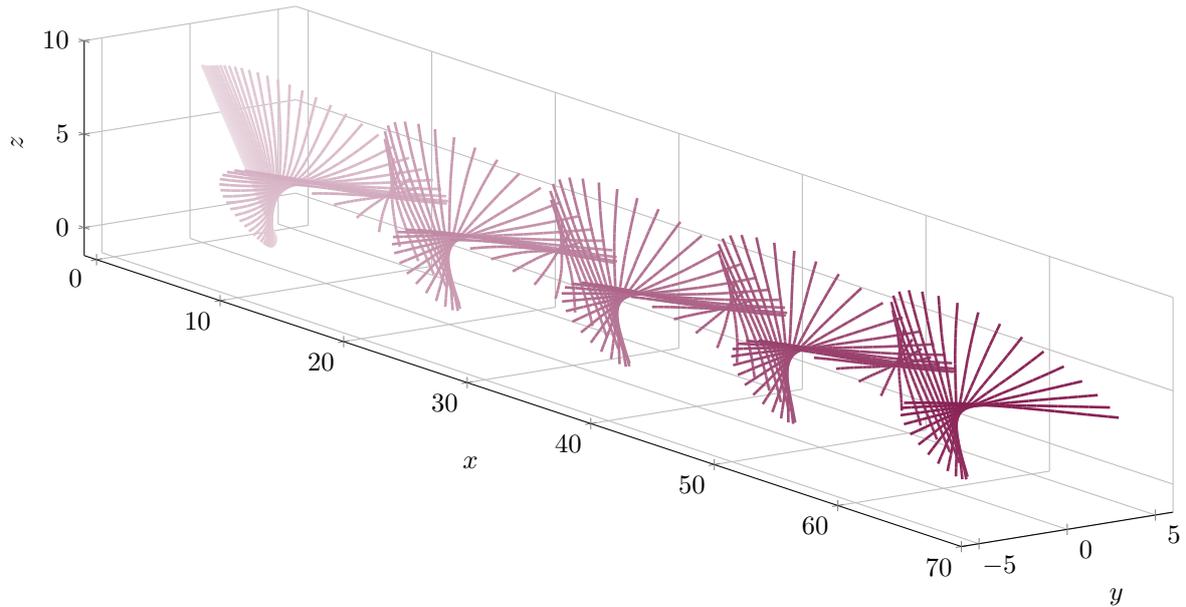}
    \caption{Flying spaghetti: Snapshots of the \textit{kayak-rowing} motion with azimuth and elevation perspective angles $(55,15)$ and the colormap for time $t\n$ with $\redgradrule \in [0, 15]$.}
    \label{fig_flying_snapshots}
\end{figure}%
\subsection{Flying spaghetti}\label{sec_example_flying}

This first benchmark problem has been initially proposed in \cite{simo_1988_dynamics,simo_1995_nonlinear} and taken up many times by, inter alia, \cite{hsiao_1999_consistent,hesse_2012_consistent,zhang_2016_quadrature,cammarata_2025_implicit,zupan_2019_conservation,marino_2019_isogeometric}. The main objective of this free flight problem is to verify the benefitial properties of the proposed beam formulation and compare the results with data from the literature.
We adopt the parameters chosen in \cite{hesse_2012_consistent,zhang_2016_quadrature}\footnote{The validity of the parameters reported by Simo and co-workers
    have been questioned in \cite{hsiao_1999_consistent,hesse_2012_consistent}.
    With the parameters employed here, the authors of \cite{hesse_2012_consistent,zhang_2016_quadrature} could reproduce the results from the original works.}, see Tab.~\ref{tab_freeflight}.
Initially at rest in an inclined configuration with $\centerline\h(s=0,t=0) = (6,0,0)$ and $\centerline\h(s=L,t=0) = (0,0,8)$,
the beam is subject to an external force and torque at $s=L$, both of which are applied for a short period of time, i.e., $ \vec{n}(s=L,t) = \vec{\nu}(t)$ and $\vec{m}(s=L,t) = \vec{\mu}(t)$ with
\begin{equation} \label{eq_1_loadfunction}
    \vec{\nu}(t) = f(t) \begin{bmatrix}
        \frac{1}{10} \\
        0            \\
        0
    \end{bmatrix}, \quad
    \vec{\mu}(t) = f(t)  \begin{bmatrix}
        0 \\
        1 \\
        \frac{1}{2}
    \end{bmatrix}, \quad
    f(t) = \begin{cases}
        80 t ,    & \text{for}\quad t \leq 2.5,     \\
        400-80t , & \text{for}\quad 2.5 < t \leq 5,
        \\
        0,        & \text{for} \quad t > 5 .
    \end{cases}
\end{equation}%
At $t=5$, the system is closed as both boundary inputs become zero and the beam continues with a free flight, where total energy, linear and angular momentum should be preserved.
The resulting motion
fits well to the depictions from \cite{cammarata_2025_implicit,hesse_2012_consistent,zupan_2019_conservation}, see snapshots in Fig.~\ref{fig_flying_snapshots}. Moreover, the conservation of total energy\footnote{
    References \cite{zhang_2016_quadrature,zupan_2019_conservation,cammarata_2025_implicit}, in which the originally reported parameters are used, report significantly larger energy plateaus but similar angular momentum levels.} can be seen from Fig.~\ref{fig_flying_energy}.
\begin{figure}[tb]
    \setlength{\plotwidth}{0.3\textwidth}
    \setlength{\plotheight}{1\plotwidth}
    \centering
    \tikzsetnextfilename{flyingspaghetti_energy}
    \begin{tikzpicture}

    \begin{axis}[
            width=\plotwidth,
            height=\plotheight,
            xlabel={$t\n$},
            xmin=0, xmax=15,
            ymin=-50, ymax=650,
            xtick distance=5,
            ytick distance=200,
            minor x tick num=1,
            minor y tick num=1,
            grid=both,
            axis lines=box,
            xtick align=center,
            xtick pos = bottom,
            ytick align=center,
            ytick pos = left,
            legend style={%
                    draw=none,                
                    fill=none,                
                    font=\small,             
                    legend columns=-1,
                    at={(axis description cs:0.5,1.01)}, 
                    anchor=south,
                },
        ]

        \addplot+[
            line width=\plotlinewidth,
        ]
        table[col sep=comma,
                x=time,
                y=total_energy] {ex01_flyingspaghetti_results.csv};
        \label{legend_total_energy}

    \end{axis}
\end{tikzpicture}
    \tikzsetnextfilename{flyingspaghetti_energy_balance}
    \begin{tikzpicture}

    \begin{axis}[
            width=\plotwidth,
            height=\plotheight,
            xlabel={$t\n$},
            xmin=0, xmax=15,
            ymin=-1, ymax=25,
            xtick distance=5,
            ytick distance=25,
            minor x tick num=1,
            minor y tick num=3,
            grid=both,
            axis lines=box,
            xtick align=center,
            xtick pos = bottom,
            ytick align=center,
            ytick pos = left,
            legend style={%
                    draw=none,                
                    fill=none,                
                    font=\small,             
                    legend columns=-1,
                    at={(axis description cs:0.5,1.01)}, 
                    anchor=south,
                },
        ]

        \addplot+[
            color2!30,
            smooth,
            line width=\thickplotlinewidth,
        ]
        table[col sep=comma,
                x=time,
                y=external_work] {ex01_flyingspaghetti_results.csv};
        \label{legend_ext_work}

        \addplot+[
            color2!80!black,
            dashed,
            smooth,
            line width=\plotlinewidth,
        ]
        table[col sep=comma,
                x=time,
                y=energy_increment] {ex01_flyingspaghetti_results.csv};
        \label{legend_energy_inc}

    \end{axis}
\end{tikzpicture}
    \tikzsetnextfilename{flyingspaghetti_energydiff}
    \begin{tikzpicture}

    \begin{axis}[
            width=\plotwidth,
            height=\plotheight,
            xlabel={$t\n$},
            xmin=0, xmax=15,
            ymin=-1.5e-11, ymax=1.5e-11,
            xtick distance=5,
            ytick distance=1e-11,
            minor x tick num=1,
            grid=both,
            axis lines=box,
            xtick align=center,
            xtick pos = bottom,
            ytick align=center,
            ytick pos = left,
            legend style={%
                    draw=none,                
                    fill=none,                
                    font=\small,             
                    legend columns=-1,
                    at={(axis description cs:0.5,1.01)}, 
                    anchor=south,
                },
        ]

        \addplot+[
            color = color3,
            line width=\plotlinewidth,
        ]
        table[col sep=comma,
                x=time,
                y=power_balance] {ex01_flyingspaghetti_results.csv};
        \label{legend_power_balance}

        \addplot[color=gray,
            dashed,
            line width = \thinplotlinewidth]
        table[row sep=crcr]{%
                0	1e-11\\
                15	1e-11\\};
        \addplot[color=gray,
            dashed,
            line width = \thinplotlinewidth]
        table[row sep=crcr]{%
                0	-1e-11\\
                15	-1e-11\\};
    \end{axis}
\end{tikzpicture}
    \caption{Flying spaghetti: Total energy \plotref{legend_total_energy} \textcolor{color1}{$\hat{H}(\state\n)$}, input work \plotref{legend_ext_work} \textcolor{color2!80!black}{$W^{\mathrm{ext}}\n$}, total energy increments \plotref{legend_energy_inc} \textcolor{color2!80!black}{$\hat{H}(\nodal{\state}\npe)-\hat{H}(\nodal{\state}\n)$}, and energy balance violation \plotref{legend_power_balance} \textcolor{color3}{$\Delta E\n$}.}
    \label{fig_flying_energy}
\end{figure}
In each time step, the work done by the external load $W^{\mathrm{ext}}\n = h \nodal{u}\npeh\transp \nodal{y}\npeh$ identically matches the total energy increase $\hat{H}(\nodal{\state}\npe)-\hat{H}(\nodal{\state}\n)$ during the loading phase. Correspondingly, at all times, the energy balance violation $\Delta E\n := \hat{H}(\nodal{\state}\npe)-\hat{H}(\nodal{\state}\n) - W^{\mathrm{ext}}\n$ remains at the level of machine precision, confirming the discrete-time power balance \eqref{energy_balance_discrete_time}.
The conservation of total linear and angular momentum after the loading phase can be seen from Fig.~\ref{fig_flying_mom}, confirming the results from Appendix~\ref{sec_ang_mom}.
\begin{figure}[tb]
    \setlength{\plotwidth}{0.4\textwidth}
    \setlength{\plotheight}{0.65\plotwidth}
    \centering
    \tikzsetnextfilename{flyingspaghetti_linmom}
    \begin{tikzpicture}

    \begin{axis}[
            width=\plotwidth,
            height=\plotheight,
            xlabel={$t\n$},
            ylabel={$(p_i)\n$},
            xmin=0, xmax=15,
            ymin=-5, ymax=60,
            xtick distance=5,
            ytick distance=25,
            minor x tick num=1,
            minor y tick num=1,
            grid=both,
            axis lines=box,
            xtick align=center,
            xtick pos = bottom,
            ytick align=center,
            ytick pos = left,
            legend style={%
                    draw=none,                
                    fill=none,                
                    font=\small,             
                    legend columns=-1,
                    at={(axis description cs:0.5,1.01)}, 
                    anchor=south,
                },
        ]

        \addplot+[
            line width=\plotlinewidth,
        ]
        table[col sep=comma,
                x=time,
                y=P_1] {ex01_flyingspaghetti_results.csv};
        \label{legend_mom_1}

        \addplot+[
            line width=\plotlinewidth,
        ]
        table[col sep=comma,
                x=time,
                y=P_2] {ex01_flyingspaghetti_results.csv};
        \label{legend_mom_2}

        \addplot+[
            line width=\plotlinewidth,
            dashed,
        ]
        table[col sep=comma,
                x=time,
                y=P_3] {ex01_flyingspaghetti_results.csv};
        \label{legend_mom_3}

    \end{axis}
\end{tikzpicture}
    \tikzsetnextfilename{flyingspaghetti_angmom}
    \begin{tikzpicture}

    \begin{axis}[
            width=\plotwidth,
            height=\plotheight,
            xlabel={$t\n$},
            ylabel={$(l_i)\n$},
            xmin=0, xmax=15,
            ymin=-50, ymax=600,
            xtick distance=5,
            ytick distance=200,
            minor x tick num=1,
            minor y tick num=1,
            grid=both,
            axis lines=box,
            xtick align=center,
            xtick pos = bottom,
            ytick align=center,
            ytick pos = left,
            legend style={%
                    draw=none,                
                    fill=none,                
                    font=\small,             
                    legend columns=-1,
                    at={(axis description cs:0.5,1.01)}, 
                    anchor=south,
                },
        ]

        \addplot+[
            line width=\plotlinewidth,
        ]
        table[col sep=comma,
                x=time,
                y=L_1] {ex01_flyingspaghetti_results.csv};

        \addplot+[
            line width=\plotlinewidth,
        ]
        table[col sep=comma,
                x=time,
                y=L_2] {ex01_flyingspaghetti_results.csv};

        \addplot+[
            dashed,
            line width=\plotlinewidth,
        ]
        table[col sep=comma,
                x=time,
                y=L_3] {ex01_flyingspaghetti_results.csv};

    \end{axis}
\end{tikzpicture}
    \caption{Flying spaghetti: Total linear momentum components $(p_i)\n = \vec{p}(\nodal{\state}\n)\transp\vec{e}_i$ and total angular momentum components $(l_i)\n = \vec{l}(\nodal{\state}\n)\transp\vec{e}_i$, where $i \in \{$\textcolor{color1}{$1$}\,\plotref{legend_mom_1}, \textcolor{color2}{$2$}\,\plotref{legend_mom_2}, \textcolor{color3}{$3$}\,\plotref{legend_mom_3}$\}$.}
    \label{fig_flying_mom}
\end{figure}
Further, the constant total linear momentum of $\vec{p}=(50,0,0)$ for $t\geq5$, matches the analytical solution exactly within machine precision.
Moreover, the load function \eqref{eq_1_loadfunction}, together with the chosen boundary conditions, allows for an analytical solution for the position of the center of mass of the beam during motion \cite{hsiao_1999_consistent,zhang_2016_quadrature}, given by
\begin{equation} \label{eq:analytic_solution}
    \vec{r} = \begin{bmatrix}
        r_1(t) \\
        0      \\
        4
    \end{bmatrix}
    , \qquad
    r_1(t) =
    \begin{cases}
        3 + \frac{2}{15} t^3  ,                    & \text{for}\quad t \leq 2.5,     \\
        \frac{43}{6} -5t + 2t^2 -\frac{2}{15}t^3 , & \text{for}\quad 2.5 < t \leq 5,
        \\
        -\frac{19}{2} + 5t,                        & \text{for} \quad t > 5 .
    \end{cases}
\end{equation}
This is closely matched by the simulated center of mass $\vec{r}\h(t) := \frac{1}{L}\int_{\domainspace} \centerline\h(s,t) \d s$ obtained with our approach, see Fig.~\ref{fig_flying_center_of_mass}. Some slight deviations $\Delta r_i := \vec{e}_i\transp\vec{r}\h - r_i$ are still present during the loading phase for the $1$-component but after $t=5$, all three curves are in the range of computer precision.
\begin{figure}[tb]
    \setlength{\plotwidth}{0.4\textwidth}
    \setlength{\plotheight}{0.65\plotwidth}
    \centering
    \tikzsetnextfilename{flyingspaghetti_com}
    \begin{tikzpicture}
    \begin{axis}[
            width=\plotwidth,
            height=\plotheight,
            xlabel={$t\n$},
            ylabel={$r_1$},
            xmin=0, xmax=15,
            ymin=0, ymax=70,
            xtick distance=5,
            ytick distance=10,
            minor x tick num=1,   
            minor y tick num=0,
            grid=both,
            axis lines=box,
            xtick align=center,
            xtick pos = bottom,
            ytick align=center,
            ytick pos = left,
            legend style={%
                    draw=none,                
                    fill=none,                
                    font=\small,             
                    legend columns=-1,
                    at={(axis description cs:0.5,1.01)}, 
                    anchor=south,
                },
            samples=200,                 
            domain=0:15,                
        ]

        \addplot+[
            color1!30,
            line width=\thickplotlinewidth,
            domain=0:2.5]{3+2/15*x^3};
        \label{legend_analytical_com}

        \addplot+[
            color1!30,
            line width=\thickplotlinewidth,
            forget plot,
            domain=2.5:5]{43/6-5*x+2*x^2-2/15*x^3};

        \addplot+[
            color1!30,
            line width=\thickplotlinewidth,
            forget plot,
            domain=5:15]{-19/2+5*x};

        \addplot[color=gray,
            line width=\thinplotlinewidth,
            forget plot,
            dashed
        ] table[row sep=crcr]{%
                2.5	0\\
                2.5	70\\};

        \addplot[color=gray,
            line width=\thinplotlinewidth,
            forget plot,
            dashed
        ] table[row sep=crcr]{%
                5	0\\
                5	70\\};

        \addplot+[color1, dashed, line width = \plotlinewidth]
        table[col sep=comma,
                x=time,
                y=center_of_mass_1] {ex01_flyingspaghetti_results.csv};
        \label{legend_numerical_com}

    \end{axis}
\end{tikzpicture}
    \tikzsetnextfilename{flyingspaghetti_comdiff}
    \begin{tikzpicture}
    \begin{axis}[
            width=\plotwidth,
            height=\plotheight,
            xlabel={$t\n$},
            ylabel={$\Delta r_i$},
            xmin=0, xmax=15,
            ymin=1e-20, ymax=1e-2,
            xtick distance=5,
            ytick={1e-20,1e-14,1e-8,1e-2},
            minor x tick num=1,   
            minor y tick num=1,
            grid=both,
            axis lines=box,
            xtick align=center,
            xtick pos = bottom,
            ytick align=center,
            ytick pos = left,
            ymode = log,
            legend style={%
                    draw=none,                
                    fill=none,                
                    font=\small,             
                    legend columns=-1,
                    at={(axis description cs:0.5,1.01)}, 
                    anchor=south,
                },
            samples=200,                 
            domain=0:15,                
        ]

        \addplot[color=gray,
            line width=\thinplotlinewidth,
            forget plot,
            dashed
        ] table[row sep=crcr]{%
                2.5	1e-1\\
                2.5	1e-2\\};

        \addplot[color=gray,
            line width=\thinplotlinewidth,
            forget plot,
            dashed
        ] table[row sep=crcr]{%
                5	1e-1\\
                5	1e-2\\};

        \addplot+[color1, line width = \plotlinewidth]
        table[col sep=comma,
                x=time,
                y=center_of_mass_difference_1] {ex01_flyingspaghetti_results.csv};
        \label{legend_com_diff_1}

        \addplot+[color2, line width = \plotlinewidth]
        table[col sep=comma,
                x=time,
                y=center_of_mass_difference_2] {ex01_flyingspaghetti_results.csv};
        \label{legend_com_diff_2}

        \addplot+[color3, line width = \plotlinewidth]
        table[col sep=comma,
                x=time,
                y=center_of_mass_difference_3] {ex01_flyingspaghetti_results.csv};
        \label{legend_com_diff_3}

    \end{axis}
\end{tikzpicture}
    \caption{Flying spaghetti: Analytical solution (\plotref{legend_analytical_com}) and numerical solution (\plotref{legend_numerical_com}) for $1$-component of the center of mass as well as deviation $\Delta r_i $ with $i \in \{$\textcolor{color1}{$1$}\,\plotref{legend_com_diff_1}, \textcolor{color2}{$2$}\,\plotref{legend_com_diff_2}, \textcolor{color3}{$3$}\,\plotref{legend_com_diff_3}$\}$.}
    \label{fig_flying_center_of_mass}
\end{figure}
Additionally, a convergence study is performed to investigate the temporal discretization accuracy for the centerline position and velocity at the free end with $s=0$ with respect to the temporal discretization. The results in Fig.~\ref{fig_flying_convergence} show that the relative error
\begin{equation}
    e_{\smallsquare} = \frac{\lvert \lvert \square\h(s=0,t=\bar{t}\,) - \square^{\mathrm{ref}}(s=0,t=\bar{t}\,) \rvert \rvert }{\lvert \lvert \square^{\mathrm{ref}}(s=0,t=\bar{t}\,) \rvert \rvert} ,
\end{equation}
for $\square \in \{\centerline, \vphi\}$,
obtained with our approach decreases with order two with increasing temporal resolution. This is in line with the expected behavior for the midpoint-type discretization. Here, the evaluation time is chosen as $\bar{t}=5$. Moreover, $\square^{\mathrm{ref}}$ is the finely discretized reference solution obtained with $h=10^{-3}$ and $\lvert \lvert \square \rvert \rvert$ denotes the Euclidean norm. In this context, the Newton tolerance was reduced for all simulations to $10^{-8}$.

Another notable advantage of the time integration method, which is made possible by the PH formulation, is its remarkable stability and robustness. In theory, the time step size is not constrained by the stability of the integrator but only by the deterioration of the Newton method due to poor initial guesses at excessively large steps. In practice, the method must adequately resolve the loading function profile \eqref{eq_1_loadfunction}. Notably, choosing $h = 2.5$ poses no difficulty for the proposed method, while other works simulating this example typically report time step sizes of $h=0.01$ or $0.02$. However, in light of the preceding convergence analysis, the accuracy of the results with larger time step sizes should definitly be questioned.

\begin{figure}[t]
    \setlength{\plotwidth}{0.45\textwidth}
    \setlength{\plotheight}{0.8\plotwidth}
    \centering
    \tikzsetnextfilename{flyingspaghetti_convergence}
    \begin{tikzpicture}

    \begin{axis}[%
            width=\plotwidth,
            height=\plotheight,
            xmode=log,
            ymode=log,
            xmin=1e-2, xmax=1,
            ymin=1e-5, ymax=10,
            ytick={1e-5,1e-4,1e-3,1e-2,1e-1,1,10},
            yticklabels={$10^{-5}$,,$10^{-3}$,,$10^{-1}$,,$10^1$},
            grid=both,
            axis lines=box,
            xtick align=center,
            xtick pos = bottom,
            ytick align=center,
            ytick pos = left,
            xlabel={$h$},
            ylabel={relative error $e_{\smallsquare}$},
            legend style={%
                    draw=none,                
                    fill=none,                
                    font=\small,             
                    legend columns=-1,
                    at={(axis description cs:0.5,1.01)}, 
                    anchor=south,
                },
        ]

        \addplot+[line width=\plotlinewidth,
            mark=*]
        table[col sep=comma, x=timestepsize, y=position_error] {ex01_flyingspaghetti_convergence_results.csv};
        \label{legend_convergence_phi}

        \addplot+[line width=\plotlinewidth, mark=square*]
        table[col sep=comma, x=timestepsize, y=velocity_error] {ex01_flyingspaghetti_convergence_results.csv};
        \label{legend_convergence_vphi}

        \coordinate (A) at (0.02,3.6e-4);
        \coordinate (B) at (0.04,1.45e-3);
        \draw (A) to (B);

        \addplot [draw = none, forget plot]
        table[row sep=crcr,y={create col/linear regression}]{%
                0.05 0.00230121175116172\\
                0.01 8.9958710916992e-05\\
            };
        \xdef\slope{\pgfplotstableregressiona}

        \draw (A) -| (B) node [pos=0.75,anchor=west] {\pgfmathprintnumber[fixed, fixed zerofill, precision=4]\slope};

    \end{axis}
\end{tikzpicture}%
    \caption{Flying spaghetti: Convergence of the relative error for centerline position \textcolor{color1}{$e_{\centerline}$} \plotref{legend_convergence_phi} and velocity \textcolor{color2}{$e_{\vphi}$} \plotref{legend_convergence_vphi} with varying time step size parameter.}
    \label{fig_flying_convergence}
\end{figure}

\subsection{Nonlinear oscillation of a cantilever}

\begin{table}[bt]
    \centering
    \caption{Nonlinear oscillation of a cantilever: Physical and numerical parameters.}
    \begin{tabular}{l l l l l l l l l l l}
        \toprule
        $h$       & $\finaltime$ & $\nel$ & $\epsilon$ & $L$ & $\rho$   & $A$           & $M_\rho^{11} = M_\rho^{22}$ & $k_{\mathrm{b1}} = k_{\mathrm{b2}}$ & $k_{\mathrm{t}}$ & $k_{\mathrm{s1}}=k_{\mathrm{s2}}=k_{\mathrm{e}}$ \\
        \midrule
        $10^{-3}$ & $0.3$        & $8$    & $10^{-12}$ & $1$ & $ 2850 $ & $\pi d^2 / 4$ & $ \rho I $                  & $ EI $                              & $G I_T$          & $\infty$                                         \\
        \bottomrule
    \end{tabular}
    \label{tab_cantilever}
\end{table}

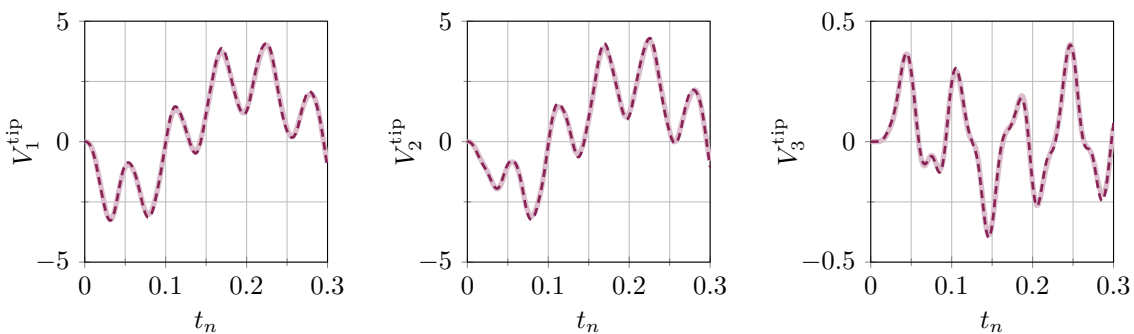
\begin{figure}[bt]
    \setlength{\plotwidth}{0.3\textwidth}
    \setlength{\plotheight}{1\plotwidth}
    \centering
    \tikzsetnextfilename{cantilever_tipvelo1}
    \begin{tikzpicture}

    \begin{axis}[
            width=\plotwidth,
            height=\plotheight,
            xlabel={$t\n$},
            ylabel={$V_1^{\mathrm{tip}}$},
            xmin=0, xmax=0.3,
            ymin=-5, ymax=5,
            xtick distance=0.1,
            ytick distance=5,
            minor x tick num=1,
            minor y tick num=1,
            grid=both,
            axis lines=box,
            xtick align=center,
            xtick pos = bottom,
            ytick align=center,
            ytick pos = left,
            ylabel style={yshift=-8pt},
            legend style={%
                    draw=none,                
                    fill=none,                
                    font=\small,             
                    legend columns=-1,
                    at={(axis description cs:0.5,1.01)}, 
                    anchor=south,
                },
        ]

        \addplot+[
            color1!30,
            line width=\thickplotlinewidth,
        ]
        table[col sep=comma,
                x=time,
                y=tip_velocity_B_1] {ex02_cantilever_reference_results.csv};

        \addplot+[
            color1,
            line width=\plotlinewidth,
            densely dashed
        ]
        table[col sep=comma,
                x=time,
                y=tip_velocity_B_1] {ex02_cantilever_results.csv};

    \end{axis}
\end{tikzpicture}
    \tikzsetnextfilename{cantilever_tipvelo2}
    \begin{tikzpicture}

    \begin{axis}[
            width=\plotwidth,
            height=\plotheight,
            xlabel={$t\n$},
            ylabel={$V_2^{\mathrm{tip}}$},
            xmin=0, xmax=0.3,
            ymin=-5, ymax=5,
            xtick distance=0.1,
            ytick distance=5,
            minor x tick num=1,
            minor y tick num=1,
            grid=both,
            axis lines=box,
            xtick align=center,
            xtick pos = bottom,
            ytick align=center,
            ytick pos = left,
            ylabel style={yshift=-8pt},
            legend style={%
                    draw=none,                
                    fill=none,                
                    font=\small,             
                    legend columns=-1,
                    at={(axis description cs:0.5,1.01)}, 
                    anchor=south,
                },
        ]

        \addplot+[
            color1!30,
            line width=\thickplotlinewidth,
        ]
        table[col sep=comma,
                x=time,
                y=tip_velocity_B_2] {ex02_cantilever_reference_results.csv};
        \label{legend_HK24}

        \addplot+[
            color1,
            line width=\plotlinewidth,
            densely dashed
        ]
        table[col sep=comma,
                x=time,
                y=tip_velocity_B_2] {ex02_cantilever_results.csv};
        \label{legend_PHMP}

    \end{axis}
\end{tikzpicture}
    \tikzsetnextfilename{cantilever_tipvelo3}
    \begin{tikzpicture}

    \begin{axis}[
            width=\plotwidth,
            height=\plotheight,
            xlabel={$t\n$},
            ylabel={$V_3^{\mathrm{tip}}$},
            xmin=0, xmax=0.3,
            ymin=-0.5, ymax=0.5,
            xtick distance=0.1,
            ytick distance=0.5,
            minor x tick num=1,
            minor y tick num=1,
            grid=both,
            axis lines=box,
            xtick align=center,
            xtick pos = bottom,
            ytick align=center,
            ytick pos = left,
            ylabel style={yshift=-8pt},
            legend style={%
                    draw=none,                
                    fill=none,                
                    font=\small,             
                    legend columns=-1,
                    at={(axis description cs:0.5,1.01)}, 
                    anchor=south,
                },
        ]

        \addplot+[
            color1!30,
            line width=\thickplotlinewidth,
        ]
        table[col sep=comma,
                x=time,
                y=tip_velocity_B_3] {ex02_cantilever_reference_results.csv};

        \addplot+[
            color1,
            line width=\plotlinewidth,
            densely dashed
        ]
        table[col sep=comma,
                x=time,
                y=tip_velocity_B_3] {ex02_cantilever_results.csv};

    \end{axis}
\end{tikzpicture}
    \caption{Nonlinear oscillation of a cantilever: Tip velocities in cross-section frame $V_i^{\mathrm{tip}}$. Results from our approach (\plotref{legend_PHMP}) and reference solution (\plotref{legend_HK24}) \cite{herrmann_2024_relativekinematic} coincide. The sign flip and different order of components compared to the reference is due to a different director frame orientation.}
    \label{fig_cantilever_tip_velos}
\end{figure}
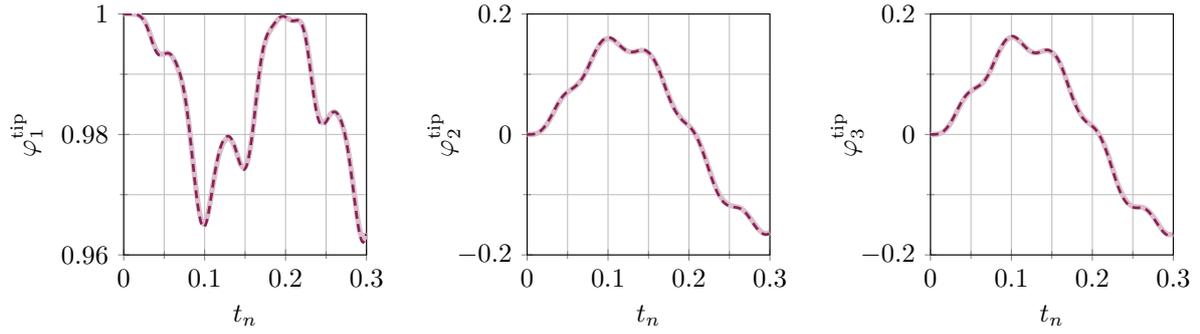
\begin{figure}[tb]
    \setlength{\plotwidth}{0.3\textwidth}
    \setlength{\plotheight}{1\plotwidth}
    \centering
    \tikzsetnextfilename{cantilever_tippos1}
    \begin{tikzpicture}

    \begin{axis}[
            width=\plotwidth,
            height=\plotheight,
            xlabel={$t\n$},
            ylabel={$\varphi_1^{\mathrm{tip}}$},
            xmin=0, xmax=0.3,
            ymin=0.96, ymax=1,
            xtick distance=0.1,
            ytick distance=0.02,
            minor x tick num=1,
            minor y tick num=1,
            grid=both,
            axis lines=box,
            xtick align=center,
            xtick pos = bottom,
            ytick align=center,
            ytick pos = left,
            legend style={%
                    draw=none,                
                    fill=none,                
                    font=\small,             
                    legend columns=-1,
                    at={(axis description cs:0.5,1.01)}, 
                    anchor=south,
                },
        ]

        \addplot+[
            color1!30,
            line width=\thickplotlinewidth,
        ]
        table[col sep=comma,
                x=time,
                y=tip_position_I_1] {ex02_cantilever_reference_results.csv};

        \addplot+[
            color1,
            line width=\plotlinewidth,
            densely dashed
        ]
        table[col sep=comma,
                x=time,
                y=tip_position_I_1] {ex02_cantilever_results.csv};

    \end{axis}
\end{tikzpicture}
    \tikzsetnextfilename{cantilever_tippos2}
    \begin{tikzpicture}

    \begin{axis}[
            width=\plotwidth,
            height=\plotheight,
            xlabel={$t\n$},
            ylabel={$\varphi_2^{\mathrm{tip}}$},
            xmin=0, xmax=0.3,
            ymin=-0.2, ymax=0.2,
            xtick distance=0.1,
            ytick distance=0.2,
            minor x tick num=1,
            minor y tick num=1,
            grid=both,
            axis lines=box,
            xtick align=center,
            xtick pos = bottom,
            ytick align=center,
            ytick pos = left,
            ylabel style={yshift=-8pt},
            legend style={%
                    draw=none,                
                    fill=none,                
                    font=\small,             
                    legend columns=-1,
                    at={(axis description cs:0.5,1.01)}, 
                    anchor=south,
                },
        ]

        \addplot+[
            color1!30,
            line width=\thickplotlinewidth,
        ]
        table[col sep=comma,
                x=time,
                y=tip_position_I_2] {ex02_cantilever_reference_results.csv};

        \addplot+[
            color1,
            line width=\plotlinewidth,
            densely dashed
        ]
        table[col sep=comma,
                x=time,
                y=tip_position_I_2] {ex02_cantilever_results.csv};

    \end{axis}
\end{tikzpicture}
    \tikzsetnextfilename{cantilever_tippos3}
    \begin{tikzpicture}

    \begin{axis}[
            width=\plotwidth,
            height=\plotheight,
            xlabel={$t\n$},
            ylabel={$\varphi_3^{\mathrm{tip}}$},
            xmin=0, xmax=0.3,
            ymin=-0.2, ymax=0.2,
            xtick distance=0.1,
            ytick distance=0.2,
            minor x tick num=1,
            minor y tick num=1,
            grid=both,
            axis lines=box,
            xtick align=center,
            xtick pos = bottom,
            ytick align=center,
            ytick pos = left,
            ylabel style={yshift=-8pt},
            legend style={%
                    draw=none,                
                    fill=none,                
                    font=\small,             
                    legend columns=-1,
                    at={(axis description cs:0.5,1.01)}, 
                    anchor=south,
                },
        ]

        \addplot+[
            color1!30,
            line width=\thickplotlinewidth,
        ]
        table[col sep=comma,
                x=time,
                y=tip_position_I_3] {ex02_cantilever_reference_results.csv};

        \addplot+[
            color1,
            line width=\plotlinewidth,
            densely dashed
        ]
        table[col sep=comma,
                x=time,
                y=tip_position_I_3] {ex02_cantilever_results.csv};

    \end{axis}
\end{tikzpicture}
    \caption{Nonlinear oscillation of a cantilever: Tip positions in inertial frame $\varphi_i^{\mathrm{tip}}$. Results from our approach (\plotref{legend_PHMP}) and reference solution (\plotref{legend_HK24}) \cite{herrmann_2024_relativekinematic}. The sign flip and different order of components compared to the reference is due to a different director frame orientation.}
    \label{fig_cantilever_tip_pos}
\end{figure}

The example of an oscillating slender aluminum beam has been recently included in \cite{herrmann_2024_relativekinematic}. The beam has a circular cross-section with diameter $d = 4 \cdot 10^{-3} $ and we assume material parameters for Young's modulus $E=7.2 \cdot 10^{10}$, Poisson ratio $\nu = 0.35$ and shear modulus $G = E / (2 (1 + \nu))$. This leads to cross-sectional moments of area $I_1 = I_2 = 1/2 I_T = \pi d^4/64$. Further details are comprised in
Tab.~\ref{tab_cantilever}.
Note that this beam is relatively slender with a high slenderness ratio of $L/d = 250$. Conforming with \cite{herrmann_2024_relativekinematic}, we thus assume an inextensible Kirchhoff beam and impose the related kinematic constraints of shear rigidity and inextensibility. To this end, we set $\vec{C}_N = \vec{0}$ in the compliance equation \eqref{eq:strain_evo_stress_resultants}, thus ensuring $\dot{\vec{\Gamma}}=\vec{0}$.
Moreover, note that we are able to choose a relatively large time step size of $h=10^{-3}$ due to the robustness of the employed energy-momentum time integration \eqref{timediscrete_EoM} (see paragraph in the first example). Compared to the original setting in \cite{herrmann_2024_relativekinematic}, where the largest possible time step size was $h = 4\cdot 10^{-5}$, we obtain results even for $h=5 \cdot 10^{-2}$.
\begin{figure}[bt]
    \setlength{\plotwidth}{0.4\textwidth}
    \setlength{\plotheight}{0.65\plotwidth}
    \centering
    \tikzsetnextfilename{cantilever_constraint1}
    \begin{tikzpicture}

    \begin{axis}[
            width=\plotwidth,
            height=\plotheight,
            xlabel={$t\n$},
            ylabel={$g\h_i $},
            xmin=0, xmax=0.3,
            ymin=-1.5e-14, ymax=1.5e-14,
            xtick distance=0.1,
            ytick distance=1e-14,
            minor x tick num=1,
            minor y tick num=0,
            grid=both,
            axis lines=box,
            xtick align=center,
            xtick pos = bottom,
            ytick align=center,
            ytick pos = left,
            legend style={%
                    draw=none,                
                    fill=none,                
                    font=\small,             
                    legend columns=-1,
                    at={(axis description cs:0.5,1.01)}, 
                    anchor=south,
                },
        ]

        \addplot+[
            line width=\plotlinewidth,
        ]
        table[col sep=comma,
                x=time,
                y=mid_constraint_1] {ex02_cantilever_results.csv};
        \label{legend_constraint_1}

        \addplot+[
            line width=\plotlinewidth,
        ]
        table[col sep=comma,
                x=time,
                y=mid_constraint_2] {ex02_cantilever_results.csv};
        \label{legend_constraint_2}

        \addplot+[
            line width=\plotlinewidth,
        ]
        table[col sep=comma,
                x=time,
                y=mid_constraint_3] {ex02_cantilever_results.csv};
        \label{legend_constraint_3}

        \addplot[color=gray,
            dashed,
            line width = \thinplotlinewidth]
        table[row sep=crcr]{%
                0	1e-14\\
                0.3	1e-14\\};
        \addplot[color=gray,
            dashed,
            line width = \thinplotlinewidth]
        table[row sep=crcr]{%
                0	-1e-14\\
                0.3	-1e-14\\};

    \end{axis}
\end{tikzpicture}
    \tikzsetnextfilename{cantilever_constraint2}
    \begin{tikzpicture}

    \begin{axis}[
            width=\plotwidth,
            height=\plotheight,
            xlabel={$t\n$},
            xmin=0, xmax=0.3,
            ymin=-1.5e-15, ymax=1.5e-15,
            xtick distance=0.1,
            ytick distance=1e-15,
            minor x tick num=1,
            minor y tick num=0,
            grid=both,
            axis lines=box,
            xtick align=center,
            xtick pos = bottom,
            ytick align=center,
            ytick pos = left,
            legend style={%
                    draw=none,                
                    fill=none,                
                    font=\small,             
                    legend columns=-1,
                    at={(axis description cs:0.5,1.01)}, 
                    anchor=south,
                },
        ]

        \addplot+[
            color4,
            line width=\plotlinewidth,
        ]
        table[col sep=comma,
                x=time,
                y=mid_constraint_4] {ex02_cantilever_results.csv};
        \label{legend_constraint_4}

        \addplot+[
            color5,
            line width=\plotlinewidth,
        ]
        table[col sep=comma,
                x=time,
                y=mid_constraint_5] {ex02_cantilever_results.csv};
        \label{legend_constraint_5}

        \addplot+[
            color6,
            line width=\plotlinewidth,
        ]
        table[col sep=comma,
                x=time,
                y=mid_constraint_6] {ex02_cantilever_results.csv};
        \label{legend_constraint_6}

        \addplot[color=gray,
            dashed,
            line width = \thinplotlinewidth]
        table[row sep=crcr]{%
                0	1e-15\\
                0.3	1e-15\\};

        \addplot[color=gray,
            dashed,
            line width = \thinplotlinewidth]
        table[row sep=crcr]{%
                0	-1e-15\\
                0.3	-1e-15\\};

    \end{axis}
\end{tikzpicture}
    \caption{Nonlinear oscillation of a cantilever: Orthonormality constraints $g\h_i$ with $i \in \{$\textcolor{color1}{$1$}\,\plotref{legend_constraint_1}, \textcolor{color2}{$2$}\,\plotref{legend_constraint_2}, \textcolor{color3}{$3$}\,\plotref{legend_constraint_3}, \textcolor{color4}{$4$}\,\plotref{legend_constraint_4}, \textcolor{color5}{$5$}\,\plotref{legend_constraint_5}, \textcolor{color6}{$6$}\,\plotref{legend_constraint_6}$\}$.}
    \label{fig_cantilever_drift}
\end{figure}
\begin{figure}[tb]
    \setlength{\plotwidth}{0.4\textwidth}
    \setlength{\plotheight}{0.65\plotwidth}
    \centering
    \tikzsetnextfilename{cantilever_drift}
    \begin{tikzpicture}

    \begin{axis}[
            width=\plotwidth,
            height=\plotheight,
            xlabel={$t\n$},
            ylabel={$\lvert\lvert \Gamma_i \rvert\rvert_\domainspace$},
            xmin=0, xmax=0.3,
            ymin=-0.5e-4, ymax=9e-4,
            xtick distance=0.1,
            ytick distance=3e-4,
            minor x tick num=1,
            minor y tick num=0,
            grid=both,
            axis lines=box,
            xtick align=center,
            xtick pos = bottom,
            ytick align=center,
            ytick pos = left,
            legend style={%
                    draw=none,                
                    fill=none,                
                    font=\small,             
                    legend columns=-1,
                    at={(axis description cs:0.5,1.01)}, 
                    anchor=south,
                },
        ]

        \addplot+[
            line width=\plotlinewidth,
        ]
        table[col sep=comma,
                x=time,
                y=Gamma_1] {ex02_cantilever_results.csv};
        \label{legend_strain_1}

        \addplot+[
            line width=\plotlinewidth,
        ]
        table[col sep=comma,
                x=time,
                y=Gamma_2] {ex02_cantilever_results.csv};
        \label{legend_strain_2}

        \addplot+[
            line width=\plotlinewidth,
        ]
        table[col sep=comma,
                x=time,
                y=Gamma_3] {ex02_cantilever_results.csv};
        \label{legend_strain_3}



    \end{axis}
\end{tikzpicture}
    \tikzsetnextfilename{cantilever_drift2}
    \begin{tikzpicture}

    \begin{axis}[
            width=\plotwidth,
            height=\plotheight,
            xlabel={$t\n$},
            ylabel={$ \lvert\lvert \Delta K_j \rvert\rvert_\domainspace$},
            xmin=0, xmax=0.3,
            ymin=-0.5e-4, ymax=9e-4,
            xtick distance=0.1,
            ytick distance=3e-4,
            minor x tick num=1,
            minor y tick num=0,
            grid=both,
            axis lines=box,
            xtick align=center,
            xtick pos = bottom,
            ytick align=center,
            ytick pos = left,
            legend style={%
                    draw=none,                
                    fill=none,                
                    font=\small,             
                    legend columns=-1,
                    at={(axis description cs:0.5,1.01)}, 
                    anchor=south,
                },
        ]

        \addplot+[
            color4,
            line width=\plotlinewidth,
        ]
        table[col sep=comma,
                x=time,
                y=drift_K_1] {ex02_cantilever_results.csv};
        \label{legend_strain_4}

        \addplot+[
            color5,
            line width=\plotlinewidth,
        ]
        table[col sep=comma,
                x=time,
                y=drift_K_2] {ex02_cantilever_results.csv};
        \label{legend_strain_5}

        \addplot+[
            color6,
            line width=\plotlinewidth,
        ]
        table[col sep=comma,
                x=time,
                y=drift_K_3] {ex02_cantilever_results.csv};
        \label{legend_strain_6}



    \end{axis}
\end{tikzpicture}
    \caption{
        Nonlinear oscillation of a cantilever:
        $L^2$-norms of the shear and dilatation strain components computed from displacement quantities, i.e.,
        $\lvert\lvert \Gamma_i\rvert\rvert_\domainspace$ with $i \in \{$\textcolor{color1}{$1$}\,\plotref{legend_strain_1}, \textcolor{color2}{$2$}\,\plotref{legend_strain_2}, \textcolor{color3}{$3$}\,\plotref{legend_strain_3}$\}$. $L^2$-norms of the drift between the curvature and torsion strain measures computed from the stress quantities and the ones computed from displacement quantities, denoted as$\lvert\lvert \Delta K_j \rvert\rvert_\domainspace$ for
        $j \in \{$
        \textcolor{color4}{$1$}\,\plotref{legend_strain_4}, \textcolor{color5}{$2$}\,\plotref{legend_strain_5}, \textcolor{color6}{$3$}\,\plotref{legend_strain_6}$\}$.
    }
    \label{fig_cantilever_strains}
\end{figure}
Initially at rest and positioned along the $\vec{e}_1$-direction, 
the beam is subject to an external force and torque at $s=L$, both of which are applied for a short period of time, i.e., $ \vec{n}(s=L,t) = \vec{\nu}(t)$ and $\vec{m}(s=L,t) = \vec{\mu}(t)$ with
\begin{equation} \label{eq_2_loadfunction}
    \vec{\nu}(t) = f(t)  \begin{bmatrix}
        0 \\
        1 \\
        1
    \end{bmatrix} \quad
    \vec{\mu}(t) = f(t) \begin{bmatrix}
        0.25 \\
        0    \\
        0
    \end{bmatrix}, \quad
    f(t) =
    \begin{cases}
        \frac{1}{2}\left(1 - \cos\left( \frac{2 \pi}{0.05} t \right) \right) , & \text{for}\quad t \leq 0.05, \\
        0,                                                                     & \text{for} \quad t > 0.05 .
    \end{cases}
\end{equation}

We can observe good agreements of the material tip velocities $V_i^{\mathrm{tip}} := \vec{e}_i\transp\rotmat(\vec{d}\h(s=L,t=t\n))\transp\vphi\h(s=L,t=t\n)$, see Fig.~\ref{fig_cantilever_tip_velos}, and spatial tip centerline positions $\varphi_i^{\mathrm{tip}} := \vec{e}_i\transp\centerline\h(s=L,t=t\n)$, see Fig.~\ref{fig_cantilever_tip_pos}, with the results from \cite{herrmann_2024_relativekinematic}. Furthermore, confirming \eqref{eq_proof_no_drift}, it can be verified that the orthonormality constraints \eqref{eq:ortho_constraint} for the directors at the center FE node $g\h_i := g_i(\vec{d}\h(s=L/2,t=t\n))$ are numerically zero during the whole simulation, see Fig.~\ref{fig_cantilever_drift}.
These results are representative also for the other nodes.
Additionally, also the kinematic accuracy can be verified.
The dilatation and shear strain components derived from the mixed variables vanish identically. When computed from the displacements by means of \eqref{eq:strain_measures_short}, these strains are also variationally consistent with the imposed kinematic constraints (see the variational identity \eqref{eq_proof_consistent_strain_new}). The curvature and torsion components are likewise approximated consistently. The results depicted in Fig. \ref{fig_cantilever_strains} are expressed in terms of the $L^2$-norms over the spatial domain, i.e., $\lvert\lvert \Gamma_i\rvert\rvert_\domainspace := (\int_0^L \Gamma_i(\vec{q}\h\n,\partial_s \vec{q}\h\n)^2  \d s)^{1/2}$ and $\lvert\lvert \Delta K_j \rvert\rvert_\domainspace := (\int_0^L ((K_j)\n - K_j(\vec{q}\h\n,\partial_s \vec{q}\h\n))^2 \d s)^{1/2}$, respectively.

\begin{figure}[tb]
    \setlength{\plotwidth}{0.33\textwidth}
    \setlength{\plotheight}{1\plotwidth}
    \centering
    \tikzsetnextfilename{cantilever_energy}
    \begin{tikzpicture}

    \begin{axis}[
            width=\plotwidth,
            height=\plotheight,
            xlabel={$t\n$},
            ylabel={$\hat{H}(\nodal{\state}\n)$},
            xmin=0, xmax=0.3,
            ymin=0, ymax=0.1,
            xtick distance=0.1,
            ytick distance=0.1,
            minor x tick num=1,
            minor y tick num=3,
            grid=both,
            axis lines=box,
            xtick align=center,
            xtick pos = bottom,
            ytick align=center,
            ytick pos = left,
            legend style={%
                    draw=none,                
                    fill=none,                
                    font=\small,             
                    legend columns=-1,
                    at={(axis description cs:0.5,1.01)}, 
                    anchor=south,
                },
        ]

        \addplot+[
            color1!30,
            line width=\thickplotlinewidth,
        ]
        table[col sep=comma,
                x=time,
                y=total_energy] {ex02_cantilever_results.csv};
        \label{legend_cantilever_energy_infty}

        \addplot+[
            color1,
            line width=\plotlinewidth,
            densely dashed
        ]
        table[col sep=comma,
                x=time,
                y=total_energy] {ex02_cantilever_dissipative_results.csv};
        \label{legend_cantilever_energy_008}

    \end{axis}
\end{tikzpicture}
    \tikzsetnextfilename{cantilever_dissipation}
    \begin{tikzpicture}

    \begin{axis}[
            width=\plotwidth,
            height=\plotheight,
            xlabel={$t\n$},
            xmin=0, xmax=0.3,
            ymin=-0.0001, ymax=0.001,
            xtick distance=0.1,
            ytick distance=0.0005,
            minor x tick num=1,
            minor y tick num=1,
            grid=both,
            axis lines=box,
            xtick align=center,
            xtick pos = bottom,
            ytick align=center,
            ytick pos = left,
            legend style={%
                    draw=none,                
                    fill=none,                
                    font=\small,             
                    legend columns=-1,
                    at={(axis description cs:0.5,1.01)}, 
                    anchor=south,
                },
        ]

        \addplot+[
            line width=\plotlinewidth,
            color6,
        ]
        table[col sep=comma,
                x=time,
                y=dissipated_energy] {ex02_cantilever_dissipative_results.csv};
        \label{legend_cantilever_dissipation}

        \addplot[color=gray,
            dashed,
            line width = \thinplotlinewidth]
        table[row sep=crcr]{%
                0	0\\
                0.3	0\\};

    \end{axis}
\end{tikzpicture}
    \tikzsetnextfilename{cantilever_energydiff}
    \begin{tikzpicture}

    \begin{axis}[
            width=\plotwidth,
            height=\plotheight,
            xlabel={$t\n$},
            xmin=0, xmax=0.3,
            ymin=-1.5e-14, ymax=1.5e-14,
            xtick distance=0.1,
            ytick distance=1e-14,
            minor x tick num=1,
            grid=both,
            axis lines=box,
            xtick align=center,
            xtick pos = bottom,
            ytick align=center,
            ytick pos = left,
            legend style={%
                    draw=none,                
                    fill=none,                
                    font=\small,             
                    legend columns=-1,
                    at={(axis description cs:0.5,1.01)}, 
                    anchor=south,
                },
        ]

        \addplot+[
            color = color3,
            line width=\plotlinewidth,
        ]
        table[col sep=comma,
                x=time,
                y=power_balance] {ex02_cantilever_dissipative_results.csv};
        \label{legend_cantilever_energy_diff}

        \addplot[color=gray,
            dashed,
            line width = \thinplotlinewidth]
        table[row sep=crcr]{%
                0	1e-14\\
                15	1e-14\\};
        \addplot[color=gray,
            dashed,
            line width = \thinplotlinewidth]
        table[row sep=crcr]{%
                0	-1e-14\\
                15	-1e-14\\};
    \end{axis}
\end{tikzpicture}
    \caption{Nonlinear oscillation of a cantilever: Total energy $\hat{H}(\state\n)$ with dissipation ( \textcolor{color1}{$\tau = 0.08$} \plotref{legend_cantilever_energy_008}) and without dissipation (\textcolor{color1}{$\tau \rightarrow \infty$} \plotref{legend_cantilever_energy_infty}). Strictly positive dissipated energy per timestep \textcolor{color6}{$\hat{D}\n$} \plotref{legend_cantilever_dissipation} and energy balance violation \textcolor{color3}{$\Delta E\n$} \plotref{legend_cantilever_energy_diff}.
    }
    \label{fig_cantilever_energy}
\end{figure}
Lastly, we consider the visco-elastic generalized Maxwell model, introduced in \Cref{sec_visco}, with $m=1$ viscous branch and $\tau_E = \tau_G = \tau = 0.08$ as relaxation time. The stiffness parameters are split up consistently on both branches, with fractions of $1/4$ in the purely elastic and $3/4$ in the viscous branch. Note that for this simulation with viscous effects, we also keep the kinematic assumptions from above, such that only the bending and torsion modes are damped and $\vec{V}_\Gamma^{-1} = \vec{0}$. Correspondingly, Fig.~\ref{fig_cantilever_energy} extends the findings from the flying spaghetti problem, see Sec.~\ref{sec_example_flying}, to verify the consistent discrete-time representation of the power balance equation. The total energy decreases over time and the dissipated energy $\hat{D}\n := h \nodal{z}(\nodal{x}\npeh)\transp \mathcal{R}\nodal{z}(\nodal{x}\npeh)$ is strictly positive. Thus, the subfigure on the right-hand side of Fig.~\ref{fig_cantilever_energy} shows that the violation of the discrete time energy balance \eqref{energy_balance_discrete_time} denoted as $\Delta E\n := \hat{H}(\nodal{\state}\npe)-\hat{H}(\nodal{\state}\n) - W^{\mathrm{ext}}\n - \hat{D}(\nodal{\state}\n)$ is numerically zero.

\subsection{Quasistatic cantilever problem}
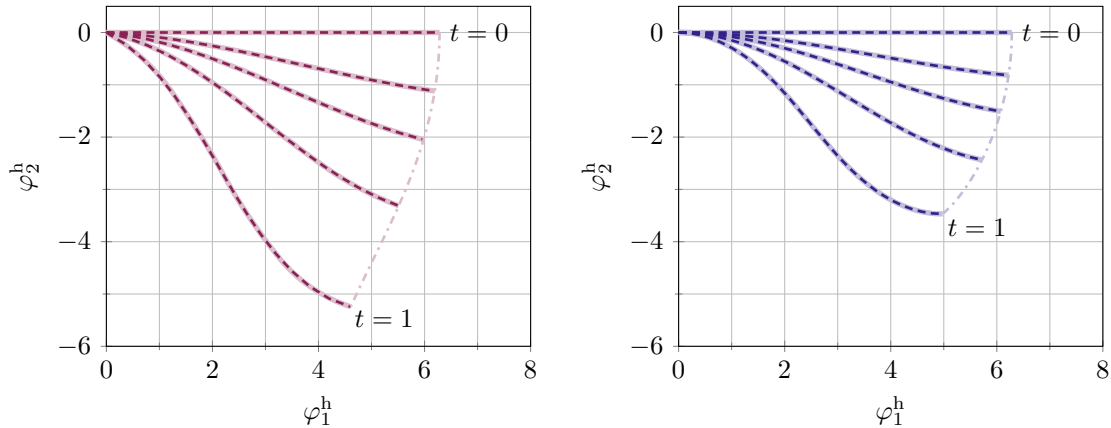
\begin{figure}[t]
    \setlength{\plotwidth}{0.45\textwidth}
    \setlength{\plotheight}{0.85\plotwidth}
    \centering
    \tikzsetnextfilename{quasistatic_unconstrained}
    \begin{tikzpicture}

    \begin{axis}[
            width=\plotwidth,
            height=\plotheight,
            xlabel={$\varphi\h_1$},
            ylabel={$\varphi\h_2$},
            xmin=0, xmax=8,
            ymin=-6, ymax=0.5,
            xtick distance=2,
            ytick distance=2,
            minor x tick num=1,
            minor y tick num=1,
            grid=both,
            axis lines=box,
            xtick align=center,
            xtick pos = bottom,
            ytick align=center,
            ytick pos = left,
            legend style={%
                    draw=none,                
                    fill=none,                
                    font=\small,             
                    legend columns=-1,
                    at={(axis description cs:0.5,1.01)}, 
                    anchor=south,
                },
        ]

        \node[anchor= west] at (axis cs:6.3,0) {$t=0$};
        \node[anchor=north west] at (axis cs:4.51,-5.14) {$t=1$};

        \addplot+[
            color1!30,
            line width=\thickplotlinewidth,
        ]
        table[col sep=comma,
                x=x,
                y=y] {unconstrained_MH_0.csv};
        \label{legend_HE25}

        \addplot+[
            color1!30,
            line width=\thickplotlinewidth,
            forget plot,
        ]
        table[col sep=comma,
                x=x,
                y=y] {unconstrained_MH_1.csv};

        \addplot+[
            color1!30,
            line width=\thickplotlinewidth,
            forget plot,
        ]
        table[col sep=comma,
                x=x,
                y=y] {unconstrained_MH_2.csv};

        \addplot+[
            color1!30,
            line width=\thickplotlinewidth,
            forget plot,
        ]
        table[col sep=comma,
                x=x,
                y=y] {unconstrained_MH_4.csv};

        \addplot+[
            color1!30,
            line width=\thickplotlinewidth,
            forget plot,
        ]
        table[col sep=comma,
                x=x,
                y=y] {unconstrained_MH_10.csv};

        \addplot+[
            color1,
            line width=\plotlinewidth,
            densely dashed,
        ]
        table[col sep=comma,
                x=x_0,
                y=y_0] {ex03_quasistatic_cantilever_results.csv};
        \label{legend_ours_unconstrained}

        \addplot+[
            color1,
            line width=\plotlinewidth,
            densely dashed,
            forget plot
        ]
        table[col sep=comma,
                x=x_1,
                y=y_1] {ex03_quasistatic_cantilever_results.csv};

        \addplot+[
            color1,
            line width=\plotlinewidth,
            densely dashed,
            forget plot
        ]
        table[col sep=comma,
                x=x_2,
                y=y_2] {ex03_quasistatic_cantilever_results.csv};

        \addplot+[
            color1,
            line width=\plotlinewidth,
            densely dashed,
            forget plot
        ]
        table[col sep=comma,
                x=x_4,
                y=y_4] {ex03_quasistatic_cantilever_results.csv};

        \addplot+[
            color1,
            line width=\plotlinewidth,
            densely dashed,
            forget plot
        ]
        table[col sep=comma,
                x=x_10,
                y=y_10] {ex03_quasistatic_cantilever_results.csv};

        \addplot+[
            color1!30,
            line width=\plotlinewidth,
            dash dot,
        ]
        table[col sep=comma,
                x=x,
                y=y] {unconstrained_MH_tip.csv};
        \label{tip_unconstrained}

    \end{axis}
\end{tikzpicture}
    \tikzsetnextfilename{quasistatic_constrained}
    \begin{tikzpicture}

    \begin{axis}[
            width=\plotwidth,
            height=\plotheight,
            xlabel={$\varphi\h_1$},
            ylabel={$\varphi\h_2$},
            xmin=0, xmax=8,
            ymin=-6, ymax=0.5,
            xtick distance=2,
            ytick distance=2,
            minor x tick num=1,
            minor y tick num=1,
            grid=both,
            axis lines=box,
            xtick align=center,
            xtick pos = bottom,
            ytick align=center,
            ytick pos = left,
            legend style={%
                    draw=none,                
                    fill=none,                
                    font=\small,             
                    legend columns=-1,
                    at={(axis description cs:0.5,1.01)}, 
                    anchor=south,
                },
        ]

        \node[anchor= west] at (axis cs:6.3,0) {$t=0$};
        \node[anchor=north west] at (axis cs:4.9,-3.4) {$t=1$};

        \addplot+[
            color3!30,
            line width=\thickplotlinewidth,
        ]
        table[col sep=comma,
                x=x,
                y=y] {constrained_analytical_0.csv};
        \label{legend_analytical}

        \addplot+[
            color3!30,
            line width=\thickplotlinewidth,
            forget plot,
        ]
        table[col sep=comma,
                x=x,
                y=y] {constrained_analytical_1.csv};
        \addplot+[
            color3!30,
            line width=\thickplotlinewidth,
            forget plot,
        ]
        table[col sep=comma,
                x=x,
                y=y] {constrained_analytical_2.csv};
        \addplot+[
            color3!30,
            line width=\thickplotlinewidth,
            forget plot,
        ]
        table[col sep=comma,
                x=x,
                y=y] {constrained_analytical_4.csv};
        \addplot+[
            color3!30,
            line width=\thickplotlinewidth,
            forget plot,
        ]
        table[col sep=comma,
                x=x,
                y=y] {constrained_analytical_10.csv};

        \addplot+[
            color3,
            line width=\plotlinewidth,
            densely dashed,
        ]
        table[col sep=comma,
                x=x_0,
                y=y_0] {ex03_quasistatic_cantilever_constrained_results.csv};
        \label{legend_ours_constrained}

        \addplot+[
            color3,
            line width=\plotlinewidth,
            densely dashed,
            forget plot
        ]
        table[col sep=comma,
                x=x_1,
                y=y_1] {ex03_quasistatic_cantilever_constrained_results.csv};

        \addplot+[
            color3,
            line width=\plotlinewidth,
            densely dashed,
            forget plot
        ]
        table[col sep=comma,
                x=x_2,
                y=y_2] {ex03_quasistatic_cantilever_constrained_results.csv};

        \addplot+[
            color3,
            line width=\plotlinewidth,
            densely dashed,
            forget plot
        ]
        table[col sep=comma,
                x=x_4,
                y=y_4] {ex03_quasistatic_cantilever_constrained_results.csv};

        \addplot+[
            color3,
            line width=\plotlinewidth,
            densely dashed,
            forget plot
        ]
        table[col sep=comma,
                x=x_10,
                y=y_10] {ex03_quasistatic_cantilever_constrained_results.csv};

        \addplot+[
            color3!30,
            line width=\plotlinewidth,
            dash dot,
        ]
        table[col sep=comma,
                x=x,
                y=y] {constrained_analytical_tip.csv};
        \label{tip_constrained}

    \end{axis}
\end{tikzpicture}
    \caption{Quasistatic cantilever problem:
        Centerline position $\varphi_i\h = \vec{e}_i\transp\centerline\h$ for $t\n \in \{0, 0.1, 0.2, 0.4, 1\}$. Good agreement between our results (\plotref{legend_ours_unconstrained}, \plotref{legend_ours_constrained}) and the ones reported in \cite{herrmann_2025_mixed} (\plotref{legend_HE25}) for the unconstrained setting and the analytical solution (\plotref{legend_analytical}) for the constrained setting. The tip centerline moves along the dash-dotted lines (\plotref{tip_unconstrained},\plotref{tip_constrained}).}
    \label{fig_cantilever_quasistatic}
\end{figure}

This numerical experiment displays the usage of the presented formulation in a quasistatic use case and has been taken from \cite{herrmann_2025_mixed,harsch_2021_finite}. We show results for an initally straight and clamped beam with length $L=2\pi$, whose initial configuration aligns with the $\vec{e}_1$-axis. Further parameters are shown in Tab.~\ref{tab_cantilever_quasistatic}. We consider a load case with point force and moment at the free end at $s=L$ such that $ \vec{n}(s=L,t) = \vec{\nu}(t)$ and $\vec{m}(s=L,t) = \vec{\mu}(t)$
with
\begin{table}[t]
    \centering
    \caption{Quasistatic cantilever problem: Physical and numerical parameters.}
    \begin{tabular}{l l l l l l l l l l l l}
        \toprule
        $h$             & $\finaltime$ & $\nel$ & $\epsilon$ & $L$    & $\rho A$ & $M_\rho^{11} = M_\rho^{22}$ & $k_{\mathrm{b1}} = k_{\mathrm{b2}}$ & $k_{\mathrm{t}}$ & $k_{\mathrm{s1}}=k_{\mathrm{s2}}$ & $
        k_{\mathrm{e}}$ & $P$                                                                                                                                                                                                                  \\
        \midrule
        $10^{-2}$       & $1$          & $8$    & $10^{-12}$ & $2\pi$ & $ 0 $    & $ 0 $                       & $ 2 $                               & $0.5$            & $1$                               & $5$ & $10 k_{\mathrm{b1}} / L^2$ \\
        \bottomrule
    \end{tabular}
    \label{tab_cantilever_quasistatic}
\end{table}
\begin{equation}
    \vec{\nu}(t) = t \begin{bmatrix}
        0  \\
        -P \\
        0
    \end{bmatrix}, \quad
    \vec{\mu}(t) = t \begin{bmatrix}
        0 \\
        0 \\
        2.5 P
    \end{bmatrix}, \qquad \text{where} \quad t \in [0,1].
\end{equation}
Time $t$ is regarded in this context as a load factor, see \cite{kinon_2025_energymomentumconsistent} for further explanations. Correspondingly, the inertial parameters are set to zero to achieve quasistatic simulations, which only take place in the $\vec{e}_1$-$\vec{e}_2$-plane. We distinguish two scenarios, one shear-deformable and extensible (\lq\lq unconstrained\rq\rq) and one shear-rigid and inextensible (\lq\lq constrained\rq\rq). In the latter one, an analytical solution for the Euler's elastica can be given \cite{harsch_2021_finite}. The resulting equilibria for both constrained and unconstrained settings are shown in Fig.~\ref{fig_cantilever_quasistatic} and compared to reference results. These confirm not only the correct implementation but also the absence of locking effects.
%

\subsection{Dynamic maneuver of a soft robotic arm}
In this last example, we simulate a soft robotic arm, that is clamped at its bottom and is actuated with pneumatic chambers, see Sec.~\ref{sec_pneumatic_tendon}. The cantilever beam has a circular cross section with diameter $d=0.03$ and we assume the material parameters corresponding to the silicone DragonSkin\texttrademark \,30, i.e., Young's modulus $E=6 \cdot 10^{5}$ and shear modulus $G = 2 \cdot 10^{5}$ \cite{eugster_2022_soft}. This corresponds to a Poisson ratio of $\nu = 0.35$. Furthermore, we have cross-sectional moments of area $I_1 = I_2 = 1/2 I_T = \pi d^4/64$. Further details are comprised in Tab.~\ref{tab_soft_arm}. The beam is actuated by means of three pneumatic chambers with index $k \in \{1, 2, 3\}$, whose lines of centroids are parallel to the beam's centerline with a distance of $r_k=6.5\cdot 10^{-2}$, for all $k$. With respect to the global $\vec{e}_1$-axis these chambers are located initially at angles of $\alpha_k \in \{ \frac{\pi}{6}, \frac{5\pi}{6}, \frac{9\pi}{6}\}$, yielding material coordinates in the cross-section frame $(\rho_k^1,\rho_k^2)= (r_k \cos\alpha_k,r_k\sin\alpha_k)$
for $k=1,2,3$. We prescribe pneumatic forces
\begin{table}[b]
    \centering
    \caption{Dynamic maneuver of a soft robotic arm: Physical and numerical parameters.}
    \begin{tabular}{l l l l l l l l l l l l}
        \toprule
        $h$    & $\finaltime$ & $\nel$ & $\epsilon$ & $L$      & $\rho$   & $A$           & $M_\rho^{11} = M_\rho^{22}$ & $k_{\mathrm{b1}} = k_{\mathrm{b2}}$ & $k_{\mathrm{t}}$ & $k_{\mathrm{s1}}=k_{\mathrm{s2}}$ & $k_{\mathrm{e}}$ \\
        \midrule
        $0.05$ & $4$          & $10$   & $10^{-11}$ & $0.1755$ & $ 1080 $ & $\pi d^2 / 4$ & $ \rho I $                  & $ EI $                              & $G I_T$          & $GA$                              & $EA$             \\
        \bottomrule
    \end{tabular}
    \label{tab_soft_arm}
\end{table}
\begin{equation} \label{pressure_forces}
    p_k (t) A_k
    =\frac{1}{2} f(t) [1 + \cos (\phi(t)-\alpha_k)]
\end{equation}
with the time-dependent functions
\begin{equation} \label{input_functions}
    f(t)=
    \begin{cases}
        \frac{1}{2}f_{\max}\Bigl(1-\cos \bigl(\frac{\pi t}{t_1}\bigr)\Bigr),
         & 0\le t \le t_1, \\
        f_{\max},
         & t_1<t\le t_2,   \\
        \frac{1}{2}f_{\max}\Bigl(1+\cos \bigl(\pi\frac{t-t_2}{T-t_2} \bigr)\Bigr),
         & t_2<t\le T,
    \end{cases}, \qquad \phi(t)=
    \begin{cases}
        0,
         & 0\le t \le t_1, \\
        \pi \Bigl(1-\cos \bigl(\pi \frac{t-t_1}{t_2-t_1}\bigr)\Bigr),
         & t_1<t\le t_2,   \\
        2\pi,
         & t_2<t\le T,
    \end{cases}
\end{equation}
\begin{figure}[tb]
    \setlength{\plotwidth}{0.3\textwidth}
    \setlength{\plotheight}{1\plotwidth}
    \centering
    \tikzsetnextfilename{actuation_input_combined}
    \begin{tikzpicture}

    \begin{axis}[
            width=\plotwidth,
            height=\plotheight,
            xlabel={$t$},
            ylabel={$\phi(t)$},
            xmin=0, xmax=4,
            ymin=-pi/4, ymax=9*pi/4,
            xtick distance=1,
            ytick={0, pi/2, pi, 3*pi/2, 2*pi},
            yticklabels={$0$, , $\pi$, , $2\pi$},
            minor x tick num=1,
            grid=both,
            axis lines=box,
            ytick align=center,
            ytick pos=left,
            ylabel style={color=color4},
            yticklabel style={color=color4},
        ]

        \addplot+[
            color4,
            line width=\plotlinewidth,
        ]
        table[
                col sep=comma,
                x=t,
                y=angle
            ] {ex04_pressure_circle_input.csv};

    \end{axis}

    \begin{axis}[
            width=\plotwidth,
            height=\plotheight,
            xmin=0, xmax=4,
            ymin=-56.25, ymax=6.25,
            minor y tick num=0,
            ytick distance=50,
            axis x line=none,          
            axis y line*=right,        
            ylabel={$f(t)$},
            ytick align=center,
            grid=none,                 
            ylabel style={color=color5,yshift=1em},
            yticklabel style={color=color5},
        ]

        \addplot+[
            color5,
            line width=\plotlinewidth,
        ]
        table[
                col sep=comma,
                x=t,
                y=amplitude
            ] {ex04_pressure_circle_input.csv};

    \end{axis}

\end{tikzpicture}
    \tikzsetnextfilename{actuation_input_forces}
    \begin{tikzpicture}

    \begin{axis}[
            width=\plotwidth,
            height=\plotheight,
            xlabel={$t$},
            ylabel={$p_k(t)A_k$},
            xmin=0, xmax=4,
            ymin=-56.25, ymax=6.25,
            xtick distance=1,
            ytick distance=50,
            minor x tick num=1,
            minor y tick num=3,
            grid=both,
            axis lines=box,
            xtick align=center,
            xtick pos = bottom,
            ytick align=center,
            ytick pos = left,
            ylabel style={yshift=-1em},
            legend style={%
                    draw=none,                
                    fill=none,                
                    font=\small,             
                    legend columns=-1,
                    at={(axis description cs:0.5,1.01)}, 
                    anchor=south,
                },
        ]

        \addplot+[
            line width=\plotlinewidth,
        ]
        table[col sep=comma,
                x=t,
                y=p1] {ex04_pressure_circle_input.csv};
        \label{legend_pressure_1}

        \addplot+[
            line width=\plotlinewidth,
        ]
        table[col sep=comma,
                x=t,
                y=p2] {ex04_pressure_circle_input.csv};
        \label{legend_pressure_2}

        \addplot+[
            line width=\plotlinewidth,
        ]
        table[col sep=comma,
                x=t,
                y=p3] {ex04_pressure_circle_input.csv};
        \label{legend_pressure_3}

    \end{axis}
\end{tikzpicture}
    \tikzsetnextfilename{actuation_time_histories}
    \begin{tikzpicture}

    \begin{axis}[
            width=\plotwidth,
            height=\plotheight,
            xlabel={$t\n$},
            ylabel={$\varphi_i^{\mathrm{tip}}$},
            xmin=0, xmax=4,
            ymin=-0.12, ymax=0.22,
            xtick distance=1,
            ytick distance=0.1,
            minor x tick num=1,
            minor y tick num=1,
            grid=both,
            axis lines=box,
            xtick align=center,
            xtick pos = bottom,
            ytick align=center,
            ytick pos = left,
            ylabel style={yshift=-1em},
            legend style={%
                    draw=none,                
                    fill=none,                
                    font=\small,             
                    legend columns=-1,
                    at={(axis description cs:0.5,1.01)}, 
                    anchor=south,
                },
        ]

        \addplot+[
            color6,
            line width=\plotlinewidth,
        ]
        table[col sep=comma,
                x=time,
                y=x_tip] {ex04_pressure_circle_results.csv};
        \label{legend_xtip_soft}

        \addplot+[
            line width=\plotlinewidth,
            color7,
        ]
        table[col sep=comma,
                x=time,
                y=y_tip] {ex04_pressure_circle_results.csv};
        \label{legend_ytip_soft}

        \addplot+[
            color8,
            line width=\plotlinewidth,
        ]
        table[col sep=comma,
                x=time,
                y=z_tip] {ex04_pressure_circle_results.csv};
        \label{legend_ztip_soft}

    \end{axis}
\end{tikzpicture}
    \caption{Dynamic maneuver of a soft robotic arm: Phase angle function, amplitude function and corresponding pneumatic forces for $k \in \{\textcolor{color1}{1}$\,\plotref{legend_pressure_1}, $\textcolor{color2}{2}$\,\plotref{legend_pressure_2}, $\textcolor{color3}{3}$\,\plotref{legend_pressure_3}$\}$. The tip position components $\varphi_i^{\mathrm{tip}}= \vec{e}_i\transp\vec{\varphi}^{\mathrm{tip}}$ with $i \in \{$
        \textcolor{color6}{$1$}\,\plotref{legend_xtip_soft}, \textcolor{color7}{$2$}\,\plotref{legend_ytip_soft}, \textcolor{color8}{$3$}\,\plotref{legend_ztip_soft}$\}$.}
    \label{fig_actuation_input}
\end{figure}
\begin{figure}[t]
    \setlength{\plotwidth}{0.4\textwidth}
    \setlength{\plotheight}{0.8\plotwidth}
    \centering
    \tikzsetnextfilename{actuation_snapshots}
    \input{snapshots.tex}
    \setlength{\plotwidth}{0.35\textwidth}
    \setlength{\plotheight}{1\plotwidth}
    \tikzsetnextfilename{actuation_path}
    \begin{tikzpicture}

    \begin{axis}[
            width=\plotwidth,
            height=\plotheight,
            xlabel={$\varphi_2^{\mathrm{tip}}$},
            ylabel={$\varphi_1^{\mathrm{tip}}$},
            xmin=-0.1, xmax=0.1,
            ymin=-0.1, ymax=0.1,
            xtick distance=0.1,
            ytick distance=0.1,
            minor x tick num=1,
            minor y tick num=1,
            grid=both,
            axis lines=box,
            xtick align=center,
            xtick pos = bottom,
            ytick align=center,
            ytick pos = left,
            legend style={%
                    draw=none,                
                    fill=none,                
                    font=\small,             
                    legend columns=-1,
                    at={(axis description cs:0.5,1.01)}, 
                    anchor=south,
                },
        ]
        \addplot+[
            color3,
            line width=\plotlinewidth,
        ]
        table[col sep=comma,
                x=x_tip,
                y=y_tip] {ex04_pressure_circle_results.csv};

        \addplot+[
            color2,
            line width=\plotlinewidth,
        ]
        table[col sep=comma,
                x=x_tip,
                y=y_tip] {ex04_pressure_heart_results.csv};
        \label{legend_heart_path};

    \end{axis}
\end{tikzpicture}
    \caption{Dynamic maneuver of a soft robotic arm: Snapshots of the centerline with azimuth and elevation perspective angles $(55,15)$ and the colormap for time $t\n$ with $\redgradrule \in [0.5, 3.5]$. The tip position $\centerline^{\mathrm{tip}}$ follows a circular path (\plotref{legend_soft_tip}). For the modified input functions \eqref{input_functions_heart} the tip follows a heart-shaped path (\plotref{legend_heart_path}).}
    \label{fig_actuation_motion}
\end{figure}%

that describe the amplitude and the phase angle, respectively. Therein,
we set $f_{\max}=-50$, $t_1=0.5$, $t_2=3.5$ and $T=4$.
These functions are visualized in Fig.~\ref{fig_actuation_input}. With a straight initial configuration along the $\vec{e}_3$-axis, the soft robotic arm first undergoes a bending motion in the $\vec{e}_1$-$\vec{e}_3$-plane until $t=t_1$. This is followed by a rotational motion, in which the tip position $\vec{\varphi}^{\mathrm{tip}} := \centerline\h(s=L,t=t\n)$ moves approximately on a constant height, see right depiction in Fig.~\ref{fig_actuation_motion}. At $t=t_2$ the arm has completed one rotation and bends back until $t=T$. The snapshots shown in Fig.~\ref{fig_actuation_motion} only track the circular motion phase for the sake of clarity.

This approach can be further exploited to let the tip follow a more complex path in the $\vec{e}_1$-$\vec{e}_2$ plane, e.g., a heart-shaped one. To this end, for $\vartheta  \in [0, 2\pi]$ a heart-shaped graph has the coordinates $x(\vartheta) = \sin^3(\vartheta)$ and $y(\vartheta) = \cos(\vartheta) - \cos(2\vartheta)$. The corresponding distance to the origin is given by $r(\vartheta) = \sqrt{x(\vartheta)^2 + y(\vartheta)^2}$ and can be normalized as $\tilde r(\vartheta) = r(\vartheta)/ \max_{\vartheta \in [0,2\pi]} r(\vartheta)$. With $\vartheta = 2\pi\frac{t}{T}$ and $\tilde\phi(\vartheta) = \operatorname{atan2}\left(y(\vartheta), x(\vartheta)\right)$, we can eventually specify the amplitude and phase angle functions in \eqref{pressure_forces} as
\begin{equation} \label{input_functions_heart}
    f(t)     = f_{\max} \tilde r\left(2\pi\frac{t}{T}\right), \qquad
    \phi(t) = \tilde\phi\left(2\pi\frac{t}{T}\right).
\end{equation}
The tip of the soft robotic arm then exhibits a heart shaped trajectory, see Fig.~\ref{fig_actuation_motion}.
More complex deformation patterns can be achieved by including tendons and pneumatic chambers, whose line of centroids do not lie parallel to the centerline of the beam, or by dividing the beam into multiple segments with individual actuation input forces.

\section{Conclusion \& outlook}
\label{sec_conclusion}
In this work, we have presented a mixed formulation and structure-preserving discretization procedure for the dynamics of spatial Cosserat rods.
The proposed framework is objective, locking-free and enjoys many favorable properties due to its port-Hamiltonian structure. Through a compliance equation also shear-rigid and inextensible beam formulations can be seamlessly simulated within the present framework.

The employed PH formalism naturally leads to a structure-preserving mixed FE method. In this setting, the resulting PH differential-algebraic equations admit a particularly simple time-integration strategy based on pure midpoint approximations. The scheme is energy- and momentum-consistent, provides second-order accuracy and is remarkably robust.
A further consequence is the accurate representation of crucial kinematic relations. Both the compliance equation and the orthonormality constraints have been differentiated with respect to time to reveal the PH structure. The proposed time integration scheme exactly integrates those relations.

For convenience, we have restricted ourselves to initially straight and stress-free beam problems with constant cross-section geometry and material parameters. An extension of the proposed framework to more general settings would merely impose more tedious algebra.
While the presented formulation features the full geometrical nonlinearity, i.e., allowing for
large displacements and rotations in three dimensions,
we have always assumed linear constitutive relations.
The present framework can be extended to general hyperelasticity along the lines of \cite{kinon_2023_porthamiltonian,kinon_2024_generalized}, requiring the inclusion of a discrete-gradient time integration approach, which goes well with the PH structure \cite{kinon_2023_discrete,kinon_2023_porthamiltonian,kinon_2026_discrete,goren-sumer_2008_gradient}.

Control-oriented actuation mechanisms and visco-elastic material behavior have been included into the PH framework in a straightforward manner, making it ideally suited for interesting application fields
\cite{albu-schaffer_2008_soft,dorlich_2018_flexible}. One future research direction would be the inclusion of our model into novel control strategies, like control by interconnection or energy shaping methods \cite{caasenbrood_2022_energyshaping,ortega_2002_interconnection,mehrmann_2023_control}.

By employing a singularity-free director parametrization, the pitfalls of numerical methods for finite rotations have been circumvented.
Despite its conceptual simplicity, the director formulation however does not reproduce the correct
rotation
inside the finite element, since the orthonormality condition on the directors is violated inside the element in general. In the future, projected unit quaternions \cite{lang_2012_numerical,harsch_2023_nonunit,wasmer_2024_projectionbased,debeurre_2025_quaternionbased} could be used to obtain a manifold-consistent discrete PH formulation for Cosserat rod dynamics, coming at the price of additional nonlinearities and configuration-dependent mass matrices.

\begin{acks}
PLK gratefully acknowledges funding by the Research Travel Grant of the Karlsruhe House of Young Scientists (KHYS).
We are thankful for fruitful discussions with Marco Herrmann (TU Eindhoven), Tianxiang Dai (University of Stuttgart) and Maximilian Herrmann (TU Munich).
\end{acks}
\begin{authcontrib}
\begin{table*}[h]
    \centering
    \vspace{-2em}
    \begin{tabular}{l p{14cm}}
        \textbf{PLK} & Conceptualization,
        methodology,
        software,
        validation,
        formal analysis,
        investigation,
        data curation,
        writing--original draft preparation,
        visualization,
        funding acquisition.              \\
        \textbf{SRE} &
        Methodology,
        validation,
        formal analysis,
        resources,
        writing--review and editing,
        supervision.                      \\
        \textbf{PB}  & Conceptualization,
        methodology,
        validation,
        formal analysis,
        resources,
        writing--review and editing,
        supervision,
        project administration,
        funding acquisition.
    \end{tabular}
\end{table*}
\vspace{-1.5em}
\end{authcontrib}
\begin{dci}
    The authors declare no known conflict of interest.
\end{dci}
\begin{code}
    The simulation results are openly available at \cite{kinon_2025_repository} and are archived under \cite{kinon_2025_18007058}.
\end{code}

\sloppy
\printbibliography

\appendix
\section*{Appendix}
\section{Director velocities and virtual director displacements}\label{app_math}

This section introduces an important relation between virtual rotations and virtual director displacements. Additionally, convenient matrix identities are provided for later use.

The directors $\{\vec{d}_i\}_{i=1}^3$ introduced in \Cref{sec_kinematic} satisfy the orthogonality condition \eqref{eq:ortho_constraint} to be a valid parametrization of finite rotations. Correspondingly, we have $\vec{d}_i = \vec{R} \vec{e}_i$ for each $i \in \{1, 2, 3\}$ such that the time derivative and the variation of a director is given by \cite{harsch_2021_finite}
\begin{align} \label{virtual_and_velocity_director}
    \delta \vec{d}_i = \delta \vec{\phi} \times \vec{d}_i, \qquad \qquad  \qquad \sdot{\vec{d}}_i = \vec{\omega} \times \vec{d}_i ,
\end{align}
where the spatial angular velocity $\vec{\omega}$ and the virtual rotation $\delta \vec{\phi}$ are defined by means of their skew-symmetric matrices
$\delta \tilde{\vec{\phi}} = \delta \vec{R} \vec{R}\transp$ and $\tilde{\vec{\omega}} = \dot{\vec{R}}\vec{R}\transp$. To invert relation \eqref{virtual_and_velocity_director}, we make use of Graßmann's identity (or also referred to as Lagrange's formula) for the vector triple product
\begin{align} \label{triple_cross_product}
    \vec{a} \times (\vec{b} \times \vec{c}) = \vec{b} (\vec{a}\transp \vec{c}) - \vec{c}(\vec{a}\transp \vec{b}) \qquad \Leftrightarrow \quad \tilde{\vec{a}}\tilde{\vec{b}} = \vec{b} \vec{a}\transp - (\vec{a}\transp \vec{b})\vec{I}
\end{align}
to show that $\vec{d}_i \times \delta \vec{d}_i = 2 \delta \vec{\phi}$, \cite{betsch_2013_consistent,harsch_2021_finite}. Consequently, the virtual rotation $\delta \vec{\phi}$ and angular velocity $\vec{\omega}$ can be expressed as
\begin{align}
    \delta \vec{\phi} & = \frac{1}{2} \vec{d}_i \times \delta \vec{d}_i = \frac{1}{2}  \tilde{\vec{d}}_i \delta \vec{d}_i , \qquad \qquad  \qquad \vec{\omega} = \frac{1}{2} \vec{d}_i \times \sdot{\vec{d}}_i = \frac{1}{2}  \tilde{\vec{d}}_i \dot{\vec{d}}_i.
\end{align}
These identities can be rewritten as
\begin{align} \label{eq_tmp2}
    \delta \vec{\phi} = \vec{T}(\vec{d})\transp \delta \vec{d}, \qquad  \qquad \vec{\omega} = \vec{T}(\vec{d})\transp \dot{\vec{d}} ,
\end{align}
where we made use of $\vec{d}=(\vec{d}_1,\vec{d}_2,\vec{d}_3) \in \mathbb{R}^9$ and the director-dependent matrix $ \vec{T}(\vec{d}) \in \mathbb{R}^{9 \times 3}$ introduced in \eqref{eq_T}.
This matrix moreover appears in the identities
\begin{equation} \label{eq_T_identities}
    \begin{aligned}
        \vec{T}(\vec{d})\transp \vec{d}    & =\vec{0}, & \quad        \vec{L}(\vec{d}_{,s})\vec{T}(\vec{d})                & = \vec{L}(\vec{d})\vec{T}(\vec{d}_{,s}),                  & \quad \vec{T}(\vec{d})\transp \vec{T}(\vec{d}) & = -\frac{1}{2}\vec{I}, \\
        \vec{G}_d(\vec{d})\vec{T}(\vec{d}) & =\vec{0}, & \quad \vec{T}(\vec{d})\transp \vec{J}_{\vec{N}}(\centerline_{,s}) & = -\frac{1}{2} \tilde{\centerline}_{,s}\vec{R}(\vec{d}) , &
    \end{aligned}
\end{equation}
which will be used at a later stage. The matrices $\vec{R}$,  $\vec{L}$, $\vec{J}_{\vec{N}}$ and $\vec{G}_d$, are defined in \eqref{eq_rotmat}, \eqref{eq_JN_L} and \eqref{constraint_jacobian}, respectively.

\section{Principle of virtual work} \label{virtual_work_appendix}

\subsection{Displacement-based principle of virtual work}
The partial differential equations governing the dynamics of the Cosserat rod can be derived from the principle of virtual work \cite{betsch_2002_frameindifferent,harsch_2021_finite}, in which the equality
\begin{equation} \label{eq:virtual_work_principle}
    \delta W^\mathrm{dyn} + \delta W^\mathrm{int} + \delta W^\mathrm{ext} + \delta W^{\mathrm{c}} = 0
\end{equation}
holds for all admissible and arbitrary virtual displacements and
where the four terms denote the virtual work contributions of dynamical, internal, external and constraint forces, respectively. The individual terms are discussed to obtain first the formulation \eqref{EoM_Antman_both}, which can be found in standard literature \cite{antman_2005_nonlinear,simo_1985_finitea}, and then the director-based formulation \cite{romero_2002_objective,betsch_2002_frameindifferent}, which we here formulated by means of operator matrices in \eqref{eq:balance_ang_momentum_1}.

\paragraph{Formulation in terms of virtual rotations}
Formulation \eqref{EoM_Antman_both} can be obtained by expressing the virtual work contributions in \eqref{eq:virtual_work_principle} in terms of virtual centerline displacements $\delta \centerline$ and virtual rotation vectors $\delta \vec{\phi}$, see relations in \Cref{app_math}.
The dynamical virtual work $\delta W^\mathrm{dyn}$ reads
\begin{equation} \label{W_dyn}
    \delta W^\mathrm{dyn} = \int_\domainspace \left( \delta\centerline\transp \rho A \ddot{\centerline} + \delta \vec{\phi}\transp (\vec{I}_\rho \dot{\vec{\omega}} + \vec{\omega} \times \vec{I}_\rho \vec{\omega} ) \right) \d s ,
\end{equation}
see \cite{eugster_2020_variational}. The spatial inertia tensor  $\vec{I}_\rho$ is related to its material counterpart $ \vec{J}_\rho = M_\rho^{\alpha \alpha} \vec{I} - M_\rho^{\alpha \beta} \vec{e}_{\alpha} \otimes \vec{e}_\beta$
via $\vec{I}_\rho = \rotmat \vec{J}_\rho \rotmat\transp$. $M_\rho^{\alpha \beta}$ for $\alpha, \beta = 1,2$ are the uniformly distributed second moments of mass of the cross section with respect to the director basis.
The internal virtual work $\delta W^\mathrm{int}$ emanates from the variation of the elastic potential \eqref{quadratic_stored_energy}
\begin{equation} \label{W_int}
    \delta W^\mathrm{int} = \delta \int_\domainspace \breve{W}(\vec{\Gamma},\vec{K}) \d s = \int_\domainspace \left( \vec{N}\transp \delta \vec{\Gamma} + \vec{M}\transp \delta\vec{K} \right) \d s,
\end{equation}
in which $\vec{N}=\partial_{\vec{\Gamma}}\breve{W}$ and $\vec{M}=\partial_{\vec{K}}\breve{W}$ are the dual quantities to the variation of the strain measures, see \cite{betsch_2002_frameindifferent,eugster_2020_variational}, given by
\begin{align} \label{virtual_strains}
    \delta \vec{\Gamma} = \vec{R}\transp (\delta \centerline_{,s} - \delta \vec{\phi} \times \centerline_{,s} ) \qquad \qquad  \qquad
    \delta \vec{K}  = \vec{R}\transp \delta \vec{\phi}_{,s} .
\end{align}
Substituting \eqref{virtual_strains} into \eqref{W_int} and using \eqref{trafo_spat_forces}, the internal virtual work takes the alternative form
\begin{equation}
    \delta W^\mathrm{int} = \int_\domainspace \left( \vec{n}\transp \left( \delta \centerline_{,s} - \delta \vec{\phi} \times \centerline_{,s} \right) + \vec{m}\transp \delta \vec{\phi}_{,s} \right) \d s .
\end{equation}
The external virtual work due to distributed forces and torques as well as Neumann boundary terms can be written as
\begin{equation} \label{W_ext}
    \delta W^\mathrm{ext} = -\int_\domainspace \left( \delta \centerline \transp \bar{\vec{n}} + \delta \vec{\phi} \transp \bar{\vec{m}} \right) \d s - \left[ \delta \centerline \transp \vec{n}_{\mathrm{N}} + \delta \vec{\phi} \transp \vec{m}_{\mathrm{N}} \right]_{\partial\domainspace_{\mathrm{N}}} .
\end{equation}
There are no constraints appearing in this formulation, hence $\delta W^{\mathrm{c}}=0$. Eventually, one proceeds by collecting all terms connected to the virtual centerline displacements and virtual rotations in \eqref{eq:virtual_work_principle}. After integration by parts and stating the arbitrariness of the virtual quantities one obtains the local balance equations \eqref{EoM_Antman_both} together with the corresponding Neumann boundary conditions.

\paragraph{Formulation in terms of directors}
Again, we begin with stating the terms appearing in the principle of virtual work \eqref{eq:virtual_work_principle}.
Accounting for the director-kinematics, the dynamical virtual work \eqref{W_dyn} has the alternative representation \cite{betsch_2002_frameindifferent,eugster_2020_variational}
\begin{equation} \label{W_dyn_director}
    \delta W^\mathrm{dyn} = \int_\domainspace \left( \delta\centerline\transp \rho A \ddot{\centerline} + \delta \vec{d}_\alpha\transp M_\rho^{\alpha \beta} \ddot{\vec{d}}_\beta  \right) \d s = \int_\domainspace \left( \delta\centerline\transp \rho A \ddot{\centerline} + \delta \vec{d}\transp \vec{M}_\rho \ddot{\vec{d}} \, \right) \d s ,
\end{equation}
where we have used the director-related mass matrix $\vec{M}_\rho$, introduced in \Cref{sec_eom}.
For the internal virtual work \eqref{W_int}, we can express the variation of the strain measures \eqref{virtual_strains} alternatively as $\delta \vec{\Gamma} = (\delta \vec{d}_k \transp \centerline_{,s} + \vec{d}_k \transp \delta \centerline_{,s})\vec{e}_k $ and $\delta \vec{K}    = \frac{1}{2} \epsilon_{ijk} (\delta \vec{d}_k \transp \vec{d}_{j,s} + \vec{d}_k \transp \delta \vec{d}_{j,s}) \vec{e}_i$,
see definitions \eqref{eq:strain_measures}.
This can be recast as
\begin{equation}
    \delta\vec{\Gamma} = \vec{J}_{\vec{N}}(\centerline_{,s})\transp \delta\vec{d} + \vec{R}(\vec{d})\transp \delta \centerline_{,s}, \qquad \qquad \delta \vec{K} = \vec{L}(\vec{d}_{,s})\delta \vec{d} - \vec{L}(\vec{d})\delta \vec{d}_{,s},
\end{equation}
conforming with \eqref{eq:strain_measures_short}.
Substituting these variations of the strain measures into \eqref{virtual_strains} and performing integration by parts yields
\begin{equation}
    \begin{aligned} \label{W_int_dir}
        \delta W^\mathrm{int} & = \int_\domainspace \left( \delta \vec{d} \transp \vec{J}_{\vec{N}}(\centerline_{,s}) \vec{N} - \delta \centerline\transp \partial_s(\vec{R}(\vec{d})\vec{N}) + \delta \vec{d}\transp \operator{J}_{\vec{M}}(\vec{d})\vec{M} \right) \d s + \left[\delta \centerline\transp \vec{n} + \delta\vec{d}\transp \vec{T}(\vec{d})\vec{m} \right]_{\partial \domainspace} .
    \end{aligned}
\end{equation}
Here, we have made use of \eqref{trafo_spat_forces}, \eqref{virtual_and_velocity_director}, \eqref{eq_tmp2} and the relation $\vec{T}(\vec{d})= -\vec{L}(\vec{d})\transp \vec{R}(\vec{d})\transp$,
which can be shown by using \eqref{triple_cross_product}. Eventually, the matrix $\vec{J}_{\vec{N}}$ as defined in \eqref{eq_JN_L} and the matrix operator $\operator{J}_{\vec{M}}$ as defined in \eqref{eq_BM_BK} could be identified.
For the external virtual work \eqref{W_ext}, relation \eqref{eq_tmp2} is employed to obtain
\begin{equation} \label{W_ext_dir}
    \delta W^\mathrm{ext} = -\int_\domainspace \left( \delta \centerline \transp \bar{\vec{n}} + \delta \vec{d} \transp \vec{T}(\vec{d}) \bar{\vec{m}} \right) \d s - \left[ \delta \centerline \transp \vec{n}_{\mathrm{N}} + \delta \vec{d} \transp \vec{T}(\vec{d}) \vec{m}_{\mathrm{N}} \right]_{\partial\domainspace_{\mathrm{N}}} .
\end{equation}
And lastly, the virtual work done by the constraint forces related to the orthonormality conditions \eqref{eq:ortho_constraint}, is given by
\begin{equation} \label{W_c}
    \delta W^{\mathrm{c}} = \delta \int_\domainspace \left( \vec{\lambda}\transp \vec{G}(\vec{d}) \right) \d s =  \int_\domainspace \left( \delta \vec{\lambda}\transp \vec{G}(\vec{d}) + \delta\vec{d}\transp \vec{G}_d(\vec{d})\transp \vec{\lambda} \right) \d s ,
\end{equation}
where the constraint gradient matrix $\vec{G}_d = \partial_{\vec{d}}\vec{g}(\vec{d})$ is specified in
\eqref{constraint_jacobian}.
Eventually, one proceeds by
collecting the different terms related to the virtual centerline displacements, virtual director displacements and virtual Lagrange multipliers. Subsequently, accounting for the arbitrariness of the virtual quantities yields the equations \eqref{eq:balance_ang_momentum_1}.

\subsection{Hellinger-Reissner type principle of virtual work} \label{sec_HR}
The weak form of the PH formulation employed in this work \eqref{weak_form_directors} can be traced back to a Hellinger-Reissner type virtual work principle.%
Let us start by modifying the internal virtual work, which lies at the core of Hellinger-Reissner approaches.
As shown in \cite{herrmann_2025_mixed}, the strain energy density function pertaining to a Hellinger-Reissner two-field principle 
is given by
\begin{equation} \label{HR_approach}
    W = W^{\mathrm{HR}}(\vec{\Gamma},\vec{K},\vec{N},\vec{M}) = - W^\star(\vec{N},\vec{M}) + \vec{N}\transp \vec{\Gamma} + \vec{M}\transp \vec{K}
\end{equation}
where the strains are functions of the displacement quantities, see \eqref{eq:strain_measures} and $W^\star$ is the complementary strain energy density \eqref{eq_comp_strain} depending on independent stress quantities. Correspondingly, the virtual internal work appearing in the principle of virtual work \eqref{eq:virtual_work_principle} assumes the alternative form
\begin{align}
    \delta W^\mathrm{int,HR} & = \delta \int_{\domainspace} W^{\mathrm{HR}} \d s  = \int_\domainspace \left( \delta \vec{\Gamma}\transp \vec{N} + \delta\vec{K}\transp \vec{M} + \delta\vec{N}\transp(\vec{\Gamma} - \vec{C}_N \vec{N}) + \delta\vec{M}\transp(\vec{K} - \vec{C}_M \vec{M})  \right) \d s , \label{HR_general}
\end{align}
where we have substituted \eqref{eq:material}. In contrast to \eqref{W_int}, the stress quantities are assumed to be independent from the displacement quantities. Correspondingly, the additional terms in \eqref{HR_general} enforce the compliance relation \eqref{eq:material} for compatibility. This expression is further rewritten similarly to \eqref{W_int_dir} as
\begin{equation} \label{W_int_HR}
    \begin{aligned}
        \delta W^\mathrm{int,HR}
         & = \int_\domainspace \left( \delta \centerline_{,s}\transp \vec{n} + \delta\vec{d}\transp \vec{J}_{\vec{N}}\vec{N} + \vec{M}\transp \operator{J}_{\vec{K}}(\vec{d})\delta\vec{d} \right. \\
         & \qquad \qquad \left.
        + \, \delta\vec{N}\transp(\vec{\Gamma} - \vec{C}_N \vec{N}) + \delta\vec{M}\transp(\vec{K} - \vec{C}_M \vec{M})  \right) \d s ,
    \end{aligned}
\end{equation}
where the differential operator matrix $\operator{J}_{\vec{K}}$ is given in \eqref{eq_BM_BK}.
Additionally, we can modify the virtual work of dynamical forces in the director formulation \eqref{W_dyn_director} to account for independent velocities \eqref{eq:kinematics} by mimicking the Hellinger-Reissner approach \eqref{HR_general} as
\begin{equation} \label{W_dyn_indy}
    \delta W^\mathrm{dyn,HR} = \int_\domainspace \left( \delta\centerline\transp \rho A \dotvphi + \delta \vec{d}\transp \vec{M}_\rho \dot{\vec{v}}_d + \rho A \delta\vphi\transp (\dot{\centerline} - \vphi) + (\vec{M}_\rho\delta\vec{v}_d)\transp(\dot{\vec{d}} - \vec{v}_d) \, \right) \d s ,
\end{equation}
where the kinematic relation \eqref{eq:kinematics} is enforced to ensure compatible independent velocities.
Collecting the terms \eqref{W_ext_dir}, \eqref{W_c}, \eqref{W_int_HR} and \eqref{W_dyn_indy}, and separating the terms corresponding to the different virtual quantities yields the weak form
\begin{subequations}
    \begin{align}
        \inner{\rho A \delta \vphi}{\dot{\centerline}- \vphi}            & =0 , \\
        \inner{\vec{M}_\rho \delta \vec{v}_d}{\sdot{\vec{d}}- \vec{v}_d} & =0 , \\
        \inner{\delta \centerline}{\rho A \dotvphi - \bar{\vec{n}}}
        +  \inner{\partial_s(\delta \centerline)}{ \rotmat\vec{N}}
        - \innerBoundary{\delta \centerline}{\vec{n}_N}                  & =0 , \\
        \inner{\delta \vec{d}}{\vec{M}_\rho \dot{\vec{v}}_d +  \vec{J}_{\vec{N}} \vec{N} + \vec{G}_d\transp \vec{\lambda}
            - \vec{T}\bar{\vec{m}}}
        + \inner{\operator{J}_{\vec{K}}\delta \vec{d}}{ \vec{M}}
        - \innerBoundary{\delta \vec{d}}{\vec{T}\vec{m}_N}               & =0 , \\
        \inner{\delta \vec{N}}{\vec{C}_N \vec{N} -
            \rotmat\transp \centerline_{,s}
        }                                                                & =0 , \\
        \inner{\delta \vec{M}}{\vec{C}_M \vec{M} -
            \vec{L}(\vec{d}_{,s})\vec{d}
        }                                                                & =0 , \\
        \inner{\delta \vec{\lambda}}{\vec{g}(\vec{d})
        }                                                                & =0 .
    \end{align}
\end{subequations}
which can be linked to the weak form \eqref{weak_form_directors} by replacing the last three relations with their time-differentiated version.

\section{Proof: Formally skew-adjoint structure operator} \label{app_skew_adjoint}

Let us here verify that the structure operator $\operator{J}$ contained in the dynamics \eqref{phDAE_compact} is formally skew-adjoint. This is equivalent with showing that \eqref{boundary_terms} holds true and the inner product $\inner{\vec{z}}{\vec{J} \vec{z}}$ only yields boundary terms. For elements of $\operator{J}$ that do not involve spatial derivatives, the formal skew-adjointness boils down to skew-symmetry and no integration by parts is required. The corresponding terms cancel each other out.
To this end, it remains
to show that
\begin{equation} \label{skew_adjoint_proof_1}
    \inner{
        \begin{bmatrix}
            \vec{v} \\ \vec{\sigma}
        \end{bmatrix}
    }{\begin{bmatrix}
            \zeros                & -\operator{J}_{\vec{v}} \\
            \operator{J}_{\sigma} & \zeros
        \end{bmatrix}
        \begin{bmatrix}
            \vec{v} \\ \vec{\sigma}
        \end{bmatrix}
    } = \inner{\vec{\sigma}}{\operator{J}_{\sigma}\vec{v}} - \inner{\vec{v}}{\operator{J}_{\vec{v}} \vec{\sigma}} = \left[\vphi \transp \vec{n} + \vec{\omega} \transp \vec{m} \right]_{\partial \domainspace}
\end{equation}
corresponding to $ \operator{J}_{\vec{v}}^* = \operator{J}_{\vec{\sigma}}$.
This condition can be expressed more detailed as
\begin{equation} \label{eq_skew_adj}
    \begin{aligned}
         & \inner{
            \begin{bmatrix}
                \vec{N} \\ \vec{M}
            \end{bmatrix}}{
            \begin{bmatrix}
                \rotmat(\vec{d})\transp \partial_s \square & \vec{J}_{\vec{N}}(\centerline_{,s})\transp   \\
                \vec{0}                                    & \operator{J}_{\vec{K}}(\vec{d},\vec{d}_{,s})
            \end{bmatrix}
        \begin{bmatrix}
                \vphi \\ \vec{v}_d
            \end{bmatrix}}       \\
         & \qquad \qquad \qquad
        -
        \inner{
            \begin{bmatrix}
                \vphi \\ \vec{v}_d
            \end{bmatrix}}{
            \begin{bmatrix}
                -\partial_s (\rotmat(\vec{d}) \square) & \vec{0}                                      \\
                \vec{J}_{\vec{N}}(\centerline_{,s})    & \operator{J}_{\vec{M}}(\vec{d},\vec{d}_{,s})
            \end{bmatrix}
            \begin{bmatrix}
                \vec{N} \\ \vec{M}
            \end{bmatrix}} = \left[\vphi \transp \vec{n} + \vec{\omega} \transp \vec{m} \right]_{\partial \domainspace} .
    \end{aligned}
\end{equation}
Here, one can identify three terms. The first one vanishes due to skew-symmetry, i.e.,
\begin{equation}
    \int_\domainspace \left( \vec{N}\transp \vec{J}_{\vec{N}}(\centerline_{,s})\transp \vec{v}_d - \vec{v}_d\transp \vec{J}_{\vec{N}}(\centerline_{,s}) \vec{N} \right) \d s = 0 .
\end{equation}
The second term can be analyzed as
\begin{align}
     & \int_\domainspace \left( \vec{N} \transp \rotmat(\vec{d})\transp \vphis +  \vphi \transp \partial_s(\rotmat(\vec{d})\vec{N}) \right) \d s =
    \int_\domainspace \left( \vphi \transp \vec{n}_{,s} + \vec{n} \transp \vphis \right) \d s =
    \left[\vphi \transp \vec{n} \right]_{\partial \domainspace} ,
\end{align}
which is equivalent to $(\partial_s(\rotmat(\vec{d}) \square ))^* = - \rotmat(\vec{d})\transp \partial_s \square $. The third and last term yields
\begin{align}
     & \int_\domainspace \left( \vec{M} \transp \operator{J}_{\vec{K}}\vec{v}_d- \vec{v}_d\transp \operator{J}_{\vec{M}} \vec{M} \right) \d s                                                                               \notag                                                                                \\
     & \qquad \qquad = \int_\domainspace \left( \vec{M} \transp \vec{L}(\vec{d}_{,s}) \vec{v}_d - \vec{M} \transp \vec{L}(\vec{d}) \vec{v}_{d,s}  - \vec{v}_d\transp  \vec{L}(\vec{d}_{,s})\transp \vec{M} - \vec{v}_d\transp \partial_s(\vec{L}(\vec{d})\transp \vec{M})  \right) \d s                           \\
     & \qquad \qquad  = \left[ -\vec{v}_d\transp \vec{L}(\vec{d})\transp \vec{M} \right]_{\partial \domainspace} =  \left[ -\vec{v}_d\transp \vec{L}(\vec{d})\transp \rotmat(\vec{d})\transp \vec{m} \right]_{\partial \domainspace} = \left[ \vec{\omega} \transp \vec{m} \right]_{\partial \domainspace} \notag
\end{align}
and shows that $\operator{J}^*_{\vec{M}} = \operator{J}_{\vec{K}} $. In the last equation, we have used \eqref{eq_tmp2} and the relation $\vec{T}(\vec{d})= -\vec{L}(\vec{d})\transp \vec{R}(\vec{d})\transp$,
which can be shown by using \eqref{triple_cross_product}.
Eventually, the claim \eqref{skew_adjoint_proof_1} follows from adding up the last three relations and comparing them with \eqref{eq_skew_adj}.

\section{Balance of total angular momentum} \label{sec_ang_mom}

The fundamental balance of total angular momentum can be shown by choosing specific test functions in the weak form \eqref{weak_form_directors}. In the present case, the total angular momentum \cite{romero_2002_objective} is given by
\begin{equation} \label{eq_total_ang_mom}
    \vec{l}(\vec{q},\vec{v}) = \int_\domainspace \left( \centerline \times \rho A \vphi + \vec{d}_\alpha \times M_\rho^{\alpha \beta} \vec{v}_\beta \right) \d s = \int_\domainspace \left( \tilde{\centerline} \rho A \vphi - 2 \vec{T}(\vec{d})\transp \vec{M}_\rho \vec{v}_d \right) \d s ,
\end{equation}
where summation over $\alpha = 1,2$ is implied and \eqref{eq_T} has been used.
To show the balance of total angular momentum, we assume pure Neumann boundary conditions such that we are allowed to choose the test functions
\begin{equation}
    \begin{aligned}
        \delta \centerline & = \vec{\xi} \times \centerline = - \tilde{\centerline} \vec{\xi} , &  & \delta \vphi      = \vec{\xi} \times \vphi = - \tilde{\vec{v}}_\varphi \vec{\xi} , \\
        \delta \vec{d}_i   & = \vec{\xi} \times \vec{d}_i = - \tilde{\vec{d}}_i \vec{\xi} ,     &  & \delta \vec{v}_i  = \vec{\xi} \times \vec{v}_i = - \tilde{\vec{v}}_i \vec{\xi} ,
    \end{aligned}
\end{equation}
in \eqref{weak_form_directors}. With constant but arbitrary $\vec{\xi} \in \R^3$, this corresponds to a rigid body rotation.
Using \eqref{eq_T}, we have
\begin{align} \label{eq_test_angmom}
    \delta \vec{d} = -2 \vec{T}(\vec{d}) \vec{\xi} , \qquad  \qquad
    \delta \vec{v}_d = -2 \vec{T}(\vec{v}_d) \vec{\xi} .
\end{align}
Next, we insert these test functions into \eqref{weak_form_directors_1} - \eqref{weak_form_directors_4}. Adding up \eqref{weak_form_directors_3} and \eqref{weak_form_directors_4} and subtracting \eqref{weak_form_directors_1} and \eqref{weak_form_directors_2} yields after some algebraic manipulations
\begin{align} \label{eq_total_angmom_balance}
    \dot{\vec{l}} = \int_\domainspace ( \centerline \times \bar{\vec{n}} + \bar{\vec{m}} ) \d s + \left[ \centerline \times \vec{n}_N + \vec{m}_N \right]_{\partial \domainspace} ,
\end{align}
Here, we have used the properties \eqref{eq_T_identities} and the arbitrariness of the constant $\vec{\xi}$.
Equation~\eqref{eq_total_angmom_balance} states that the rate of change of total angular momentum is equal to the sum of all external moments acting on the beam, both distributed and concentrated at the boundaries. In case of closed systems, the total angular momentum is a conserved quantity.

\begin{remark}
    Note that in a similar fashion one may show the balance of total linear momentum and of total energy. The latter has already been demonstrated in a more appealing way by means of the PH structure, see \eqref{power_balance_cont}.
\end{remark}

\paragraph{Discrete time angular momentum balance.}
It can be shown that the time stepping scheme \eqref{timediscrete_EoM} retains the balance equation for total angular momentum \eqref{eq_total_angmom_balance}.
In the bulk part of this work we have derived the time integration method \eqref{timediscrete_EoM} by means of a sequential space- and time approximation of the continuous weak form \eqref{weak_form_directors}. In fact, the first four lines of \eqref{timediscrete_EoM} can be equivalently rewritten in a discrete weak form given by
\begin{subequations} \label{weak_form_directors_discrete}
    \begin{align}
        \inner{\rho A \delta \vphi\h}{\centerline\h\npe - \centerline\h\n - h (\vphi\h)\npeh}                                                                                                                                    & =0 , \label{weak_form_directors_discrete_1} \\
        \inner{\vec{M}_\rho \delta \vec{v}_d\h}{\vec{d}\h\npe - \vec{d}\h\n - h (\vec{v}_d\h)\npeh}                                                                                                                              & =0 , \label{weak_form_directors_discrete_2} \\
        \inner{\delta \centerline\h}{\rho A (\vphi\h)\npe - \rho A (\vphi\h)\n - h \bar{\vec{n}}\npeh}
        +  \inner{\partial_s(\delta \centerline\h)}{h \rotmat(\vec{d}\h\npeh) \vec{N}\h\npeh}   \qquad                                                                                                                           & \notag                                      \\
        - \innerBoundary{\delta \centerline\h}{ h (\vec{n}_N)\npeh}                                                                                                                                                              & =0 , \label{weak_form_directors_discrete_3} \\
        \inner{\delta \vec{d}\h}{\vec{M}_\rho (\vec{v}_d\h)\npe - \vec{M}_\rho (\vec{v}_d\h)\n + h \vec{J}_{\vec{N}}(\vec{\centerline}\h\npeh) \vec{N}\h\npeh + h \vec{G}_d(\vec{d}\h\npeh)\transp \vec{\lambda}\h\npeh}  \qquad & \notag                                      \\
        - \inner{\delta \vec{d}\h}{h \vec{T}(\vec{d}\h\npeh) \bar{\vec{m}}\npeh}
        + \inner{\operator{J}_{\vec{K}}(\vec{d}\h\npeh)\delta \vec{d}\h}{h \vec{M}\h\npeh} \qquad                                                                                                                                & \notag                                      \\
        - \innerBoundary{\delta \vec{d}\h}{h \vec{T}(\vec{d}\h\npeh)(\vec{m}_N)\npeh}                                                                                                                                            & =0 . \label{weak_form_directors_discrete_4}
    \end{align}
\end{subequations}
Due to the bilinear form of the total angular momentum \eqref{eq_total_ang_mom}, we can apply the findings from Remark~\ref{rem_meanvalue} to the difference from one time step to another, such that
\begin{equation} \label{eq_discrete_difference_angmom}
    \begin{aligned}
        \vec{l}(\vec{q}\h\npe,\vec{v}\h\npe) - \vec{l}(\vec{q}\h\n,\vec{v}\h\n) & = \int_\domainspace \left( \tilde{\centerline}\h\npeh \rho A ((\vphi\h)\npe - (\vphi\h)\n) - 2 \vec{T}(\vec{d}\h\npeh)\transp \vec{M}_\rho ((\vec{v}_d)\h\npe - (\vec{v}_d)\h\n) \right. \\
                                                                                & \qquad \quad \left. + \, (\tilde{\centerline}\h\npe - \tilde{\centerline}\h\n) \rho A (\vphi\h)\npeh - 2 \vec{T}(\vec{d}\h\npe - \vec{d}\h\n)\transp \vec{M}_\rho (\vec{v}_d)\h\npeh
        \right) \d s .
    \end{aligned}
\end{equation}
Similar to above, we choose test functions
\begin{equation}
    \delta \centerline = - \tilde{\centerline}\h\npeh \vec{\xi}, \quad \delta \vphi = - (\tilde{\vec{v}}_\varphi\h)\npeh \vec{\xi}, \quad \delta \vec{d} = -2 \vec{T}(\vec{d}\h\npeh) \vec{\xi} , \quad
    \delta \vec{v}_d = -2 \vec{T}((\vec{v}_d\h)\npeh) \vec{\xi}
\end{equation}
in the discrete-time weak forms \eqref{weak_form_directors_discrete_1}-\eqref{weak_form_directors_discrete_4}. As in the continuous case, adding up \eqref{weak_form_directors_discrete_3} and \eqref{weak_form_directors_discrete_4} and subtracting \eqref{weak_form_directors_discrete_1} and \eqref{weak_form_directors_discrete_2} eventually yields for arbitrary $\vec{\xi}$
\begin{equation} \label{eq_discrete_time_angmom_balance}
    \begin{aligned}
         & \vec{l}(\vec{q}\h\npe,\vec{v}\h\npe) - \vec{l}(\vec{q}\h\n,\vec{v}\h\n) =                                                                                                                                                 \\
         & \qquad \qquad  h \int_\domainspace ( \centerline\h\npeh \times \bar{\vec{n}}\npeh + \bar{\vec{m}}\npeh ) \d s  +  h \left[ \centerline\h\npeh \times (\vec{n}_N)\npeh + (\vec{m}_N)\npeh \right]_{\partial \domainspace},
    \end{aligned}
\end{equation}
where we have also substituted \eqref{eq_discrete_difference_angmom}. This calculation further involves the identities \eqref{eq_T_identities}.

\section{System matrices after spatial discretization} \label{sec_appendix_system_matrices}
Here, we display detailed definitions of the matrices contained in the PH system \eqref{compliance_form_discrete} after spatial discretization.
Consistently with \eqref{eq_ansatz}, we write the FE ansatz as
\begin{equation}
    \begin{aligned}
        \begin{bmatrix}
            \centerline\h(s,t) \\
            \vec{d}\h(s,t)
        \end{bmatrix} & = \begin{bmatrix}
                              \vec{\Phi}_\varphi(s) & \zeros          \\
                              \zeros                & \vec{\Phi}_d(s)
                          \end{bmatrix}
        \begin{bmatrix}
            \hat{\centerline}(t) \\
            \hat{\vec{d}}(t)
        \end{bmatrix} ,
                              & \quad
        \begin{bmatrix}
            \vphi\h(s,t) \\
            \vec{v}_d\h(s,t)
        \end{bmatrix}
                              & = \begin{bmatrix}
                                      \vec{\Phi}_\varphi(s) & \zeros          \\
                                      \zeros                & \vec{\Phi}_d(s)
                                  \end{bmatrix}
        \begin{bmatrix}
            \nodevphi(t) \\
            \hat{\vec{v}}_d(t)
        \end{bmatrix} ,
        \\
        \begin{bmatrix}
            \vec{N}\h(s,t) \\
            \vec{M}\h(s,t)
        \end{bmatrix}     & =
        \begin{bmatrix}
            \vec{\Psi}_N(s) & \zeros          \\
            \zeros          & \vec{\Psi}_M(s)
        \end{bmatrix}
        \begin{bmatrix}
            \hat{\vec{N}}(t) \\
            \hat{\vec{M}}(t)
        \end{bmatrix} ,
                              & \quad  \vec{\lambda}\h(s,t)                & = \vec{\Xi}(s) \hat{\vec{\lambda}}(t) ,
    \end{aligned}
\end{equation}
together with the same ansatz for the respective test functions.
Plugging these relations into the weak form \eqref{weak_form_directors} and accounting for arbitrary test functions yields the semi-discrete set of differential-algebraic equations
\begin{equation}
    \begin{aligned}
        \dot{\nodecenterline}                  & = \nodevphi                      ,                                                                                                                                                                                                                    \\
        \dot{\hat{\vec{d}}}                    & = \nodal{v}_d                                                 ,                                                                                                                                                                                       \\
        \rho A \nodal{M}_\varphi  \nodedotvphi & = - \nodal{G}_{\varphi N}(\nodal{q})\transp \nodal{N} + \nodal{M}_\varphi \nodal{\bar{n}} + \nodal{B}_{\partial,n} \vec{n}_N                 ,                                                                                                        \\
        \nodal{M}_d \dot{\nodal{v}}_d          & =- \nodal{G}_{d N}(\nodal{q})\transp\nodal{N} - \nodal{G}_{d M}(\nodal{q})\transp \nodal{M} - \nodal{G}_{d}(\nodal{d})\transp\nodal{\lambda}  + \nodal{B}_{\domainspace,m}(\nodal{q}) \nodal{\bar{m}} + \nodal{B}_{\partial,m}(\nodal{q}) \vec{m}_N , \\
        \nodal{C}_N \dot{\nodal{N}}            & = \nodal{G}_{\varphi N}(\nodal{q}) \nodevphi + \nodal{G}_{d N}(\nodal{q}) \nodal{v}                                                                                                ,                                                                  \\
        \nodal{C}_M \dot{\nodal{M}}            & = \nodal{G}_{d M}(\nodal{q}) \nodal{v}_d                        ,                                                                                                                                                                                     \\
        \vec{0}                                & = \nodal{G}_{d}(\nodal{d}) \nodal{v}_d ,
    \end{aligned}
\end{equation}
which are equivalent to the finite-dimensional PH system \eqref{compliance_form_discrete}. To show this in depth, consider the abbreviations
\begin{align}
     & \nodal{M}_o = \int_\domainspace \vec{\Phi}\transp \vec{M}_o \vec{\Phi} \d s = \diag(\rho A \nodal{M}_\varphi, \nodal{M}_d) , \quad \nodal{C} = \int_\domainspace \vec{\Psi}\transp \vec{C} \vec{\Psi} \d s = \diag(\nodal{C}_N, \nodal{C}_M) ,                                          \\
     & \nodal{G}(\nodal{q}) = \begin{bmatrix}
                                  \vec{0} & \nodal{G}_d
                              \end{bmatrix},                                                                                                                                                        \quad \nodal{J}_{\vec{\sigma}\vec{v}}(\nodal{q}) = \begin{bmatrix}
                                                                                                                                                                                                                                                           \nodal{G}_{\varphi N} & \nodal{G}_{d N} \\
                                                                                                                                                                                                                                                           \zeros                & \nodal{G}_{d M}
                                                                                                                                                                                                                                                       \end{bmatrix}, \quad
    \nodal{B}_{\domainspace}(\nodal{q})  = \begin{bmatrix}
                                               \nodal{M}_\varphi & \zeros                     \\
                                               \zeros            & \nodal{B}_{\domainspace,m}
                                           \end{bmatrix} , \quad \nodal{B}_{\partial}(\nodal{q}) = \begin{bmatrix}
                                                                                                       \nodal{B}_{\partial,n} & \zeros                 \\
                                                                                                       \zeros                 & \nodal{B}_{\partial,m}
                                                                                                   \end{bmatrix} \notag
\end{align}
with the mass and compliance matrices
\begin{equation*}
    \nodal{M}_\varphi = \int_{\domainspace} \vec{\Phi}_\varphi\transp \vec{\Phi}_\varphi \d s , \quad \nodal{M}_d    =  \int_{\domainspace} \vec{\Phi}_d \transp \vec{M}_\rho \vec{\Phi}_d \d s , \quad \nodal{C}_N           = \int_{\domainspace} \vec{\Psi}_N \transp \vec{C}_N \vec{\Psi}_N \d s , \quad
    \nodal{C}_M                        =  \int_{\domainspace} \vec{\Psi}_M \transp \vec{C}_M \vec{\Psi}_M\d s ,
\end{equation*}
and discrete elements of the structure matrix
\begin{equation}
    \begin{aligned}
        \nodal{G}_{\varphi N}(\nodal{q}) & = \int_{\domainspace} \vec{\Psi}_N \transp \rotmat(\vec{d}\h)\transp \partial_s \vec{\Phi}_\varphi \d s ,                             \\
        \nodal{G}_{d N}(\nodal{q})       & = \int_{\domainspace} \vec{\Psi}_N \transp \vec{J}_N(\centerline\h_{,s}) \vec{\Phi}_d \d s ,                                          \\
        \nodal{G}_{d M} (\nodal{q})      & = \int_{\domainspace} \vec{\Psi}_M \transp (\vec{L}(\vec{d}\h_{,s})\vec{\Phi}_d - \vec{L}(\vec{d}\h) \partial_s \vec{\Phi}_d ) \d s ,
    \end{aligned}
\end{equation}
Concerning the treatment of the constraint equation, consider the approximation of \eqref{weak_form_directors_7} using $\nnode$ delta functions $\bar{\delta}(s-s_i)$, cf.~\cite[Chap.~1.11]{hartmann_1997_signals}, given by
\begin{equation} \label{sifting}
    \inner{\delta \vec{\lambda}\h}{\vec{G}_d(\vec{d}\h) \vec{v}_d\h
    }                       = \inner{\sum_{i=1}^{\nnode} \bar{\delta}(s-s_i)\delta\vec{\lambda}_i }{\vec{G}_d(\vec{d}\h) \vec{v}_d\h}       =  \sum_{i=1}^{\nnode} \delta \vec{\lambda}_i\transp        \vec{G}_d(\nodal{d}_i) \nodal{v}_{d,i}                          =0 .
\end{equation}
Therein, $\delta \vec{\lambda}_i \in \mathbb{R}^6$ are the test functions corresponding to the Lagrange multipliers that enforce the orthonormality of the directors located at the FE node at $s=s_i$, comprised in $\nodal{d}_i$. The related director velocities are comprised in $\vec{v}_{d,i} \in \mathbb{R}^9$. The last identity in \eqref{sifting} has made use of the sifting property of dirac functions \cite[Eq.~1.159]{hartmann_1997_signals}, i.e., $\inner{f(s)}{\bar{\delta}(s-a)} = f(a)$ for some $f: \domainspace \to \mathbb{R}$ and $a \in \domainspace$.
Eventually, we may rewrite \eqref{sifting} in terms of the system vectors $\delta\nodal{\lambda}$ and $\nodal{v}_d$, such that
\begin{equation} \label{sifting_2}
    \delta \nodal{\lambda}\transp \nodal{G}_{d}(\nodal{d}) \nodal{v}_d = 0, \qquad \text{where} \quad \nodal{G}_{d}(\nodal{d}) = \begin{bmatrix}
        \vec{G}_d(\nodal{d}_1) & \zeros                 &        & \zeros                        \\
        \zeros                 & \vec{G}_d(\nodal{d}_2) &        & \zeros                        \\
                               &                        & \ddots &                               \\
        \zeros                 & \zeros                 &        & \vec{G}_d(\nodal{d}_{\nnode})
    \end{bmatrix}
\end{equation}
has therefore full rank.
Futhermore, the port matrices are defined as
\begin{align}
    \nodal{B}_{\domainspace,m}(\nodal{q})  =         \int_{\domainspace} \vec{\Phi}_d \transp \vec{T}(\vec{d}\h)         \vec{\Phi}_\varphi \d s ,
    \quad
    \nodal{B}_{\partial,n}      =  \vec{\Phi}_\varphi\transp|_{\partial\domainspace_N}           , \quad
    \nodal{B}_{\partial,m}(\nodal{q})       = (\vec{\Phi}_d\transp \vec{T}(\vec{d}\h))|_{\partial\domainspace_N}     .
\end{align}
Consequently, the output relations in \eqref{discrete_ph} can be written in a more detailed manner as
\begin{align}
    \nodal{y}_\domainspace & = \nodal{B}_{\domainspace}\transp \nodal{v} =  \begin{bmatrix}
                                                                                \nodal{M}_\varphi & \zeros                            \\
                                                                                \zeros            & \nodal{B}_{\domainspace,m}\transp
                                                                            \end{bmatrix} \begin{bmatrix}
                                                                                              \nodevphi(t) \\
                                                                                              \hat{\vec{v}}_d(t)
                                                                                          \end{bmatrix}
    = \begin{bmatrix}
          \int_\domainspace \vec{\Phi}_\varphi\transp \vphi\h \d s \\
          \int_\domainspace \vec{\Phi}_d\transp \vec{\omega}\h \d s
      \end{bmatrix}
\end{align}
and
\begin{align}
    \nodal{y}_\partial     = \nodal{B}_{\partial}\transp \nodal{v} = \begin{bmatrix}
                                                                         \nodal{B}_{\partial,n}\transp & \zeros                        \\
                                                                         \zeros                        & \nodal{B}_{\partial,m}\transp
                                                                     \end{bmatrix} \begin{bmatrix}
                                                                                       \nodevphi(t) \\
                                                                                       \hat{\vec{v}}_d(t)
                                                                                   \end{bmatrix}
    = \begin{bmatrix}
          \vec{\Phi}_\varphi |_{\partial\domainspace_N}\nodevphi(t) \\
          (\vec{T}(\vec{d}\h)\transp \vec{\Phi}_d)|_{\partial\domainspace_N} \hat{\vec{v}}_d(t)
      \end{bmatrix}
    = \begin{bmatrix}
          \vphi\h \\
          \vec{\omega}\h
      \end{bmatrix}_{\partial\domainspace_N} .
\end{align}
While the discrete output can be regarded as the projection of the centerline velocity and angular velocity onto the FE nodes, the boundary output can be traced back to the simple evaluation of the approximated centerline velocity and angular velocity at the Neumann boundaries. Both results ensure that the discrete power balance equation \eqref{energy_balance_discrete_space} yields
\begin{equation}
    \dot{\nodal{H}} = \begin{bmatrix}
        \nodal{y}_\domainspace \\
        \nodal{y}_\partial
    \end{bmatrix}\transp
    \begin{bmatrix}
        \nodal{u}_\Omega \\
        \vec{u}_\partial
    \end{bmatrix} = \int_\domainspace (\bar{\vec{n}}\transp\vphi\h + \bar{\vec{m}}\transp  \vec{\omega}\h) \d s + \vec{u}_\partial\transp\nodal{y}_\partial ,
\end{equation}
which is the consistent approximation of the infinite-dimensional power balance \eqref{power_balance_cont}, see also \eqref{boundary_terms}.

\end{document}